\documentclass[preprint,sort&compress,12pt]{elsarticle}
\usepackage{bbm}
\usepackage{amstext}
\usepackage{amsmath}
\usepackage{amssymb}
\usepackage{mathrsfs}
\usepackage{stmaryrd}
\usepackage{bm}
\usepackage{pstricks}
\usepackage{pst-coil}
\usepackage{pst-3d}
\usepackage{amsfonts} 
\usepackage{amsthm}
\usepackage{xcolor}
\newcommand{\mm}{\mathrm}

\newcommand{\be}{\begin{equation}}
\newcommand{\bea}{\begin{equation}\begin{aligned}}
\newcommand{\beas}{\begin{equation*}\begin{aligned}}
\newcommand{\eeas}{\end{aligned}\end{equation*}}
\newcommand{\eea}{\end{aligned}\end{equation}}
\newcommand{\ee}{\end{equation}}

\renewcommand{\div}{{\rm div }}

\begin{document}
\begin{frontmatter}
\title{On Inhibition of Rayleigh--Taylor Instability by  Horizontal \\
Magnetic Field in an Inviscid MHD Fluid with Velocity Damping}

\author[sJ]{Fei Jiang}
\ead{jiangfei0951@163.com}
\author[FJ]{Song Jiang}
 \ead{jiang@iapcm.ac.cn}
\author[sJ]{Youyi Zhao\corref{cor1}}
\ead{zhaoyouyi957@163.com}
\cortext[cor1]{Corresponding author. }
\address[sJ]{College of Mathematics and
Computer Science, Fuzhou University, Fuzhou, 350108, China.}
\address[FJ]{Institute of Applied Physics and Computational Mathematics,
 Beijing, 100088, China.}

\begin{abstract}
It is still an open problem whether the inhibition phenomenon of Rayleigh--Taylor (RT) instability by horizontal magnetic field can be mathematically proved in a  non-resistive magnetohydrodynamic (MHD) fluid in a two-dimensional (2D) horizontal slab domain, since it had been roughly verified by a 2D linearized motion equations in 2012 \cite{WYC}. In this paper, we find that this inhibition phenomenon can be rigorously verified in the inhomogeneous, incompressible, inviscid case with velocity damping. More precisely, there exists a critical number $m_{\mm{C}}$ such that if the strength $|m|$ of horizontal magnetic field is bigger than $m_{\mm{C}}$, then the small perturbation solution around the magnetic RT equilibrium state is exponentially stable in time.
Our result is also the first mathematical one based on the nonlinear motion equations for the proof of inhibition of flow instabilities by a horizontal magnetic field in a horizontal slab domain.  In addition, we also provide a nonlinear instability result for the case $|m|\in [0,m_{\mm{C}})$. Our instability result presents that horizontal magnetic field can not inhibit the RT instability, if it's strength is to small.
\end{abstract}
\end{frontmatter}
\newtheorem{thm}{Theorem}[section]
\newtheorem{lem}{Lemma}[section]
\newtheorem{pro}{Proposition}[section]
\newtheorem{concl}{Conclusion}[section]
\newtheorem{cor}{Corollary}[section]
\newproof{pf}{Proof}
\newdefinition{rem}{Remark}[section]
\newtheorem{definition}{Definition}[section]

\section{Introduction}\label{introud}
\numberwithin{equation}{section}
The equilibrium of the heavier fluid on top of the lighter one
and both subject to the gravity is unstable. In this case, the equilibrium state is unstable to sustain small disturbances, and this unstable disturbance will grow and lead to the release of potential energy, as the heavier fluid moves down under the gravity force, and the lighter one is displaced upwards. This phenomenon was first studied by Rayleigh \cite{RLIS} and then Taylor \cite{TGTP}, is called therefore the RT instability.
In the last decades, this phenomenon has been extensively investigated from mathematical,
physical and numerical aspects, see \cite{CSHHSCPO,WJH} for examples. It has been also widely investigated how other physical factors, such as elasticity \cite{JFJWGCOSdd,FJWGCZXOE}, rotation \cite{CSHHSCPO,BKASMMHRJA},
internal surface tension \cite{GYTI2,WYJTIKC2014ARMA,JJTIWYJ2016CMP}, magnetic field \cite{JFJSWWWOA,JFJSJMFM,JFJSJMFMOSERT,JFJSWWWN,WYC,WYJ2019ARMA}
and so on, influence the dynamics of RT instability.

In this paper, we are interested in the phenomenon of inhibition of RT instability by magnetic field. This topic goes back to the theoretical work of Kruskal and Schwarzchild \cite{KMSMSP}. They analyzed the effect of an impressed horizontal magnetic field on the growth of RT instability in the horizontally periodic motion  of a completely ionized plasma (with zero resistance) in three dimensions in 1954, and pointed out that the curvature of the magnetic lines can influence the development of instability, but can not inhibit the growth of RT instability. The inhibition of RT instability by an impressed vertical magnetic field was first verified for inhomogeneous, incompressible, non-resistive magnetohydrodynamic (MHD) fluids in three dimensions by Hide \cite{CSHHSCPO,HRWP}. In 2012, Wang noticed that an impressed horizontal magnetic field can inhibit RT instability in a non-resistive MHD fluid  in two dimensions \cite{WYC}. Later Jiang--Jiang further found that magnetic field always inhibit RT instability, if the condition is satisfied that the non-slip velocity boundary-value condition is imposed in the direction of impressed magnetic field. In this paper, we call such condition the ``fixed condition" for the sake of simplicity.  It should be noted that all the  results stated previously are based on the linearized non-resistive MHD equations.
For reader's convenience, we also summarize the known linear results as follows:
\begin{center}
Whether an impressed horizontal/vertical magnetic field can inhibit RT \\ instability in a slab domain with a non-slip boundary-value condition of velocity?\\[0.4em]
\begin{tabular}{|c|c|c|}
\cline{2-2}  \hline    & \vspace{0.1em}{horizontal} & vertical  \\
   \hline
  2D & Yes & Yes \\
   \hline
  3D & No & Yes \\
  \hline
\end{tabular}
\end{center}

Recently, Jiang--Jiang further established a so-called magnetic inhibition theory in viscous non-resistive MHD fluids, which reveals the physical effect of fixed condition in magnetic inhibition phenomena \cite{JFJSARMA2019}. Roughly specking, let us consider an element line along an impressed field in the rest state of a MHD fluid, then the element line can be regarded as an elastic string. Thus, the bent element line will restore to its initial location under the magnetic tension, the fixed condition, as well as viscosity.  By the magnetic inhibition mechanism of a non-resistive MHD fluid, the assertions in the table above seem to be obvious. However, rigorous mathematical proofs  are not easy.

Thanks to the multi-layers method developed in the proof of well-posed problem of surface waves \cite{GYTIAED}, recently the inhibition phenomenon of RT instability by a magnetic field had been rigorously proved based on the (nonlinear) non-resistive MHD equations under a fixed condition, for example, Wang verified the inhibition phenomenon by an impressed non-horizontal magnetic field in the stratified incompressible viscous MHD fluid in a 3D slab domain \cite{WYJ2019ARMA}; moreover, he also proved that the impressed horizontal magnetic field can not inhibit the RT instability for the horizontally periodic motion. Similar results can be also found in other magnetic inhibition phenomena, see \cite{JFJSSETEFP} for Parker instability and \cite{JFJSOUI} for thermal instability.

However at present it is still an open problem whether the  phenomenon of inhibition of RT instability by a horizontal magnetic field can be rigorously proved in a non-resistive MHD fluid in a 2D slab domain. To our knowledge, there are also not any available mathematical proof for the inhibition of other flow instabilities by a horizontal magnetic field in a horizontal slab domain in the both 2D and 3D cases. The purpose of this paper is to move a step in this direction. Fortunately we find that this inhibition phenomenon can be mathematically verified in the inhomogeneous, incompressible, inviscid, non-resistive MHD fluid with velocity damping in two-dimensions. More precisely, there exists a critical number $m_{\mm{C}}$ such that if the strength $|m|$ of a horizontal magnetic field is bigger than $m_{\mm{C}}$, then the small perturbation solution around the magnetic RT equilibrium state is exponentially stable in time, i.e., RT instability can be inhibited by a horizontal magnetic field in a 2D slab domain.
Next we mathematically formulate our result.

\subsection{Mathematical formulation for the magnetic RT problem}\label{subsec:01}

The system of motion equations of an inhomogeneous, incompressible, inviscid, non-resistive MHD fluid with velocity damping in the presence
of a gravitational field in a two-dimensional domain $\Omega$ reads as follows:
\begin{equation}
\label{0101}
\begin{cases}
\rho_t+ v\cdot \nabla \rho=0,\\[1mm]
\rho v_{t}+ \rho v\cdot\nabla v
+\nabla \left(P+\lambda|M|^2/2\right)+a\rho v =\lambda M\cdot\nabla M-\rho g \mathbf{e}_2, \\
 {M}_{t}+ v\cdot\nabla {M}=M\cdot\nabla v, \\
\div  v =\mm{div}M=0.
\end{cases}
\end{equation}
Next we shall explain the mathematical notations in system \eqref{0101}.

The unknowns $\rho:=\rho(x,t)$, $v:={v}(x,t)$,
$M:= {M}(x,t)$ and $P:=P(x,t)$ denote the density, velocity,
magnetic field and kinetic pressure of incompressible MHD fluids, resp.. $x\in \Omega\subset\mathbb{R}^2$ and $t>0$ are spacial variables and time variables, resp.. The constants $\lambda>0$, $g>0$ and $a\geqslant 0$ stand for the permeability of vacuum, the gravitational constant and the velocity damping coefficient, resp.. $\mathbf{e}_2=(0,1)^{\mm{T}}$ represents the vertical unit vector, and  $-\rho g \mathbf{e}_2$ the gravitational force, where the
 superscript ${\mm{T}}$ denotes the transposition.

Since we consider the horizontally periodic motion solution of \eqref{0101}, we define the  horizontally periodic domain
\begin{align}\label{0101a}
\Omega:= 2\pi\mathbb{T} \times(0,h),
\end{align}
where $\mathbb{T}=\mathbb{R}/\mathbb{Z} $. For the horizontally periodic domain $\Omega$, the 1D periodic domain  $2\pi\mathbb{T}\times \{0,h\}$, denoted by $\partial\Omega$, customarily is regarded as the boundary of $\Omega$. For the well-posedness of the system \eqref{0101}, we shall pose the following initial-boundary value conditions:
\begin{align}
& (\rho,v,M)|_{t=0}=(\rho^0,v^0,M^0),\label{20210031303} \\
& \label{0101b}
v|_{\partial\Omega}\cdot\vec{\mathbf{n}}=0 ,
\end{align}
where $\vec{\mathbf{n}}$ denotes the outward unit normal vector on  $\partial\Omega$. Here and in what follows, we always use the superscript $0$ to emphasize the initial data.

Now we choose a RT density profile $\bar{\rho}:=\bar{\rho}(x_2)$, which is independent of $x_1$ and satisfies
\begin{align}
&\label{0102}
\bar{\rho}\in {C^4}(\overline{\Omega}),\ \inf_{ x\in {\Omega}}\bar{\rho}>0,\\[1mm]
&\label{0102n}\bar{\rho}'|_{x_2=y_{2}}>0
\mbox{ for some } y_{2}\in \{x_2~|~(x_1,x_2)^{\mm{T}}\in \Omega\},
\end{align}
where $\bar{\rho}':=\mm{d}\bar{\rho}/\mm{d}x_2$ and $\overline{\Omega}:=\mathbb{R}\times [0,h]$.
We remark that the second condition in \eqref{0102}
prevents us from treating vacuum, while  the condition in \eqref{0102n} is called RT condition,
which assures that there is at least a region in which the density is larger with increasing height $x_2$ and leads to the classical RT instability, see \cite[Theorem 1.2]{JFJSO2014}.

With RT density profile in hand, we further define a magnetic RT equilibria $r_M:=(\bar{\rho}, 0, \bar{M})$, where $\bar{M} =(m,0)^{\mm{T}}$ and $m$ is a constant. We often call $\bar{M}$ an impressive horizontal magnetic field,  while the pressure profile $\bar{P}$ under the equilibrium state is determined by the relation
\begin{equation}
\nabla \bar{P}=-\bar{\rho}g  {e}_2\mbox{ in } {\Omega}.
\label{equcomre}
\end{equation}

Denoting the perturbation around the magnetic RT equilibria by
$$\varrho=\rho -\bar{\rho},\ v= v- {0},\ N=M-\bar{M} ,$$
 and then using the relation \eqref{equcomre}, we obtain the system of perturbation equations from \eqref{0101}:
\begin{equation}\label{0103} \begin{cases}
\varrho_t+{  v}\cdot\nabla (\varrho+\bar{\rho})=0, \\[1mm]
(\varrho+\bar{\rho}){  v}_t+(\varrho+\bar{\rho}){  v}\cdot\nabla
{ v}+\nabla  \beta+a \rho v \\
=  \lambda (N+\bar{M})\cdot \nabla  N - \varrho  g \mathbf{e}_2,\\[1mm]
N_t+ v\cdot\nabla  N =(N+\bar{M})\cdot \nabla v,\\[1mm]
 \mm{div}v= \mathrm{div} N=0,\end{cases}\end{equation}
where $\beta:= P-\bar{P}+\lambda  (| M |^2-|\bar{M}|^2)/2$, and we call $ \beta $ the total perturbation pressure.
The corresponding initial-boundary value conditions read as follows:
\begin{align} \label{c0104}
&(\varrho,v, {N} )|_{t=0}=(\varrho^0,v^0,N^0) ,  \\
& \label{0105}
v|_{\partial\Omega}\cdot\vec{\mathbf{n}}=0 .
\end{align}
We call the initial-boundary value problem \eqref{0103}--\eqref{0105} magnetic RT problem for the sake of simplicity. Obviously, to mathematically prove the inhibition of RT instability by a horizontal magnetic field in a 2D slab domain, it suffices to verity the stability in time of magnetic RT problem with \emph{some non-trivial initial datum}.

We mention that the well-posedness problem of inviscid fluids with velocity damping had been widely investigated, see \cite{TZWY2013JDE,PRHZK2009JDE,WWKYT2001JDE,STCTBWDH2012CPDE,HFMPRHARMA2003,ZHHARMA2017,ZHHARMA2021,LTZHHCPAM2016} for examples.
Recently, some authors had further studied the well-posedness problem of the motion equations of incompressible inviscid,  non-resistive MHD fluids, i.e., taking $\rho =1$ and $g=0$ in \eqref{0101}. For examples, Wu--Wu--Xu first given the existence of unique global(-in-time) solutions with algebraic decay-in-time for the 2D Cauchy problem with small initial perturbation  \cite{WJHWYFXXJG}, and Du--Yang--Zhou obtained the existence of unique global solutions with exponential decay-in-time for the initial-boundary value problem in a 2D slab domain  with small initial perturbation around some non-trivial equilibria \cite{DYYWZYOSJMA}. It should be noted that the mathematical methods adopted in \cite{WJHWYFXXJG,DYYWZYOSJMA} for the well-posedness problem not applied to our stability problem, therefore  next we shall reformulate magnetic RT problem  in Lagrangian coordinates as in \cite{WYJ2019ARMA,JFJSSETEFP,JFJSOUI}.

\subsection{Reformulation in Lagrangian coordinates}\label{subsec:02}

Let the flow map $\zeta$  be the solution to the initial-value problem
\begin{equation}
\label{201806122101}
            \begin{cases}
\partial_t \zeta(y,t)=v(\zeta(y,t),t)&\mbox{in }\Omega,
\\[1mm]
\zeta(y,0)=\zeta^0(y)&\mbox{in }\Omega,
                  \end{cases}
\end{equation}
where the
invertible mapping $\zeta^0:=\zeta^0(y)$ maps $\Omega$ to $\Omega$, and satisfies
\begin{align}
&\label{zeta0inta}
J^0:=\det \nabla \zeta^0=1 \mbox{ in } {\Omega},\\[1mm]
&\label{zeta0inta0}
\zeta^0\cdot\vec{\mathbf{n}} =y\cdot\vec{\mathbf{n}}\mbox{ on } \partial\Omega.
\end{align}
Here and in what follows, ``$\det$" denotes a determinant of matrix.

We denote the Eulerian coordinates by $(x,t)$ with
$x =\zeta(y,t)$ and the Lagrangian coordinates by $(y, t)\in\Omega\times \mathbb{R}_0^+$, where $\mathbb{R}_0^+:=[0,\infty)$.
We further assume that, for each fixed $t>0$,
\begin{align}
& \zeta|_{y_2=i}   : \mathbb{R}\to \mathbb{R}\mbox{ is a } C^2(\mathbb{R})\mbox{-diffeomorphism mapping for }i=0,\ h,\label{20210301715x}\\
&\zeta   :  \overline{\Omega}\to \overline{\Omega} \mbox{ is a } C^2(\overline{\Omega})\mbox{-diffeomorphism mapping}.\label{20210301715}
\end{align}

Since $v$ satisfies the divergence-free condition,
and non-slip boundary-value condition \eqref{0105}, we can deduce from \eqref{201806122101}--\eqref{zeta0inta0} that
\begin{align}
&\nonumber
J:=\det\nabla \zeta=1\mbox{ in } \Omega,\\
&\nonumber
\zeta\cdot\vec{\mathbf{n}}=y\cdot\vec{\mathbf{n}} \mbox{ on } \partial\Omega.
\end{align}

We define the matrix $\mathcal{A}:=(\mathcal{A}_{ij})_{2\times 2}$  via
\begin{align}\nonumber
\mathcal{A}^{\mm{T}}=(\nabla\zeta)^{-1}:=
(\partial_j \zeta_i)^{-1}_{2\times 2}.
\end{align}
Then we further define
the differential operators $\nabla_{\mathcal{A}}$, $\mm{div}_{\mathcal{A}}$ and $\mm{curl}_{\mathcal{A}}$ as follows: for a scalar function $f$ and a  vector function $X:=(X_1,X_2)^{\mm{T}}$,
\begin{align}
&\nabla_{\mathcal{A}}f:=(\mathcal{A}_{1k}\partial_kf,
\mathcal{A}_{2k}\partial_kf)^{\mm{T}},\ \mm{div}_{\mathcal{A}}(X_1,X_2)^{\mm{T}}:=\mathcal{A}_{lk}\partial_k X_l, \ \mm{curl}_{\mathcal{A}}X:=\mathcal{A}_{1k}\partial_{k}X_2-\mathcal{A}_{2k}\partial_{k}X_1,
 \nonumber \end{align}
where we have used the Einstein convention of summation over repeated indices, and $\partial_k:=\partial_{y_k}$. In particular, $\mm{curl}X:=\mm{curl}_IX$, where $I$ represents an identity matrix. \emph{In addition, we will denote $(\mm{curl}_{\mathcal{A}} X^1,\ldots, \mm{curl}_{\mathcal{A}}X^n)^{\mm{T}}$ by $\mm{curl}_{\mathcal{A}}(X^1,\ldots X^n)^{\mm{T}}$ for simplicity, where $X^i=(X^i_1,X^i_2)^{\mm{T}}$ is a vector function for $1\leqslant i \leqslant  n$}.

Defining the Lagrangian unknowns
\begin{equation*}
(  \vartheta, u ,Q, B)(y,t)=(\rho,v,P+\lambda|M|^2/2,M)(\zeta(y,t),t) \mbox{ for } (y,t)\in \Omega \times\mathbb{R}^+_0,
\end{equation*}
then in Lagrangian coordinates, the initial-boundary value problem of \eqref{0101}, \eqref{20210031303} and \eqref{0101b} is rewritten as follows:
\begin{equation}\label{01dsaf16asdfasf00}
\begin{cases}
\zeta_t=u ,  \ \vartheta_t=0 ,\
\div_\mathcal{A}u=0 ,    \\[1mm]
\vartheta u_t+\nabla_{\mathcal{A}}Q+a\vartheta u=\lambda B\cdot\nabla_{\mathcal{A}}B-\vartheta g \mathbf{e}_2 , \\[1mm]
B_t=B\cdot \nabla_{\mathcal{A}}u   ,  \
\div_\mathcal{A}B=0 ,   \\[1mm]
(\zeta,\vartheta,u, B)|_{t=0}=(\zeta^0,\vartheta^0,u^0, B^0) , \\[1mm]
(\zeta-y,u)|_{ \partial\Omega}\cdot\vec{\mathbf{n}}=0  ,
\end{cases}
\end{equation}
where $(\vartheta^0,  v^0,B^0 )=
 (\rho^0(\zeta^0),v^0(\zeta^0), M^0(\zeta^0)) $. In addition, the relation \eqref{equcomre} in Lagrangian coordinates reads as follows:
 \begin{equation}
\label{dstist01}
\nabla_{\mathcal{A}}\bar{P}(\zeta_2)
=-\bar{\rho}(\zeta_2)g\mathbf{e}_2 .
\end{equation}

Let $\eta=\zeta-y$, $\eta^0=\zeta^0-y$, $q=Q-\bar{P}(\zeta_2)-\lambda|\bar{M}|^2/2$, $\mathcal{A}=(\nabla \eta+I)^{\mm{T}}$ and
\begin{align}
\label{202009130836}
{G}_{\eta}:=\bar{\rho}(\eta_2(y,t)+y_2)-\bar{\rho}(y_2).
\end{align}
If $\zeta^0$, $\vartheta^0$ and $B^0$ satisfy
 \begin{align}
\nonumber
B^0= m\partial_1\zeta^0 \mbox{ and }
\vartheta^0  =\bar{\rho}(y_2),
\end{align}
then the initial-boundary value problem \eqref{01dsaf16asdfasf00}, together with the relation \eqref{dstist01}, implies that
  \begin{equation}\label{01dsaf16asdfasf}
                              \begin{cases}
\eta_t=u ,\\[1mm]
\bar{\rho}u_t+\nabla_{\mathcal{A}} q+a\bar{\rho} u=\lambda m^2\partial_1^2\eta+g{G}_{\eta}\mathbf{e}_2  ,\\[1mm]
\div_{\mathcal{A}} u=0  , \\[1mm]
(\eta,u)|_{t=0}=(\eta^0,u^0)  , \\[1mm]
(\eta,u)|_{\partial\Omega}\cdot\vec{\mathbf{n}}=0
\end{cases}
\end{equation}
and
\begin{align}\label{202012280945}
\vartheta=\bar{\rho}(y_2),\
B=m\partial_1\zeta,
\end{align}
please refer to \cite{JFJSJMFMOSERT} for the derivation.

It should be noted that \eqref{01dsaf16asdfasf}, together with \eqref{202012280945}, also implies \eqref{01dsaf16asdfasf00}. In addition, noting that $q$ is the sum of the perturbation pressure and perturbation magnetic pressure in Lagrangian coordinates, however we still call $q$ the perturbation pressure for the sake of simplicity. From now on, we call the  initial-boundary value problem \eqref{01dsaf16asdfasf} the transformed MRT problem. Obviously the stability problem of magnetic RT problem reduces to investigate the stability of the transformed MRT problem.

\subsection{Notations}\label{subsec:04}

Before stating our main results on  the transformed MRT problem, we shall  introduce simplified notations throughout this paper.

\begin{enumerate}[(1)]
  \item Simplified basic notations:  $I_a:=(0,a)$ denotes a time interval, in particular, $I_\infty=\mathbb{R}^+$.  $\overline{S}$ denotes the closure of the set $S\subset \mathbb{R}^n$ with $n\geqslant 1$, in particular, $\overline{I_T}:=[0,T]$.  $\int:=\int_\Omega=\int_{(0,2\pi)\times (0,h)}$. $(u)_{\Omega}$ denotes the mean value of $u$ in a periodic box.   $a\lesssim b$ means that $a\leqslant cb$.  If not stated explicitly,  the positive constant $c$ may depend on $g$, $a$,  $\lambda$, $m$, $\bar{\rho}$ and $\Omega$ in the transformed MRT problem, and may vary from one place to other place. Sometimes we use $c_i$ for $i\geqslant 1$ to replace $c$, and to emphasize that $c_i$ is fixed value. $\alpha$ always denotes the multiindex with respect to the variable $y$,  $|\alpha|=\alpha_1+\alpha_2$ is called the order of multiindex,  $\partial^{\alpha}:=\partial_{1}^{\alpha_1} \partial_{2}^{\alpha_2}$,
$[\partial^{\alpha},\phi]\varphi:=\partial^{\alpha}(\phi\varphi)-\phi\partial^{\alpha}\varphi$
 and \begin{align}
[\partial^{\alpha}\mm{curl}_{\partial_1^j\mathcal{A}},\phi]\chi:=
\partial^{\alpha}\mm{curl}_{\partial_1^j\mathcal{A}}(\phi\chi)-
\phi\partial^{\alpha}\mm{curl}_{\partial_1^j\mathcal{A}}\chi
\mbox{ for }j=0,\ 1.\nonumber
\end{align}  ``$\nabla^{i}f\in \mathbb{X}$" represents that $\partial^\alpha f\in \mathbb{X}$ for any multiindex $\alpha$ satisfying $|\alpha|=i$, where $\mathbb{X}$ denotes some set of functions.
  \item  Simplified Banach spaces, norms and semi-norms:
  \begin{align}
&L^p:=L^p (\Omega)=W^{0,p}(\Omega),\
{H}^i:=W^{i,2}(\Omega ),\ \underline{H}^{i}:=\{w\in H^{i}~|~(w)_{\Omega}=0\}, \nonumber \\[1mm]
&
  {H}^{1,i}:=\{w\in H^i~|~\partial_1w\in H^i\},\ H^{j}_{\mathrm{s}}:=\{w\in {H}^{j}~|~w|_{\partial\Omega}\cdot\vec{\mathbf{n}}=0\}, \nonumber
  \\
& H^{j}_\sigma:=\{w\in {H}^{j}_{\mathrm{s}}~|~\div w=0\}, \ H^{1,j}_{\mathrm{s}}:= H^{j}_{\mathrm{s}}\cap H^{1,j},\ L^p_TH^i:=L^p(I_T,H^i),\nonumber \\
&\textstyle \|\cdot \|_i :=\|\cdot \|_{H^i},\
\|\cdot\|_{k,i}:= \|\partial_{1}^{k}\cdot\|_{i},
 \ \|\cdot\|_{\underline{k},i}:=\sqrt{\sum_{0\leqslant   l \leqslant k}\|\cdot\|_{l,i}^2}, \nonumber
\end{align}
where $1\leqslant p\leqslant \infty$,
and  $i$, $k\geqslant 0$,  $j \geqslant 1$.

In addition, for simplicity, we denote $\sqrt{\sum_{1\leqslant k\leqslant j}\|f^k\|_{\mathcal{X}}^2}$ by $\|(f^1,\ldots,f^j)\|_{\mathcal{X}}$, where $\|\cdot\|_{\mathcal{X}}$ represents a norm or a semi-norm, and $f^k$ may be a  scalar function, a vector or a matrix for $1\leqslant k\leqslant j$.
 \item  simplified function classes:  for integer  $j \geqslant  1 $,
\begin{align}
&{H}^{j}_1:=\{w\in H^{j}~|~\det(\nabla w+I)=1\},
\  H^{1,j}_{1,\mathrm{s}}:= H^{1,j}_{\mathrm{s}}\cap  H_1^{j}, \nonumber \\ & {{C}}^0_{B,\mm{weak}}(\overline{I_T} ,L^2):= L^\infty_TL^2\cap  C^0(\overline{I_T}, L^2_{\mm{weak}}),   \nonumber\\
& H^j_{*}:=\{\xi\in H^j~|~ \xi(y)+y \mbox{ satisfies the diffeomorphism}\nonumber \\
&\qquad \qquad  \mbox{ properties  as } \zeta\mbox{ in \eqref{20210301715x} and \eqref{20210301715}}\}, \nonumber\\
&\mathfrak{C}^0( \overline{I_T},{H}^{1,j}_{ \mathrm{s}} ):=\{\eta\in C(\overline{I_T} ,{H}^{j}_{ \mathrm{s}}) ~|~\nabla^j\partial_1\eta\in   {C}^0_{B  ,\mm{weak}}(\overline{I_T} ,L^2)  \},\nonumber \\
&   {\mathfrak{H}}^{1,4}_{1,*,T}:=\{\eta\in \mathfrak{C}^0(\overline{I_T}, H^{1,4}_{\mm{s}} )~|~ \eta(t) \in  H^4_{*}\mbox{ for each }t\in \overline{I_T}\},\nonumber\\
& {\mathfrak{U}}_T^4:=  \{u\in C^0(\overline{I_T}, H^3_{\mm{s}}) ~|~
\nabla^4 u\in {C}^0_{ B ,\mm{weak}}( \overline{I_T} ,L^2) , \nonumber   \\
&\qquad \qquad\qquad \qquad\qquad \quad
u_t\in C^0(\overline{I_T}, H^2_{\mathrm{s}}) ,\  u_t\in L^\infty_TH^3  \}
\nonumber,\\
&  \mathfrak{Q}_T^4:=\{q\in C(\overline{I_T} ,\underline{H}^3)\cap L^\infty_TH^4~|~q_t \in L^\infty_TH^3 \}.\nonumber
\end{align}
\item Energy integral: for any given $w\in H^1$,
\begin{align}\nonumber
E(w):=\int g\bar{\rho}'w_2^2\mm{d}y-\lambda \| m\partial_1w\|_0^2 .
\end{align}
  \item Energy and dissipation functionals:
\begin{align}
& \mathcal{E}:= \|(\eta,\partial_1\eta,u)\|_4^2+\|(u_t,\nabla q)\|_3^2,\quad \mathcal{E}_{\mm{p}}:= \|(\eta, \partial_1 \eta,u)\|_{1,3}^2{+\|(u,u_t,\nabla q)\|_{3}^2}\nonumber,\\[1mm]
& \mathcal{D}:= \|( u,\partial_1\eta)\|_4^2+\|(u_t,\nabla q)\|_3^2,\quad  \mathcal{D}_{\mm{p}}:= \|(u, \partial_1\eta) \|_{1,3}^2{+\|(u,u_t,\nabla q)\|_{3}^2}.\nonumber
\end{align}
We call $\mathcal{E}$, resp. $\mathcal{D}$ the total energy, resp. dissipation functionals, and $\mathcal{E}_{\mm{p}}$, resp. $\mathcal{D}_{\mm{p}}$ the partial energy, resp. dissipation  functionals.  In addition, we use the following notation for simplicity
\begin{align}I^0:=\|(\eta^0,\partial_1 \eta^0,u^0)\|_4^2. \label{202104121056}
\end{align}
\end{enumerate}
\subsection{Main results}\label{subsec:04xx}
Now we introduce the stability result for the transformed MRT problem.
\begin{thm}[Stability]\label{thm2}
 Assume  $\bar{\rho}$ satisfies \eqref{0102}--\eqref{0102n} and \begin{align}
\label{2020102241504}
|m|>m_{\mm{C}}:=\sqrt{\sup_{ { w}\in H_{\sigma}^1}\frac{g\int\bar{\rho}' { w}_2^2\mm{d}y}
{\lambda\|\partial_1  w \|^2_0}}.
 \end{align}
We further assume $a>0$, $(\eta^0,u^0)\in (H^{1,4}_{1,\mathrm{s}}\cap H^4_*) \times  H^4_{\mathrm{s}}$
and $\mm{div}_{\mathcal{A}^0}u^0=0$,
where $\mathcal{A}^0:=(\nabla \eta^0+I)^{-\mm{T}}$.
Then there exist  a sufficiently small constant $\delta >0$  such that, for any $(\eta^0,u^0)$  satisfying
$I^0\leqslant\delta^2$,
the transformed MRT problem \eqref{01dsaf16asdfasf} admits a unique global classical solution $(\eta,u,q)$ in the function class $\mathfrak{H}^{1,4}_{1,*,\infty}\times \mathfrak{U}_\infty^4\times   \mathfrak{Q}_\infty^4$. Moreover,  the solution enjoys the estimate \eqref{omessetsim122n} and the following properties:
\begin{enumerate}[(1)]
  \item the energy inequality in differential form: for a.e. $t>0$,
  \begin{align}\label{1.200xyx}
  \frac{\mm{d}}{\mm{d}t}\tilde{\mathcal{E}}+\mathcal{D}
\leqslant 0 \end{align}
for some functional $\tilde{\mathcal{E}}$, which belongs to $W^{1,\infty}(\mathbb{R}^+)${\footnote{Since $\tilde{\mathcal{E}}\in W^{1,\infty}(\mathbb{R}^+)$,  there exists a function $\underline{\tilde{\mathcal{E}}}\in W^{1,\infty}(\mathbb{R}^+)\cap AC(\mathbb{R}_0^+)$ such that $\underline{\tilde{\mathcal{E}}}= {\tilde{\mathcal{E}}}$ and $\underline{\tilde{\mathcal{E}}}'= {\tilde{\mathcal{E}}}'$ for a.e. $t\in \mathbb{R}^+$. }} and  is equivalent to $\mathcal{E} $ a.e. $t>0$.
  \item the stability estimate of total energy:  for a.e. $t>0$,
   \begin{align}\label{1.200} \mathcal{E}(t)+\int_0^t\mathcal{D}(\tau)\mm{d} \tau
\lesssim I^0.\end{align}
  \item  the exponential stability estimates:  for a.e. $t>0$,
    \begin{align}
   &\label{1.200n0}
  e^{c_1 t}(\|\eta_2(t)\|_3^2+\mathcal{E}_{\mm{p}}(t))+\int_0^t e^{c_1 \tau}\mathcal{D}_{\mm{p}}(\tau)\mm{d}t
\lesssim I^0, \\
&e^{c_1 t }\|\eta_1(t)-\eta^\infty_1\|_2\lesssim  \sqrt{I^0} \label{1.200xx}
\end{align}
for some positive constant $c_1$, where $\eta_1^\infty\in H^2$ only depends on $y_2$.
\end{enumerate}
\end{thm}
\begin{rem}
If the assumptions of \eqref{0102n} and \eqref{2020102241504} are replaced by
$$ \bar{\rho}'\leqslant 0\mbox{ in }\Omega\mbox{ and }|m|>0,$$
then the conclusions  in Theorem \ref{thm2} also hold. In addition, Theorem \ref{thm2} also holds for the case $g=0$.
\end{rem}
\begin{rem}
By the assumptions of $\bar{\rho}$, it is easy to check that
$$0<m_{\mm{C}}\leqslant \frac{h  }{ \pi }\sqrt{\frac{g\|\bar{\rho}'\|_{L^\infty}}{\lambda}},$$
please refer to (4.25) and Lemma 4.6 in \cite{JFJSARMA2019}. Thus, in view of Theorem \ref{thm2}, we see that a horizontal magnetic field can inhibit the  RT instability, if the strength of magnetic field is properly large. \emph{However it is not clear to the authors that whether non-horizontal magnetic fields can also inhibit the  RT instability based on the system \eqref{0101}}
\end{rem}
\begin{rem}
For the case $\Omega=\mathbb{R}\times (0,h)$, we also obtain a similar result, where the exponential stability estimates in \eqref{1.200n0} and \eqref{1.200xx} should be replaced by algebraic stability estimates. \emph{We will verify this assertion in a forthcoming paper.}
\end{rem}
\begin{rem}
For each fixed $t\in \mathbb{R}_0^+$, the solution $\eta(y,t)$ in Theorem \ref{thm2} belongs to $\in H^3_{*}$. Let $\zeta=\eta+y$, then $\zeta $  satisfies \eqref{20210301715x} and \eqref{20210301715} for each $t\in \mathbb{R}_0^+$. We denote the inverse transformation of $\zeta$ by $\zeta^{-1}$, and then define that
$$(\varrho, v,N,\beta)(x,t):=(\bar{\rho}(y_2)-\bar{\rho}(\zeta_2), u(y,t), m\partial_1\eta(y,t),q(y,t))_{y=\zeta^{-1}(x,t)}.$$
Consequently  $(\varrho,v,N,\beta)$ is a classical solution of the magnetic RT problem  \eqref{0103}--\eqref{0105} and enjoys stability estimates, which are similar to \eqref{1.200}--\eqref{1.200n0} for sufficiently small $\delta$.
\end{rem}

Next we roughly sketch the proof of Theorem \ref{thm2}, and the details will be presented in Section \ref{sec:global}. The key proof for the existence of global small solutions is to derive an  \emph{a priori} energy inequality of differential form \eqref{1.200xyx} for some energy functional $\tilde{\mathcal{E}}$, which is equivalent to $\mathcal{E}$. To this purpose, let $(\eta,u)$ be a solution to \eqref{01dsaf16asdfasf}, and satisfy,  for  some $T>0$,
\begin{align}
& \det(\nabla \eta+I)=1\mbox{ in }\Omega\times \overline{I_T},\label{aprpiosasfesnew}\\
&\sup_{t\in\overline{I_T}} \|(\eta,\partial_1\eta,u)(t)\|_4\leqslant {\delta} \in (0,1] . \label{aprpiosesnew}
\end{align}
For sufficiently small $\delta$, we first derive the horizontal-type energy inequality, and then the $\mm{cur}$-type energy inequality, see \eqref{202008250856} and \eqref{202005021632}.
 Summing up the both energy inequalities of  horizontal-type and $\mm{cur}$-type, we can arrive at the total energy inequality
\begin{align}\label{for:0202n}
\frac{\mm{d}}{\mm{d}t}\tilde{\mathcal{E}}+ \mathcal{D}\lesssim \sqrt{\mathcal{E}}
 \mathcal{D}
\end{align}
for some (total) energy functional $\tilde{\mathcal{E}}$, which is equivalent to $\mathcal{E}$ under \emph{the stability condition \eqref{2020102241504}}. In particular, \eqref{for:0202n} further implies \eqref{1.200xyx}, which  yields the \emph{priori} stability estimate \eqref{1.200}. Thanks to the \emph{priori} estimate \eqref{1.200} and  the unique local(-in-time) solvability  of the transformed MRT problem \eqref{01dsaf16asdfasf} in Proposition \ref{202102182115}, we immediately get the unique global solvability  of the problem  \eqref{01dsaf16asdfasf}.

In addition, similarly to \eqref{1.200xyx}, we can also establish the partial energy inequality
\begin{align}\label{202012262105}
\frac{\mm{d}}{\mm{d}t}\tilde{\mathcal{E}}_{\mm{p}}+ \mathcal{D}_{\mm{p}}\leqslant 0,
\end{align}
where $\tilde{\mathcal{E}}_{\mm{p}} $ is equivalent to $\mathcal{E}_{\mm{p}}$.
Thanks to the observation that $\|\eta\|_{1,3}\lesssim\| \eta\|_{2,3}$ by the horizontal periodicity, we immediately see that
\begin{align}
\mathcal{E}_{\mm{p}} \mbox{ is equivalent to }\mathcal{D}_{\mm{p}}. \nonumber
\end{align}
Thus  the  partial  energy inequality  \eqref{202012262105}, together with the above equivalence, immediately implies the exponential stability of partial  energy \eqref{1.200n0}. Finally, \eqref{1.200xx} can be easily deduced from   \eqref{1.200n0} by an asymptotic analysis method.

Recently Pan--Zhou--Zhu  investigated the well-posdeness of the equations of a viscous MHD fluid in a 3D periodic domain \cite{PRHZYZYARMA2018} under the assumption of that the initial data satisfies some odevity condition.
 Motivated by Pan--Zhou--Zhu's result, we have the following exponential stability of total energy.
\begin{cor}\label{cor1}
If additionally the initial data $(\eta^0,u^0)$ in Theorem \ref{thm2} satisfies the odevity conditions
\begin{align}
&\label{202005011004}
(\eta_1^0,u_1^0)(y_1,y_2)=-(\eta_1^0,u_1^0)(-y_1,y_2),\\[1mm]
&\label{202005011005}
(\eta_2^0,u_2^0)(y_1,y_2)=(\eta_2^0,u_2^0)(-y_1,y_2).
\end{align}
then the solution
$(\eta,u,q)$ in Theorem \ref{thm2}  also satisfies the odevity conditions
\begin{align}
&\label{202005012204}
(\eta_1,u_1)(y_1,y_2,t)=-(\eta_1,u_1)(-y_1,y_2,t),\
(\eta_2,u_2)(y_1,y_2,t)=(\eta_2,u_2)(-y_1,y_2,t),
\end{align}
and  enjoys the following  exponential stability of total energy: for a.e. $t>0$,
\begin{equation}\label{1.200n}
 {e}^{c_2 t}\mathcal{E}(t)+\int_0^t {e}^{c_2 \tau}\mathcal{D}(\tau)\mm{d}\tau
\lesssim I^0
\end{equation}for some positive constant $c_2$.
\end{cor}

The key idea to prove Corollary \ref{cor1}  is that $\mathcal{E}$ is equivalent to $\mathcal{D}$ under the  odevity conditions, and thus we immediately get the exponential stability \eqref{1.200n}. The detailed derivation  will be provided in Section \ref{202103171455}.

We can not expect the stability result for transformed MRT problem under the condition $|m|\in (0,m_{\mm{C}})$. In fact, this condition results in RT instability.
\begin{thm}[Instability]\label{thm1}
Let  $\bar{\rho}$ satisfy \eqref{0102}--\eqref{0102n} and $a\geqslant 0$.
If $|m|\in (0,m_{\mm{C}})$, the equilibria $(\bar{\rho},0,\bar{M})$
is unstable in the Hadamard sense, that is, there are positive constants $m_0$, $\epsilon$,  $\delta_0$, and $((\tilde{\eta}^0,\eta^\mm{r}),(\tilde{u}^0,u^\mm{r}))\in  H^{4}_{\mm{s}}$
such that for any $\delta\in (0,\delta_0]$ and the initial data
 $$ (\eta^0, u^0):=\delta(\tilde{\eta}^0,\tilde{u}^0)
 +\delta^2(\eta^\mm{r},u^\mm{r})\in (H^5_{\mm{s}}\cap H^5_{1}\cap H^5_*)\times H^5_{\mm{s}}, $$
there is a  unique solution $(\eta,u,q) $ to the transformed MRT problem \eqref{01dsaf16asdfasf}, where $(\eta,u,q) \in   {\mathfrak{H}}^{1,4}_{1,*,\tau}\times   \mathfrak{U}_{\tau}^{4}\times  \mathfrak{Q}_{\tau}^4$ for any $\tau\in I_{T^{\max}} $, $T^{\max}$ denotes the maximal time of existence of the solution, and $\mm{div}_{\mathcal{A}^0}u^0=0$ with $\mathcal{A}^0:=(\nabla \eta^0+I)^{-\mm{T}}$.
However, for $1\leqslant i$, $j\leqslant 2$, and $k=0$, $1$,
\begin{equation*}
\|\partial_j^k\chi_i(T^\delta)\|_{L^1} \geqslant \epsilon
\end{equation*}
for some escape time $T^\delta:=\frac{1}{\Lambda}\mm{ln}\frac{2\epsilon}{m_0\delta}\in
I_{T^{\max}}$, where $\chi$ can be taken by $ \eta $ or $u $.
\end{thm}

The proof of Theorem \ref{thm1} is based on a so-called bootstrap instability method. The bootstrap instability method has its origin in \cite{GYSWIC,GYSWICNonlinea}. Later, various versions of the bootstrap approach were presented by many authors, see \cite{FSSWVMNA,GYHCSDDC,JFJSZWC} for examples. In particular, recently Jiang--Jiang--Zhan proved the existence of the RT instability solution  under  $L^1$-norm for the stratified viscous,  non-resistive MHD fluids \cite{JFJSZWC} . In this paper, we adapt the version of the bootstrap instability method in \cite{JFJSZWC} to prove Theorem \ref{thm1}. It should be noted that the authors in \cite{JFJSZWC} considered the RT instability in viscous fluids. However our problem is the inviscid case, thus there exist some details, which are different to the Ref. \cite{JFJSZWC}, in the proof of Theorem \ref{thm1}. In particular, the absence of the strong continuity of highest order of spacial derivatives of $(\eta,u)$ results in some troubles. Fortunately, these troubles can be overcome by making use of the stability of local(-in-time) solutions \eqref{202012212151} and  the weak continuity, i.e. $ \nabla^4( \partial_1\eta,u)\in C^0_{B, \mm{weak}}(\overline{ I_{\tau}},L^2)$ for any $\tau\in I_{T^{\max}}$.

We mention that Jang--Guo ever proved the RT instability of inviscid fluids in a 2D periodic domain \cite{HHJGY}. It is not clear to authors that whether Jang--Guo's result's can be extended to the slab domain. In other word, \emph{it is not clear that whether Theorem \ref{thm1} also holds for the case $m=0$, such case will be further investigated in future}.

The rest of this paper is organized as follows.
In Sections \ref{sec:global}--\ref{sec:instable}, we provide the proof of  Theorem \ref{thm2}, Corollary \ref{cor1} and  Theorem \ref{thm1} in sequence. In Section \ref{202102241211}, we will establish the local well-posdedness result for the transformed MRT problem \eqref{01dsaf16asdfasf}.  Finally, in \ref{sec:09}, we list some well-known mathematical results, which will be used in Sections \ref{sec:global}--\ref{202102241211}.

\section{Proof of Theorem \ref{thm2}}\label{sec:global}
This section is devoted to the proof of Theorem \ref{thm2}.
The key step is to derive the  total energy inequality \eqref{1.200xyx} and the partial energy inequality \eqref{202012262105} for the transformed MRT problem \eqref{01dsaf16asdfasf} by  \emph{a priori} estimates. To this end, let $(\eta,u,q)$ be a solution to \eqref{01dsaf16asdfasf} and satisfy  \eqref{aprpiosasfesnew} and \eqref{aprpiosesnew},
where $ {\delta}$ is sufficiently small, and the smallness of
$\delta$ depends on $g$, $a$, $\lambda$, $m$, $\bar{\rho}$ and $\Omega$. It should be noted that $a$, $m$ and $\bar{\rho}$ satisfy the assumptions in Theorem \ref{thm2}. Next we proceed with  \emph{a priori} estimates.

\subsection{Preliminary estimates}

To being with, we shall establish some preliminary estimates involving $(\eta,u)$.
\begin{lem}
\label{201805141072}
For any given $t\in \overline{I_T}$, we have
\begin{enumerate}[(1)]
  \item the estimates for $\mathcal{A}$ and $\mathcal{A}_t$:
\begin{align}
&\label{aimdse}
\|\mathcal{A}\|_{C^0(\overline{\Omega})}+ \|\mathcal{A}\|_3 \lesssim  1,\\
&\label{06142100x} \|\mathcal{A}_t\|_{i,j} \lesssim   \|   u\|_{i,j+1} \ \mbox{for}\ 0\leqslant i+j\leqslant 3.
\end{align}
\item the estimates $\tilde{\mathcal{A}}$:
\begin{align}
&  \label{06041533fwqg}
\|\tilde{\mathcal{A}}\|_{i}\lesssim   \| \eta\|_{ i+1} \mbox{ for }0\leqslant i\leqslant 3,
\\
&  \label{06041533fwqgn}
\|\tilde{\mathcal{A}}\|_{i,j}\lesssim   \|\eta\|_{i,j+1} \mbox{ for }1\leqslant i+j\leqslant 4\mbox{ and }i \geqslant 1.
\end{align}
Here and in what follows $\tilde{\mathcal{A}}:=\mathcal{A}-I$.
\item  the estimate of $\mm{div}{ {u}}$: for $0\leqslant i\leqslant 3$ ,
\begin{align}
&\label{201808181500} \|\mm{div} u\|_{i}\lesssim \|(\eta,u)\|_{4}\|(\eta,u)\|_{1,i}.
\end{align}
\item the estimates involving gravity term: for sufficiently small,
\begin{align}
&\label{omessetsim122}
\|G_{\eta}\|_{3} \lesssim  \|\eta_2\|_{3},\\
&\left\| \mathcal{G}\right\|_{k,0}\lesssim
\begin{cases}
\|\eta_2\|_2\|\eta_2\|_{0}   &\mbox{for }k=0; \\
\|\eta_2\|_{3}\|\partial_1\eta_2\|_{\underline{k-1},0}&\mbox{for }1\leqslant k\leqslant 4,
\end{cases}\label{2022011130957}
\end{align}
where $\mathcal{G}:=g \left(G_{\eta}-\bar{\rho}'\eta_2\right)$.
\end{enumerate}
\end{lem}
\begin{pf}
(1)--(3)
Recalling \eqref{aprpiosasfesnew}   and the definitions of $\mathcal{A}$ and $\tilde{\mathcal{A}}$, we can compute out that
\begin{align}
\mathcal{A}=
\left(\begin{array}{c}
\partial_2\eta_2+1\quad-\partial_1\eta_2\\[1mm]
-\partial_2\eta_1\quad\partial_1\eta_1+1,
        \end{array}\right) \nonumber
\end{align}
and thus
\begin{align}
\tilde{\mathcal{A}}=
\left(\begin{array}{c}
\partial_2\eta_2 \quad-\partial_1\eta_2\\[1mm]
-\partial_2\eta_1\quad\partial_1\eta_1
        \end{array}\right). \label{202103181309}
\end{align}

Making use of \eqref{01dsaf16asdfasf}$_1$, \eqref{aprpiosesnew},  the embedding inequality \eqref{esmmdforinfty} and the product estimate \eqref{fgestims}, we can easily deduce \eqref{aimdse}--\eqref{201808181500} from \eqref{01dsaf16asdfasf}$_3$ and the expressions of  $\mathcal{A}$ and $\tilde{\mathcal{A}}$.

(4) We turn to verifying \eqref{omessetsim122} and \eqref{2022011130957}.
  By \eqref{aprpiosesnew} and Lemma \ref{pro:1221},  $\zeta:=\eta+y$ satisfies the diffeomorphism  properties \eqref{20210301715x} and \eqref{20210301715} for sufficiently small $\delta$.  Thus $\bar{\rho}^{(l)}(y_2+\eta_2)$ for any $y\in \overline{\Omega}$ makes sense, and
\begin{align}
\label{20200830asdfa2114}
\bar{\rho}^{(l)}(y_2+\eta_2)-\bar{\rho}^{(l)}(y_2) =  \int_{0}^{\eta_2} \bar{\rho}^{(l+1)}(y_2+z)\mm{d}z
\end{align}
for  $ 0\leqslant l\leqslant  3$.
Moreover, for  any given $t\in \overline{I_T}$,
\begin{align}
&\label{esmmdforiasfanfty}
\sup_{y\in \overline{\Omega}}
\left|\left.  \bar{\rho}^{(l+1)}\right|_{y_2=y_2+\eta_2}\right|\lesssim 1,\\
&\label{esmmdforiasdfaasfanfty}
 \sup_{y\in \overline{\Omega}} \sup_{z\in \Psi}\left| \bar{\rho}^{(l+1)}(y_2+z)\right| \lesssim 1,
\end{align}
where  $\Psi:= \{\tau~|~0\leqslant \tau \leqslant \eta_2\}$ for $\eta_2\geqslant 0$ and $:=(\eta_2,0]$ for $\eta_2< 0$.
Making use of \eqref{aprpiosesnew}, \eqref{20200830asdfa2114}--\eqref{esmmdforiasdfaasfanfty} and the embedding inequality \eqref{esmmdforinfty}, we easily deduce \eqref{omessetsim122} from the definition of $G_\eta$.

Since
\begin{align}\nonumber
\bar{\rho}^{(m)}(y_2+\eta_2)-\bar{\rho}^{(m)}(y_2) =\bar{\rho}^{(m+1)} (y_2)\eta_2+{ \int_{0}^{\eta_2}\left(\eta_2(y,t) -
z\right)\bar{\rho}^{(m+2)}(y_2+z)\mm{d}z}
\end{align}
 for $0\leqslant m\leqslant 2$, it is easy to see that
\begin{align}\nonumber
\left\|\mathcal{G}\right\|_0 \lesssim\|\eta_2\|_2\|\eta_2\|_{0}
\end{align}
and
\begin{align}
\left\|\partial_1^{k}\mathcal{G}\right\|_0
&=\left\|\partial_1^{k-1}\big(\bar{\rho}'(y_2+\eta_2)\partial_1\eta_2\big)-\bar{\rho}'\partial_1^{k}\eta_2\right\|_0\nonumber\\[1mm]
&\lesssim\left\|\big(\bar{\rho}'(y_2+\eta_2)-\bar{\rho}'\big)\partial_1^{k}
\eta_2\right\|_0
+\left\|[\partial_1^{k-1},\bar{\rho}'(y_2+\eta_2)]\partial_1\eta_2\right\|_0\nonumber\\[1mm]
&\lesssim\|\eta_2\|_{3}\|\partial_1\eta_2\|_{\underline{k-1},0}\mbox{ for } 1\leqslant k\leqslant 4.  \nonumber
\end{align}
Putting the above two estimates together yields \eqref{2022011130957}.
\hfill $\Box$
\end{pf}

\begin{lem}\label{lem:202012242115}
We have
\begin{enumerate}[(1)]
\item  the estimate of $\mm{div} { {\eta}}$: for $0\leqslant i\leqslant 3$,
\begin{align}
&\label{improtian1}
\|\mm{div} { {\eta}}\|_{j,i} \lesssim
\begin{cases}
\|\eta\|_{3} \|\eta\|_{1,i}&\mbox{for } j=0 ;\\[1mm]
\|\eta\|_{3} \|\eta\|_{1,i+1}&\mbox{for }  j=1 ,
\end{cases}\\
&\label{improtifdsafan1n0}
\|\mm{div} { {\eta}}\|_{j,i} \lesssim
\|\eta\|_{3} \|\eta\|_{2,2} \mbox{ for }i+j=4\mbox{ and } j\geqslant 2.
\end{align}
\item the estimate of $\eta_2$:
\begin{align}
&\label{omessetsim122n}
\|\eta_2\|_j\lesssim\begin{cases}
\|\eta\|_{1,0}&\mbox{for } j=0 ,\ 1;\\[1mm]
 \|\eta\|_{1,j-1}&\mbox{for }  2\leqslant j\leqslant 4.
\end{cases}
\end{align}
\end{enumerate}
\end{lem}
\begin{pf}
(1) Recalling \eqref{aprpiosasfesnew}, we can compute out that
\begin{align}
&\mm{div}\eta=\partial_1\eta_2\partial_2\eta_1-\partial_1\eta_1\partial_2\eta_2.
\label{20202101131222}
\end{align}
Exploiting   the product estimate \eqref{fgestims}, we can easily deduce \eqref{improtian1}--\eqref{improtifdsafan1n0} from  the relation  above.

 (2) Noting that
\begin{align}
\eta_2|_{\partial\Omega}=0 \mbox{ and }\partial_2\eta_2=\mm{div}\eta-\partial_1\eta_1,
\label{2020101131413}
\end{align} thus, using \eqref{2020101131413} and \eqref{poinsafdcasfaressadf1}, we get
\begin{align}
&\nonumber \|\eta_2\|_0\lesssim\| (\partial_1\eta_1,\mm{div}\eta)\|_{0},\\
& \nonumber
\|\eta_2\|_{j}\lesssim\|\eta_2\|_{0}+\|\nabla\eta_2\|_{j-1}\lesssim\|(\partial_1\eta
,\mm{div}\eta)\|_{j-1}\mbox{ for }1\leqslant j\leqslant 4.
\end{align}
We immediately get \eqref{omessetsim122n} from the two estimates above and \eqref{improtian1}.
\hfill $\Box$
\end{pf}
\begin{lem}\label{lem:07291200} Let the multiindex $\alpha$ satisfy $|\alpha|\leqslant 3$ and
\begin{align}
\label{2021011281929}
{W}^{\alpha} = \partial^{\alpha}(\mm{curl}_{\mathcal{A}_{t}}\left(\bar{\rho}u\right)
-\lambda m^2\mm{curl}_{\partial_1\mathcal{A}} \partial_1\eta) .
\end{align}
\begin{enumerate}[(1)]
  \item
Then we have
\begin{align}
&\label{201907291200}
 \|{W}^{\alpha} \|_0\lesssim\|(\partial_1\eta,u)\|_{4}^2,\\
&\label{201907291432}
 \|\partial^{\alpha} (\operatorname{curl}_{\mathcal{A}} (\bar{\rho}\chi )
-\operatorname{curl} (\bar{\rho}\chi) ) \|_0\lesssim
\|\eta\|_4\|\chi\|_4,
\end{align}
where $\chi=\eta$, $\partial_1\eta$ and $u$.
  \item  For  $\alpha_2\neq 2$,
we have
\begin{align}\label{201907291200n0}
& \|{W}^{\alpha} \|_0\lesssim\|(\partial_1\eta,u)\|_{4}
{ {\|(\partial_1\eta,u)\|_{\underline{1},3}}},\\
&\label{201907291432n0}
 \|\partial^{\alpha} (\operatorname{curl}_{\mathcal{A}}  (\bar{\rho}u )
-\operatorname{curl} (\bar{\rho}u ) ) \|_0\lesssim
\|(\eta,u)\|_4(\| \eta \|_{2,2}+\| u\|_{1,3}  ) ,\\[1mm]
&\label{202005031815n0}
 \|\partial^{\alpha} (\operatorname{curl}_{\mathcal{A}}
(\bar{\rho}\partial_1^i \eta )-\operatorname{curl}(\bar{\rho}\partial_1^i \eta ) ) \|_0\lesssim
\|\eta\|_{4}\| \eta\|_{2,3}\mbox{ for }i=0,\ 1.
\end{align}
\end{enumerate}
\end{lem}
\begin{pf}
The estimates \eqref{201907291200} and \eqref{201907291200n0} can be easily derived by using  \eqref{06142100x}, \eqref{06041533fwqgn} and the product estimate
 \eqref{fgestims}. Noting that
$$\operatorname{curl}_{\mathcal{A}}w-\operatorname{curl}w=
\tilde{\mathcal{A}}_{1k}\partial_{k}w_2-\tilde{\mathcal{A}}_{2k}\partial_{k}w_1,$$ then
\eqref{201907291432} and \eqref{201907291432n0}--\eqref{202005031815n0} can be easily derived by using   \eqref{06041533fwqg}, \eqref{202103181309}, \eqref{fgestims} and \eqref{202012241002}.
\hfill$\Box$
\end{pf}

\subsection{Horizontal-type energy inequalities}\label{subsec:Horizon}
This section is devoted to establishing the total/partial horizontal-type energy inequalities.
Let $0\leqslant k\leqslant 4$, then we apply $\partial_1^{k}$  to \eqref{01dsaf16asdfasf} to get that
\begin{equation}\label{01dsaf16asdfasf03n}
                              \begin{cases}
\partial_1^{k}\eta_t=\partial_1^{k}u ,\\[1mm]
\bar{\rho}\partial_1^{k}u_t+\partial_1^{k}\nabla_{\mathcal{A}} q+a\bar{\rho} \partial_1^{k}u
=\lambda m^2 \partial_1^{k+2}\eta+g\bar{\rho}'\partial_1^{k}\eta_2\mathbf{e}_2
+\partial_1^{k}\mathcal{G}\mathbf{e}_2 ,\\[1mm]
\partial_1^{k}\div_{\mathcal{A}} u=0, \\[1mm]
(\partial_1^{k}\eta,\partial_1^{k}u)|_{\partial\Omega}\cdot\vec{\mathbf{n}}=0 .
\end{cases}
\end{equation}
Thus we can derive from \eqref{01dsaf16asdfasf03n} the following horizontal spatial estimates of $\left(\eta, u\right)$.
\begin{lem}\label{lem:08241445}
For $0\leqslant k\leqslant 4$,
\begin{align}
&
\frac{\mm{d}}{\mm{d}t}\left(\int\bar{\rho}\partial_1^{k}\eta\cdot\partial_1^{k} u\mm{d}y
+\frac{a}{2}\|\sqrt{\bar{\rho}}\partial_1^{k} \eta\|_{0}^2\right)
-E(\partial_1^{k}\eta)
\lesssim\|\sqrt{\bar{\rho}}\partial_1^{k}u\|_{0}^2
+{ \sqrt{\mathcal{E}}\mathcal{D}_{\mm{p}}},
\label{202008241446}
\\[1mm]
&\label{202008241448}
 \frac{\mm{d}}{\mm{d}t}\left(\|\sqrt{ \bar{\rho} }  u\|^2_{k,0}
-E(\partial_1^{k} \eta) \right)+ c\|\sqrt{\bar{\rho}}  u\|_{k,0}^2
\lesssim \sqrt{\mathcal{E}}\mathcal{D}_{\mm{p}}.
\end{align}
\end{lem}
\begin{pf}
Multiplying \eqref{01dsaf16asdfasf03n}$_2$ by $\partial_1^{k}\eta$, resp. $\partial_1^{k}u$ in $L^2$,
then, using the integrating by parts and \eqref{01dsaf16asdfasf03n}$_1$,
we have that
\begin{align}
&\frac{\mm{d}}{\mm{d}t}\left(\int\bar{\rho}\partial_1^{k}\eta\cdot\partial_1^{k} u\mm{d}y
+\frac{a}{2}\|\sqrt{\bar{\rho}}\partial_1^{k} \eta\|_{0}^2\right)
-E(\partial_1^{k}\eta)\nonumber \\[1mm]
&=\|\sqrt{\bar{\rho}}\partial_1^{k}u\|_{0}^2
+\int\partial_1^{k} \mathcal{G}\partial_1^{k} \eta_2\mm{d}y-\int \partial_1^{k}\nabla_{\mathcal{A}}q\cdot\partial_1^{k} \eta\mm{d}y,\label{202008241510}
\end{align}
 resp.
\begin{align}
&
\frac{1}{2}\frac{\mm{d}}{\mm{d}t}\left(\|\sqrt{\bar{\rho}} u\|_{k,0}^2
-E(\partial_1^{k}\eta)\right)+a\|\sqrt{\bar{\rho}} u\|_{k,0}^2\nonumber \\[1mm]
&=\int\partial_1^{k} \mathcal{G}\partial_1^{k} u_2\mm{d}y-\int \partial_1^{k}\nabla_{\mathcal{A}}q \cdot\partial_1^{k}u\mm{d}y
.\label{202008241510n}
\end{align}

By \eqref{2022011130957}, \eqref{omessetsim122n} and \eqref{202012241002}, we can estimate that
\begin{align}
&\label{202008241546}
\int\partial_1^{k} \mathcal{G}\partial_1^{k} \eta_2\mm{d}y\leqslant
\|  \mathcal{G}\|_{k,0}\|  \eta_2\|_{k,0}\lesssim
\sqrt{\mathcal{E}}\mathcal{D}_{\mm{p}},\\[1mm]
&\label{202008241624}
\int\partial_1^{k} \mathcal{G} \partial_1^{k} u_2\mm{d}y\leqslant
\|  \mathcal{G}\|_{k,0}\|  u_2\|_{k,0}\lesssim
\sqrt{\mathcal{E}}\mathcal{D}_{\mm{p}}.
\end{align}
Next we estimate for the integrals involving the pressure in \eqref{202008241510}--\eqref{202008241510n} by two cases.

(1) We first consider the case $k=0$.

By the integration by parts and the boundary-value condition \eqref{01dsaf16asdfasf03n}$_4$ with $k=0$,
\begin{align}\nonumber
-\int\nabla_{\mathcal{A}}q\cdot\eta\mm{d}y
 =\int q\mm{div}\eta\mm{d}y-\int\nabla_{\tilde{\mathcal{A}}}q\cdot\eta\mm{d}y.
\end{align}
Next we estimate for the two integrals on the right hand of the above identity.

Let $K:=\eta_1 (         - \partial_2\eta_2  ,    \partial_1\eta_2 )^{\mm{T}}$. Then the identity \eqref{20202101131222} can be rewritten as follows
\begin{align}\label{202008280955nq}
\mm{div}\eta=\mm{div}K.
\end{align}
By \eqref{omessetsim122n}, the  embedding inequality \eqref{esmmdforinfty} and \eqref{202012241002}, we have
\begin{equation}\label{201910040902n00}
\|K\|_0 { \lesssim\|\eta\|_2 \| \eta_2\|_{1}} \lesssim\|\eta\|_2 \| \eta\|_{2,1} ,
\end{equation}
Exploiting  \eqref{202008280955nq}  and \eqref{201910040902n00}, we have
\begin{align}
\int q\mm{div} \eta\mm{d}y =   -\int K\cdot\nabla   q\mm{d}y
\leqslant\|K\|_0\|\nabla q\|_0
\lesssim \sqrt{\mathcal{E}}{ \mathcal{D}_{\mm{p}}}. \label{201910040902n01}
\end{align}

Noting that
$$\nabla_{\tilde{\mathcal{A}}}q\cdot\eta
=\eta_1(\partial_2\eta_2\partial_1q-\partial_1\eta_2\partial_2q)
+\eta_2(\partial_1\eta_1\partial_2q-\partial_2\eta_1\partial_1q ),$$
thus, using \eqref{omessetsim122n}, \eqref{fgestims}  and \eqref{202012241002}, we deduce that
\begin{equation}\label{201910040902n02}
\left|\int \nabla_{\tilde{\mathcal{A}}}q\cdot\eta\mm{d}y\right|
\lesssim \|\eta\|_2 \| \eta\|_{1,1}\|\nabla q\|_0 \lesssim \sqrt{\mathcal{E}}{ \mathcal{D}_{\mm{p}}}.
\end{equation}
Combining with \eqref{201910040902n01} and \eqref{201910040902n02}, we get
\begin{align}\label{202008241635n03}
\left|\int\nabla_{\mathcal{A}}q\cdot \eta\mm{d}y\right|\lesssim \sqrt{\mathcal{E}}{ \mathcal{D}_{\mm{p}}}.
\end{align}
Putting \eqref{202008241546} with $k=0$ and \eqref{202008241635n03} into \eqref{202008241510} yields  \eqref{202008241446} with $k=0$.

In addition, by the integration by parts, \eqref{01dsaf16asdfasf03n}$_3$ and \eqref{01dsaf16asdfasf03n}$_4$ with $k=0$, we have
\begin{align}\label{202008241635}
-\int\nabla_{\mathcal{A}}q\cdot u\mm{d}y=\int q\mm{div}_{\mathcal{A}}u\mm{d}y=0.
\end{align}
Thus putting \eqref{202008241624} with $k=0$ and \eqref{202008241635} into \eqref{202008241510n} with $k=0$ yields \eqref{202008241448} with $k=0$.

(2) Now we further consider the case $k\neq0$.

Making use of \eqref{improtian1}, \eqref{improtifdsafan1n0} and \eqref{202012241002},  we can   estimate that
\begin{align}
-\int\partial_1^{k}\nabla_{\mathcal{A}}q \cdot\partial_1^{k} \eta\mm{d}y
=&\int\partial_1^{k-1}\nabla_{\tilde{\mathcal{A}}}q \cdot\partial_1^{k+1} \eta\mm{d}y
+\int\partial_1^{k}q \mm{div}\partial_1^{k} \eta\mm{d}y\nonumber \\[1mm]
\leqslant&
\|\partial_1^{k-1}\nabla_{\tilde{\mathcal{A}}}q\|_0\|\partial_1^{k+1} \eta\|_{0}
+\|\partial_1^{k}q\|_0\|\partial_1^{k}\mm{div}\eta\|_0\nonumber \\[1mm]
\lesssim&\sqrt{\mathcal{E}}\mathcal{D}_{\mm{p}}. \label{202008241550}
\end{align}
Putting \eqref{202008241546} and \eqref{202008241550} into \eqref{202008241510} yields \eqref{202008241446} for $k\neq 0$.

Exploiting \eqref{01dsaf16asdfasf03n}$_3$, we have
 \begin{align} \label{202008241550n}
&-\int\partial_1^{k}\nabla_{\mathcal{A}}q\cdot\partial_1^{k}u\mm{d}y
=\int\partial_1^{k}q\mm{div}_{\mathcal{A}}\partial_1^{k}u\mm{d}y
-\int[\partial_1^{k}, \mathcal{A}]\nabla q\cdot\partial_1^{k}u\mm{d}y \nonumber \\[1mm]
&=\int\partial_1^{k}q\partial_1^{k} \mm{div}_{\mathcal{A}}u \mm{d}y
-\int\partial_1^{k}q[\partial_1^{k},\mathcal{A}^{\mm{T}}]:\nabla u\mm{d}y
-\int[\partial_1^{k}, \mathcal{A}]\nabla q\cdot\partial_1^{k}u\mm{d}y\nonumber \\[1mm]
&=-\int\partial_1^{k}q[\partial_1^{k},\mathcal{A}^{\mm{T}}]:\nabla u\mm{d}y
-\int[\partial_1^{k}, \mathcal{A}]\nabla q\cdot\partial_1^{k}u\mm{d}y \lesssim\sqrt{\mathcal{E}}\mathcal{D}_{\mm{p}}.
\end{align}
Consequently, putting \eqref{202008241624} and \eqref{202008241550n} into \eqref{202008241510n} yields \eqref{202008241448}  with $k\neq 0$.
The proof is completed.
\hfill $\Box$
\end{pf}

Now we shall establish stabilizing estimates for $E(\partial_1^{k}\eta)$ appearing in \eqref{202008241446} and \eqref{202008241448}.
\begin{lem}\label{lem:08250749}
It holds that
\begin{align}
\label{202008250745}
\|\eta\|_{{i+1},0}^2
\lesssim-E(\partial_1^{i}\eta)+\|\eta\|_4\| \eta\|_{2,3}^2
\mbox{ for any }0\leqslant i\leqslant  4.
\end{align}
\end{lem}
\begin{pf}
By the definition of $m_{\mm{C}}$, it is easy to see that
\begin{align*}
-\int g\bar{\rho}'w_2^2\mm{d}y\geqslant-\lambda m_{\mm{C}}^2\|\partial_1w\|_0^2\mbox{ for any }w\in H_{\sigma}^1,
\end{align*}
 which, together with the stability condition $|m|>m_{\mm{C}}$, implies that
\begin{align}\label{202008250814}
\| w\|_{1,0}^2\lesssim\lambda (m^2-m_{\mm{C}}^2)\|\partial_1w\|_0^2\lesssim-E(w)
 \mbox{ for any } w\in H_{\sigma}^1.
\end{align}

Let us consider the Stokes problem
\begin{equation}\label{202008250820}
                              \begin{cases}
-\Delta\tilde{\eta}+\nabla\varpi=0,\  \div \tilde{\eta}=\mm{div}\eta  &\mbox{in } \Omega, \\[1mm]
\tilde{\eta}=0  &\mbox{on } \partial\Omega.
\end{cases}
\end{equation}
By the existence theory of Stokes problem, there exist a unique solution $(\tilde{\eta},\varpi)\in H^4\times \underline{H}^3$
to \eqref{202008250820}. Moreover, $\partial_1^j( \tilde{\eta},\varpi)$ is
also the solution of Stokes problem for $1\leqslant j\leqslant 3$ and
\begin{align}\label{11260850}
\| \tilde{\eta}\|_{j,2}
\lesssim\| \mm{div}{\eta}\|_{j,1}
\lesssim
\|\eta\|_4\| \eta\|_{2,3}  ,
\end{align}
where we have used \eqref{improtian1}, \eqref{improtifdsafan1n0} and \eqref{202012241002} in the  last inequality above.

Now we use $\partial_1^i(\eta- \tilde{\eta})$
to rewrite $ {E}(\partial_1^{i}\eta)$ as follows:
\begin{align}\label{201912191626}
 {E}(\partial_1^{i}\eta)=
 {E}( \partial_1^{i}(\eta-\tilde{\eta}) )
+{E}(\partial_1^{i}\tilde{\eta} )
-I_i,
\end{align}
where
\begin{equation*}\label{dabddfsffsmms}
\begin{aligned}
I_i:=
2\lambda m^2\int \partial_1^{i+1}\eta\cdot \partial_1^{i+1}\tilde{\eta}\mm{d}y
-2g\int \bar{\rho}'\partial_1^{i}\eta_2\partial_1^{i}\tilde{\eta}_2\mm{d}y.
\end{aligned}
\end{equation*}
Note that $\partial_1^{i}(\eta- \tilde{\eta})\in H_{\sigma}^1$, thus, we use \eqref{202008250814} to get
\begin{align}\nonumber
\| \eta- \tilde{\eta} \|_{i+1,0}^2
\lesssim-E(\partial_1^{i}(\eta-\tilde{\eta})),
\end{align}
which, together with   \eqref{201912191626} and Young's inequality, yields
\begin{align}\label{202008250814nn}
\|   \eta  \|_{i+1,0}^2
\lesssim{E}(\partial_1^{i}\tilde{\eta} ) -E(\partial_1^{i}\eta)
-I_i +\|   \tilde{\eta} \|_{i+1,0}^2.
\end{align}

Making use of  \eqref{omessetsim122n}, \eqref{11260850} and  \eqref{202012241002}, we can estimate that
\begin{align}
 {E}(\partial_1^{i}\tilde{\eta} )
-I_i
+\| \tilde{\eta}\|_{i+1,0}^2   \lesssim
\|\eta\|_4 \|\eta\|_{2,3}^2 \mbox{ for }0\leqslant i\leqslant 4.  \nonumber
\end{align}
Finally, putting the above estimate into  \eqref{202008250814nn} yields
\eqref{202008250745}. This completes the proof.
\hfill $\Box$
\end{pf}

Thanks to Lemmas \ref{lem:08241445}--\ref{lem:08250749} and \eqref{omessetsim122n},  we easily get the total/partial horizontal-type energy inequalities.
\begin{lem}\label{pro0902}
There exist  two functionals
 $\underline{\mathcal{E}}$ and $\underline{\mathfrak{E}}$ of $(\eta,u)$  such that
\begin{align}
&\label{202008250856}
\frac{\mm{d}}{\mm{d}t}
\underline{\mathcal{E}}
+c \| (u,\partial_1\eta)\|_{\underline{4},0}^2\lesssim\sqrt{\mathcal{E}}\mathcal{D}_{\mm{p}},
\\[1mm]
&\label{202008250856n0}
\frac{\mm{d}}{\mm{d}t}
\underline{\mathfrak{E}}
+c({ \| u\|_{0}^2}+ \|\partial_1(u,\partial_1\eta)\|^2_{\underline{3},0}
) \lesssim\sqrt{\mathcal{E}}\mathcal{D}_{\mm{p}},\\
&\|(\eta, \partial_1\eta, u )\|_{\underline{4},0}^2- {c} \| \eta\|_3\|\eta\|_{2,3}^2
\lesssim \underline{\mathcal{E}} \lesssim
\|(\eta, \partial_1\eta, u )\|_{\underline{4},0}^2   , \label{20103231012} \\
&\| u \|_{0}^2+ \|\partial_1(\eta, \partial_1\eta, u )\|_{\underline{3},0}^2-  {c} \| \eta\|_4\|\eta\|_{2,3}^2 \lesssim  \underline{\mathfrak{E}}  \lesssim  \|  u \|_{0}^2+ \|\partial_1(\eta, \partial_1\eta, u )\|_{\underline{3 },0}^2  .
\label{2022103261057}
\end{align}
 \end{lem}

\subsection{Curl-type energy inequality}\label{subsec:Vor}

This section is devoted to establishing curl-type energy inequalities.
 By \eqref{dstist01}, we have
$$\mm{curl}_{\mathcal{A}}\left(gG_{\eta}\mathbf{e}_2\right)
=\mm{curl}_{\mathcal{A}}\left(-g\bar{\rho}\mathbf{e}_2\right)
=\mathcal{A}_{1j}\partial_j\left(-g\bar{\rho}\right)= g\bar{\rho}'\partial_1\eta_2 .
$$
Thus applying $\mm{curl}_{\mathcal{A}}$ to \eqref{01dsaf16asdfasf}$_2$ yields that
\begin{align}\label{201910072117}
\partial_t\mm{curl}_{\mathcal{A}}\left(\bar{\rho}u\right)+a\mm{curl}_{\mathcal{A}}\left(\bar{\rho}u\right)
=\lambda m^2\mm{curl}_{\mathcal{A}}\left( {\bar{\rho}}^{-1}\partial_1^2\left(\bar{\rho}\eta\right)\right)
+ g\bar{\rho}'\partial_1 \eta_2 +\mm{curl}_{\mathcal{A}_{t}}\left(\bar{\rho}u\right).
\end{align}

Let the multiindex $\alpha$ satisfy $|\alpha|\leqslant 2$.
Recalling the definition of ${W}^{\alpha}$ in \eqref{2021011281929}, thus applying $\partial^{\alpha}$  to \eqref{201910072117} yields
\begin{align}
&\partial_t\partial^{\alpha}\mm{curl}_{\mathcal{A}}\left(\bar{\rho}u\right)
+a\partial^{\alpha}\mm{curl}_{\mathcal{A}}\left(\bar{\rho}u\right)\nonumber \\
&= \partial_1({\lambda m^2}({\bar{\rho}}^{-1}\partial^{\alpha}\mm{curl}_{\mathcal{A}} \left(\bar{\rho}\partial_1\eta\right)
+ [\partial^{\alpha}\mm{curl}_{\mathcal{A}},
{\bar{\rho}}^{-1}]\left(\bar{\rho}\partial_1\eta\right))
+ \partial^{\alpha}\left(g\bar{\rho}'\eta_2\right))
+{W}^{\alpha}.\label{202005021542}
\end{align}

Now we derive the following energy estimates for  $\mm{curl}u$ and $\mm{curl}\partial_1\eta$.
\begin{lem}\label{2019100216355nnn}
For multiindex $\alpha$ satisfying $|\alpha|\leqslant 3$,
we have
\begin{align}
&\label{202005021600}
\frac{\mm{d}}{\mm{d}t}
\int\left(2\partial^{\alpha}\mm{curl}_{\mathcal{A}}\left(\bar{\rho}\eta\right)
 \partial^{\alpha}\mm{curl}_{\mathcal{A}} \left(\bar{\rho}u\right)
+ {a} |\partial^{\alpha}\mm{curl}_{\mathcal{A}} \left(\bar{\rho}\eta\right) |^2\right)\mm{d}y
 +c\| \partial^{\alpha}\mm{curl}_{\mathcal{A}}
\left(\bar{\rho}\partial_1\eta\right)\|^2_{0}\nonumber\\
&\lesssim\|\partial^{\alpha}\mm{curl}_{\mathcal{A}}\left(\bar{\rho}u\right)\|^2_{0}
+\| \lambda m^2[\partial^{\alpha}\mm{curl}_{\mathcal{A}}, {\bar{\rho}}^{-1}]\left(\bar{\rho}\partial_1\eta\right)
+\partial^{\alpha}\left(g\bar{\rho}'\eta_2\right) \|_{0}^2 \nonumber
\\
&\qquad +
           \begin{cases}
\sqrt{\mathcal{E}}\mathcal{D};\\
 \sqrt{\mathcal{E}}\mathcal{D}_{\mm{p}}\mbox{ for }\alpha_1\geqslant1,
           \end{cases} \\
&\label{202005021610}
 \frac{\mm{d}}{\mm{d}t}
 \left(\|\partial^{\alpha}\mm{curl}_{\mathcal{A}}\left(\bar{\rho}u\right) \|^2_0
+\lambda\| m\sqrt{\bar{\rho}^{-1}} \partial^{\alpha}\mm{curl}_{\mathcal{A}}
\left(\bar{\rho}\partial_1\eta\right)\|^2_0\right)
 +c\|\partial^{\alpha}\mm{curl}_{\mathcal{A}}\left(\bar{\rho}u\right)\|^2_{0}\nonumber\\
&\lesssim
\| \lambda m^2 [\partial^{\alpha}\mm{curl}_{\mathcal{A}}, {\bar{\rho}}^{-1}]\left(\bar{\rho}\partial_1^2\eta\right)
+\partial^{\alpha}\left(g\bar{\rho}'\partial_1\eta_2\right) \|_{0}^2
 +
           \begin{cases}
\sqrt{\mathcal{E}}\mathcal{D}; &\\
 \sqrt{\mathcal{E}}\mathcal{D}_{\mm{p}}\mbox{ for }|\alpha|\leqslant1\mbox{ or }\alpha_1\geqslant 1.
           \end{cases}
\end{align}
\end{lem}
\begin{pf}
(1) Multiplying \eqref{202005021542} by $\partial^{\alpha}\mm{curl}_{\mathcal{A}}\left(\bar{\rho}\eta\right)$ in $L^2$ yields that
\begin{align}
&\frac{\mm{d}}{\mm{d}t}
\int\left(\partial^{\alpha}\mm{curl}_{\mathcal{A}}\left(\bar{\rho}\eta\right)
 \partial^{\alpha}\mm{curl}_{\mathcal{A}} \left(\bar{\rho}u\right)
+\frac{a}{2}|\partial^{\alpha}\mm{curl}_{\mathcal{A}}\left(\bar{\rho}\eta\right)|^2\right)\mm{d}y
+{\lambda}m^2\| {\sqrt{\bar{\rho}}^{-1}}\partial^{\alpha}\mm{curl}_{\mathcal{A}} \left(\bar{\rho}\partial_1 \eta\right)\|^2_{0}\nonumber\\
&=I_1^{\alpha}+I_2^{\alpha} +\|\partial^{\alpha}\mm{curl}_{\mathcal{A}} \left(\bar{\rho}u\right)\|^2_{0}\nonumber \\
&\quad -\int\left(\lambda m^2[\partial^{\alpha}\mm{curl}_{\mathcal{A}}, {\bar{\rho}}^{-1}]\left(\bar{\rho}\partial_1\eta\right)
+\partial^{\alpha}\left(g\bar{\rho}'\eta_2\right)\right)
 \partial^{\alpha}\mm{curl}_{\mathcal{A}}\left(\bar{\rho}\partial_1\eta\right)\mm{d}y
,\label{202005031720}
\end{align}
where
\begin{align}
I_1^{\alpha}:=&
\int\big(\partial^{\alpha}\mm{curl}_{\mathcal{A}_t} \left(\bar{\rho}\eta\right)
 \partial^{\alpha}\mm{curl}_{\mathcal{A}}  \left(\bar{\rho}u\right)
+{W}^{\alpha}
 \partial^{\alpha}\mm{curl}_{\mathcal{A}}  \left(\bar{\rho}\eta\right)\big)\mm{d}y\nonumber\\
&-
\int\left(\partial^{\alpha}\left(g\bar{\rho}'\eta_2\right)
+\lambda m^2 \partial^{\alpha}\mm{curl}_{\mathcal{A}} \partial_1\eta \right)
 \partial^{\alpha}\mm{curl}_{\partial_1\mathcal{A}} \left(\bar{\rho}\eta\right)\mm{d}y, \nonumber  \\
I_2^{\alpha}:=&a \int  \partial^{\alpha}\mm{curl}_{\mathcal{A}_t}  \left(\bar{\rho}\eta\right)\partial^{\alpha}\mm{curl}_{\mathcal{A}} \left(\bar{\rho}\eta\right)\mm{d}y. \nonumber
\end{align}

Making use of  \eqref{aimdse}, \eqref{06142100x}, \eqref{omessetsim122n}, \eqref{201907291200}, \eqref{201907291200n0}, \eqref{fgestims} and \eqref{202012241002}, we can estimate that
\begin{align}
&\label{202008231718}
I_1^{\alpha}\lesssim
           \begin{cases}
\sqrt{\mathcal{E}}\mathcal{D}; &\\
 \sqrt{\mathcal{E}}\mathcal{D}_{\mm{p}}\mbox{ for }\alpha_1\geqslant1
           \end{cases}
\end{align}
and
\begin{align}
I_2^\alpha \lesssim
 \sqrt{\mathcal{E}}\mathcal{D}_{\mm{p}}\mbox{ for }\alpha_1\geqslant1.
\end{align}

In addition, noting the structure
$$
  \mm{curl}_{\mathcal{A}_t}  \left(\bar{\rho}\eta\right)
= \partial_t\tilde{\mathcal{A}}_{1j}
\partial_j\left(\bar{\rho}\eta_2\right)-\partial_t\tilde{\mathcal{A}}_{2j}
\partial_j\left(\bar{\rho}\eta_1\right)
$$
thus, similarly to \eqref{202008231718} with further using integral formula by parts, \eqref{01dsaf16asdfasf}$_1$ and \eqref{202103181309}, we easily estimate that
\begin{align}\label{202012251428}
I_2^\alpha
\lesssim
\sqrt{\mathcal{E}}\mathcal{D}.
\end{align}

Finally, putting \eqref{202008231718}--\eqref{202012251428} into \eqref{202005031720}, and then using  Young's inequality and the lower-bound condition $\inf_{y\in \overline{\Omega}}\bar{\rho}>0$, we get
\eqref{202005021600} immediately.

(2)
Multiplying \eqref{202005021542} by $\partial^{\alpha}\mm{curl}_{\mathcal{A}}\left(\bar{\rho}u\right)$ in $L^2$ yields that
\begin{align}
&\frac{1}{2}\frac{\mm{d}}{\mm{d}t}
(\|\partial^{\alpha}\mm{curl}_{\mathcal{A}} \left(\bar{\rho}u\right)\|^2_0
+{\lambda }\|m\sqrt{{\bar{\rho}}^{-1}}\partial^{\alpha}\mm{curl}_{\mathcal{A}}\left(\bar{\rho}\partial_1\eta\right)\|^2_0
) +a\|\partial^{\alpha}\mm{curl}_{\mathcal{A}}\left(\bar{\rho}u\right)\|^2_{0}\nonumber \\
&=\int\left(\lambda m^2\big[\partial^{\alpha}\mm{curl}_{\mathcal{A}},{\bar{\rho}}^{-1}]\left(\bar{\rho}\partial_1^2\eta\right)
+\partial^{\alpha}\left(g\bar{\rho}'\partial_1\eta_2\right)\right)
 \partial^{\alpha}\mm{curl}_{\mathcal{A}}\left(\bar{\rho}u\right)\mm{d}y
+I_3^{\alpha},\label{202005031742}
\end{align}
where we have defined that
\begin{align}
I_3^{\alpha}:=&\lambda m^2
\int {\bar{\rho}}^{-1}\partial^{\alpha}\mm{curl}_{\mathcal{A}} \left(\bar{\rho}\partial_1\eta\right)
\partial^{\alpha}\left(\mm{curl}_{\mathcal{A}_t} \left(\bar{\rho}\partial_1\eta\right)
- \mm{curl}_{\partial_1\mathcal{A}}  \left(\bar{\rho}u\right)\right)\mm{d}y\nonumber \\
& +\int
(\lambda m^2 [\partial^{\alpha}\mm{curl}_{\partial_1\mathcal{A}},
{\bar{\rho}}^{-1}]\left(\bar{\rho}\partial_1\eta\right)
+ {W}^{\alpha}) \partial^{\alpha}\mm{curl}_{\mathcal{A}} \left(\bar{\rho}u\right)\mm{d}y.  \nonumber
\end{align}
Similarly to \eqref{202008231718},
it is easy to estimate that
\begin{align}\label{202008231718n}
I_3^{\alpha}\lesssim
           \begin{cases}
\sqrt{\mathcal{E}}\mathcal{D} ;\\
 \sqrt{\mathcal{E}}\mathcal{D}_{\mm{p}}\mbox{ for }|\alpha| \leqslant1\mbox{ or }\alpha_1\geqslant1.
           \end{cases}
\end{align}
Plugging \eqref{202008231718n} into \eqref{202005031742}
and then using  Young's inequalities, we obtain \eqref{202005021610}.
This completes the proof.
\hfill$\Box$
\end{pf}

 Now we use Lemmas \ref{2019100216355nnn} to further derive the total/partial vorticity-type energy inequalities.
\begin{pro}\label{pro1632}
There exist two functionals
 $\mathcal{E}^{\mm{cul}}_{\mm{p}}$ and $\mathcal{E}^{\mm{cul}}$ of $(\eta,u)$ such that
\begin{align}
\label{202005021632}
&\frac{\mm{d}}{\mm{d}t}
\mathcal{E}^{\mm{cul}}
 +c \|( u ,
  \partial_1\eta )\|_{4}^2
\lesssim\| (u, \partial_1^3u,\partial_1^5\eta)
\|_{ 0}^2
+\sqrt{\mathcal{E}}\mathcal{D},\\
\label{202012251506}
&\frac{\mm{d}}{\mm{d}t}
\mathcal{E}^{\mm{cul}}_{\mm{p}}
 +c (
 \|(u,\partial_1\eta) \|_{1,3}^2+\| u \|_3^2)
\lesssim \| (u, \partial_1^3u,\partial_1^5\eta)
\|_{ 0}^2
+\sqrt{\mathcal{E}}\mathcal{D}_{\mm{p}},
\\
&
 \|\mm{curl} (\bar{\rho}\eta,
 \bar{\rho}\partial_1\eta ,
 \bar{\rho}u)\|_3^2-c \|\eta\|_4 \|  (\eta,\partial_1\eta ,u)\|_4^2\lesssim \mathcal{E}^{\mm{cul}}\lesssim
 \|  (\eta,\partial_1\eta ,u)\|_4^2, \label{202103261032} \\
&  \|\mm{curl} (\bar{\rho}u)\|_2^2+\|\mm{curl}(
 \bar{\rho}u,
 \bar{\rho} \partial_1\eta )\|_{{1},2}^2 -c \|(\eta,u)\|_4   \| (u,\partial_1\eta) \|_{1,3}^2  \nonumber \\
&\lesssim \mathcal{E}^{\mm{cul}}_{\mm{p} }\lesssim \| ( u,\partial_1(
\eta,   \partial_1\eta ,u))\|_3^2\label{2021032610321} .
\end{align}
\end{pro}
 \begin{pf} Let $1\leqslant i\leqslant 3$ and $f=\partial_1\eta$ or $\partial_1^2\eta$.
Then, for any multiindex $\beta$ satisfying $|\beta|=3-i$, it is easy to check
 $$\|[\partial_1^i\partial^{\beta}\mm{curl}_{\mathcal{A}}, {\bar{\rho}}^{-1}](\bar{\rho}f)\|_0\lesssim
 \|f\|_{\underline{i},2-i}+ \begin{cases} 0&\mbox{for }i=1;\\
\|\eta\|_{3,2 }\|f\|_1&\mbox{for }i=2,\ 3.
\end{cases}$$
In addition, for any multiindex $\alpha$ satisfying $|\alpha|\leqslant 3$, it obviously holds that
 $$\|[ \partial^{\alpha}\mm{curl}_{\mathcal{A}}, {\bar{\rho}}^{-1}](\bar{\rho}f)\|_0\lesssim
 \|f\|_{|\alpha|}.$$

Exploiting  \eqref{202012241002} and the two estimates above,   we can derive from \eqref{202005021600} for $\alpha_1 \geqslant 1$ and \eqref{202005021610} with the both cases of $\alpha_1 \geqslant 1$ and $|\alpha|\leqslant 2$ that
  \begin{align}
&\frac{\mm{d}}{\mm{d}t}
\mathcal{E}^{\mm{cul}}_{\mm{p},i}
 +c \| \mm{curl}_{\mathcal{A}}( \bar{\rho}u ,
 \bar{\rho}\partial_1\eta)\|_{i,3-i}^2 \lesssim \|\eta\|_{i+2,3-i}^2
+\sqrt{\mathcal{E}}\mathcal{D}_{\mm{p}}, \label{202012251602}\\
&\label{202005021610n0nm}
\frac{\mm{d}}{\mm{d}t}\mathcal{E}^{\mm{cul}}_{\mm{p},4}
+c\|\mm{curl}_{\mathcal{A}}\left(\bar{\rho}u\right)\|^2_2  \lesssim\| \eta\|_{2,2}^2
+\sqrt{\mathcal{E}}\mathcal{D}_{\mm{p}}
\end{align}
for some functionals $\mathcal{E}^{\mm{cul}}_{\mm{p},i}$ and $\mathcal{E}^{\mm{cul}}_{\mm{p},4}$, which are equivalent to $ \|\mm{curl}_{\mathcal{A}}( \bar{\rho}\eta ,
\bar{\rho} \partial_1 \eta, \bar{\rho}u )\|_{ i,2-i}^2  $ and  $\|\mm{curl}_{\mathcal{A}} ( \bar{\rho}u $, $
 \bar{\rho}\partial_1 \eta  )\|_2^2$, resp..

Similarly, we can easily derive from \eqref{omessetsim122n}, \eqref{202005021600} and \eqref{202005021610} that
\begin{align}
&\frac{\mm{d}}{\mm{d}t}
\tilde{\mathcal{E}}^{\mm{cul}}
+c\|\mm{curl}_{\mathcal{A}}\left( \bar{\rho}u ,
 \bar{\rho}\partial_1\eta \right)\|_3^2
\lesssim \| \eta\|_{2,3}^2
+\sqrt{\mathcal{E}}\mathcal{D},
\label{202005021632n}
\end{align}
for some functional
$\tilde{\mathcal{E}}^{\mm{cul}} $, which is equivalent to
 $\|\mm{curl}_{\mathcal{A}}( \bar{\rho}\eta, \bar{\rho}\partial_1\eta,\bar{\rho}u)\|_{3}^2 $.

 Making use of \eqref{201808181500}, \eqref{improtian1}, \eqref{improtifdsafan1n0} and Hodge-type elliptic estimate \eqref{202005021302}, we can deduce  that, for  $0\leqslant j\leqslant 3$,
\begin{align}
\|(u,\partial_1\eta)\|_{j,4-j}\lesssim &\| (u,\partial_1\eta)\|_{j,0}
+\|(\mm{div}u,\mm{curl} u,\mm{div}\partial_1\eta,\mm{curl}\partial_1 \eta)\|_{j,3-j}\nonumber \\
 \lesssim& \| (u,\partial_1\eta, \mm{curl} (\bar{\rho} u), \mm{curl}(\bar{\rho}\partial_1 \eta)   )\|_{j,2-j}\nonumber \\
 &+\|(\eta,u)\|_4\|( \eta,u)\|_{1,3} +  \begin{cases}
\|\eta\|_{3} \|\eta\|_{1,2}&\mbox{for }  j=0 ,
\\
\|\eta\|_{3} \|\eta\|_{2,2} &\mbox{for } 1\leqslant j\leqslant 3 \end{cases}\label{202008240900} \end{align}
and
\begin{align}
\|u\|_3^2
 \lesssim \| (u,\mm{curl} (\bar{\rho} u) )\|_2 +\|(\eta,u)\|_4\|( \eta,u)\|_{1,2} . \label{202008safa240900} \end{align}

Making use of \eqref{201907291432n0}, \eqref{202005031815n0}, \eqref{202008240900}, \eqref{202008safa240900}, the interpolation inequality \eqref{201807291850} and \eqref{202012241002},  we further derive from
 \eqref{202012251602}--\eqref{202005021632n} that
\begin{align}
&\frac{\mm{d}}{\mm{d}t}
\mathcal{E}^{\mm{cul}}_{\mm{p},i}
+c\|\left(  u ,
  \partial_1\eta \right)\|_{ {i},4-i}^2
\lesssim \|(u, \partial_1^2\eta)\|_{i,3-i}^2
+\sqrt{\mathcal{E}}\mathcal{D}_{\mm{p}}\mbox{ for }1\leqslant i\leqslant 3,
\nonumber \\
&\nonumber
\frac{\mm{d}}{\mm{d}t}\mathcal{E}^{\mm{cul}}_{\mm{p},4}
+c\|  u \|^2_3 \lesssim  \|u\|_0^2+\|  \eta \|_{2, 2}^2
+\sqrt{\mathcal{E}}\mathcal{D}_{\mm{p}}\\
&\frac{\mm{d}}{\mm{d}t}
\tilde{\mathcal{E}}^{\mm{cul}}
 +c \|  (  u ,
 \partial_1\eta)\|_4^2 \lesssim \|(u, \partial_1^2\eta)\|_3^2
+\sqrt{\mathcal{E}}\mathcal{D}  \nonumber .
\end{align}
Consequently, we immediately
get \eqref{202005021632} and  \eqref{202012251506} from the three estimates above by using \eqref{201807291850} and \eqref{202012241002}, where $\mathcal{E}^{\mm{cul}}:=\tilde{\mathcal{E}}^{\mm{cul}}+c \mathcal{E}^{\mm{cul}}_{\mm{p}}$ and  $\mathcal{E}^{\mm{cul}}_{\mm{p}}= \mathcal{E}^{\mm{cul}}_{\mm{p},1}+ c\sum_{j=2}^4\mathcal{E}^{\mm{cul}}_{\mm{p},j}$ for some constants $c$; moreover, exploiting \eqref{201907291432}, \eqref{201907291432n0}, \eqref{202005031815n0} and \eqref{fgestims}, we easily see that \eqref{202103261032} and \eqref{2021032610321} hold. This completes the proof.
\hfill$\Box$
\end{pf}

\subsection{Equivalent estimates}\label{subsec:equivalent}

This section is devoted to establishing the following equivalent estimates.
\begin{lem}\label{201612132242nx}
For sufficiently small $\delta$,  we have
\begin{align}
&\label{202012252005}\mathcal{E}\mbox{ is equivalent to }\|(\eta,\partial_1\eta,u)\|_4^2,\\
&\label{2017020614181721nm}
\mathcal{D}\mbox{ is equivalent to }
 \| (u,\partial_1\eta)\|_4^2,\\
&\label{202012252005nnm}\mathcal{E}_{\mm{p}},\ \mathcal{D}_{\mm{p}} \mbox{ and }\| ( u, \partial_1(u,\partial_1\eta))\|_3^2\mbox{ are equivalent},
\end{align}
where the  equivalent coefficients in \eqref{202012252005}--\eqref{202012252005nnm} are independent of $\delta$.
\end{lem}
\begin{pf}
To obtain \eqref{202012252005}--\eqref{202012252005nnm}, the key step is to establish the estimates of $\nabla q$ and $u_t$.
By  \eqref{01dsaf16asdfasf}$_2$, we have
\begin{equation*}
\begin{cases}
-\mm{div}\left(\nabla q/\bar{\rho}\right) = f^1  &\mbox{in } \Omega,\\[1mm]
\nabla q/\bar{\rho}\cdot\vec{\mathbf{n}}=f^2\cdot\vec{\mathbf{n}}  &\mbox{on } \partial\Omega,
\end{cases}
\end{equation*}
where $$
\begin{aligned}
f^1:=&\mm{div}\left(u_t+a u- {\lambda m^2}{\bar{\rho}}^{-1}\partial_1^2\eta
+ {\bar{\rho}}^{-1}\left(\nabla_{\tilde{\mathcal{A}}}q-gG_{\eta}\mathbf{e}_2\right)\right),
\\
f^2:=&-\nabla_{\tilde{\mathcal{A}}}q/\bar{\rho}.
\end{aligned} $$

Note that
\begin{align*}
\int  f^1\mm{d}y+\int_{\partial\Omega} f^2\cdot\vec{\mathbf{n}}\mathbf{d}y_1=0,
\end{align*}
thus applying the elliptic  estimate  \eqref{neumaasdfann1n} yields
\begin{align}
\|\nabla q\|_3 &\lesssim\|f^1\|_2 +\|f^2\|_3 \nonumber \\[1mm]
&\lesssim\|\left(\mm{div}_{\mathcal{A}_t}u,\mm{div}_{\tilde{\mathcal{A}}}u_t,\mm{div}_{\tilde{\mathcal{A}}}u
\right)\|_2
+\|(\partial_1^2\eta,\nabla_{\tilde{\mathcal{A}}}q,G_{\eta})\|_3.\nonumber \end{align}
Making use of \eqref{06142100x}, \eqref{06041533fwqg}, \eqref{omessetsim122} and the product estimate \eqref{fgestims}, we further derive that, for sufficiently small $\delta$,
\begin{align}
&\|\nabla q\|_3\lesssim\|(\eta_2,\partial_1^2\eta)\|_3
+\|(\eta,u)\|_4\|(u, u_t)\|_3.\label{202008241115}
\end{align}

Similarly, dividing \eqref{01dsaf16asdfasf}$_2$ by $\bar{\rho}$, and then applying $\|\cdot\|_3$ to the resulting identity, we get
\begin{align}
\|u_t\|_3 & =\|(\lambda m^2\partial_1^2\eta+gG_{\eta}\mathbf{e}_2 -\nabla_{\mathcal{A}} q-a\bar{\rho} u)/\bar{\rho}\|_3 \nonumber  \\[1mm]
&\lesssim\|(\eta_2,\partial_1^2\eta,u,\nabla q)\|_3,\nonumber
\end{align}
which, together with \eqref{202008241115}, implies, for sufficiently small $\delta$,
\begin{align}
\|(u_t,\nabla q)\|_3\lesssim\|(\eta_2,\partial_1^2\eta,u)\|_3.\label{202008asdf241115}
\end{align}
Thanks to \eqref{omessetsim122n}, \eqref{202008asdf241115}  and \eqref{202012241002}, we easily see that \eqref{202012252005}--\eqref{202012252005nnm} hold. \hfill$\Box$
\end{pf}

\subsection{\emph{A priori} stability estimate}

Now we are in a position to establish the \emph{a priori} stability estimates under the \emph{a priori} assumptions \eqref{aprpiosasfesnew} and \eqref{aprpiosesnew}.

Similarly to \eqref{202008240900}, we have,   for sufficiently small $\delta$,
\begin{align}
\|\eta\|_4
 \lesssim
  \|(\eta ,\mm{curl}(\bar{\rho}\eta))\|_3 .\label{202005021940}
\end{align}
Making use of  \eqref{20103231012}, \eqref{202103261032}, \eqref{202008240900} with $i=0$, \eqref{2017020614181721nm},  \eqref{202005021940} and  the interpolation inequality \eqref{201807291850}, we can derive from  \eqref{202008250856} and \eqref{202005021632} that, for sufficiently small $\delta$,
\begin{align}\label{201910050940}
\frac{\mm{d}}{\mm{d}t}\tilde{\mathcal{E}}+c  \mathcal{D} \lesssim  \sqrt{\mathcal{E}} \mathcal{D},
\end{align}
where  $\tilde{\mathcal{E}}:=\mathcal{E}^{\mm{cul}}+c\underline{\mathcal{E}}$ for some constant $c$ and
\begin{align}
&\label{201908081232}
\tilde{\mathcal{E}}  \mbox{ is equivalent to }\|(\eta,\partial_1\eta,u)\|^2_4  .
\end{align}
Exploiting \eqref{aprpiosesnew} and \eqref{202012252005},
 we further deduce from \eqref{201910050940} that
there exists a positive constant $\delta_1$, such that, for any $\delta\in (0,\delta_1]$,
\begin{align}\label{20191005asdf0940}
\frac{\mm{d}}{\mm{d}t}\tilde{\mathcal{E}}+c  \mathcal{D} \leqslant 0.
\end{align}
In particular,  by \eqref{201908081232}, we easily get from \eqref{20191005asdf0940} that,
 for some $c_3\geqslant 1$,
\begin{align}
\label{estemalas}
\mathcal{E}(t)+\int_0^t\mathcal{D}(\tau)\mm{d}\tau\leqslant c_3 \|(\eta^0,\partial_1\eta,u^0)\|_4^2.
\end{align}

Similarly to \eqref{20191005asdf0940},   making use of \eqref{2022103261057},  \eqref{2021032610321}, \eqref{202008240900}, \eqref{202008safa240900},  \eqref{202012252005nnm},   the interpolation inequality \eqref{201807291850} and  \eqref{202012241002},
we derive from \eqref{202008250856n0}  and \eqref{202012251506} that,  for sufficiently small $\delta$,
\begin{align}\label{201910050940nnm}
\frac{\mm{d}}{\mm{d}t}\tilde{\mathcal{E}}_{\mm{p}}+c{\mathcal{D}_{\mm{p}}} \leqslant \sqrt{\mathcal{E}} \mathcal{D}_{\mm{p}} ,
\end{align}
 where $\tilde{\mathcal{E}}_{\mm{p}}:=
\mathcal{E}^{\mm{cul}}_{\mm{p}}+c \underline{\mathcal{E}}$ for some constant $c$ and
\begin{align}
&\label{201908081232nnm}
\tilde{\mathcal{E}}_{\mm{p}},\ {\mathcal{E}_{\mm{p}}}, \ {\mathcal{D}_{\mm{p}}}\mbox{ and }
\|(u,\partial_1(\eta,\partial_1\eta,u))\|_3^2 \mbox{ are equivalent}.
\end{align}
Thus, exploiting \eqref{202012252005} and \eqref{201908081232nnm},  we further deduce from  \eqref{201910050940nnm} that
\begin{align}\nonumber
\frac{\mm{d}}{\mm{d}t}\tilde{\mathcal{E}}_{\mm{p}}+2 c_1\tilde{\mathcal{E}}_{\mm{p}}\leqslant  0,
\end{align}
which, together with \eqref{omessetsim122n} and \eqref{201908081232nnm}, further implies that
\begin{align}    \label{estemalasn0}
e^{c_1 t}(\|\eta_2(t)\|_4^2+\mathcal{E}_{\mm{p}}(t))+\int_0^te^{c_1 \tau} \mathcal{D}_{\mm{p}}(\tau)\mm{d}\tau\lesssim \|(u^0,\partial_1(\eta^0,\partial_1\eta^0,u^0))\|^2_{1,3}.
\end{align}
This completes the derivation of the \emph{a priori} stability estimates \eqref{1.200} and \eqref{1.200n0}.

\subsection{Proof  of Theorem \ref{thm2}}\label{subsec:08}

Now we state the local well-posedness result for the transformed MRT problem \eqref{01dsaf16asdfasf}.
\begin{pro}\label{202102182115}
Let $b >0$, $a\geqslant 0$ be constants and $\iota>0$ be the  constant in Lemma \ref{pro:1221}. We assume that  $\bar{\rho}$ satisfy \eqref{0102},   $(  \eta^0,u^0)\in (H^{1,4}_{\mm{s}}\cap H^4_*) \times H^4_{\mm{s}}$,
$ \|(u^0,\partial_1\eta^0)\|_4\leqslant b $  and $\mm{div}_{\mathcal{A}^0}u^0=0$, where $\mathcal{A}^0:=(\nabla\zeta^0)^{-\mm{T}}$ and $ \zeta^0 = \eta^0+y$.
Then there exist a sufficiently small  constant $\delta_2\leqslant \iota/2$, such that  if $ \eta^0 $ satisfies
\begin{align}
\| \eta^0\|_4\leqslant \delta_2, \nonumber
\end{align}
 the transformed MRT  problem \eqref{01dsaf16asdfasf} admits a unique local-in-time classical  solution
$( \eta, u,q)\in \mathfrak{C}^0(\overline{I_{T_1}},{H}^{1,3}_{\mm{s}} )\times   \mathfrak{U}_T^3  \times \mathfrak{Q}_T^3 $  for some local existence time $T>0$.  Moreover, $
 \eta $ satisfies
$$ \sup\nolimits_{t\in  {I_T}} \| \eta\|_4\leqslant 2\delta_2 \footnote{\emph{
Since $\sup_{t\in \overline{I_T}}\| \eta\|_4\leqslant  \iota$,   we have, by Lemma \ref{pro:1221},
$$\inf\nolimits_{(y,t)\in \mathbb{R}^2\times  {I_T}} \det(\nabla \eta+I)\geqslant 1/4.$$}}$$ and \begin{align}
\label{202012212151}
  \sup\nolimits_{t\in  {I_T}}\| (\eta,    \partial_1 \eta ,  u)\|_4 \leqslant c_4 \sqrt{I^0}
\end{align}
 for some positive constant $c_4\geqslant 1 $. It should be noted that $\delta_2$ and $T$ may depend on $g$, $a$,  $\lambda$, $m$, $\bar{\rho}$ and $\Omega$; moreover $T$ further depends on $b$.
\end{pro}
\begin{pf} The proof of Proposition \ref{202102182115} will be provided in Section \ref{202102241211}.
\hfill $\Box$
\end{pf}

Thanks to the \emph{priori} estimate \eqref{estemalas}  and Proposition  \ref{202102182115}, we can easily establish
Theorem \ref{thm2}. Next we briefly give the proof.

Let  $(\eta^0,u^0)$ satisfies the assumptions in Theorem \ref{thm2}, and
\begin{align}
I^0\leqslant\delta^2,\mbox{ where }
\delta:=\min\left\{ \delta_1, \delta_2  \right\}/ \sqrt{{c}_3}c_4^2,\nonumber
\end{align}
  where $c_3$ and $c_4$ are the constants in \eqref{estemalas} and \eqref{202012212151}.
By Proposition  \ref{202102182115} and   Lemma \ref{pro:1221}, there exists a unique local solution $( \eta, u,q)$ to the  transformed MRT  problem \eqref{01dsaf16asdfasf} with a maximal existence time $T^{\max}$, which satisfies that
\begin{itemize}
  \item for any $\tau\in  I_{T^{\max}}$,
the solution $(\eta,u,q)$ belongs to ${\mathfrak{H}}^{1,4}_{1,*,\tau}\times \mathfrak{U}_\tau^4  \times   \mathfrak{Q}_\tau^4$ and
$$ \sup\nolimits_{t\in  {I_\tau}} \| \eta\|_4\leqslant 2\delta_2;$$
  \item $\limsup_{t\to T^{\max} }\|  \eta( t)\|_4 > \delta_2$ or $\limsup_{t\to T^{\max} }\|(u,\partial_1\eta)( t)\|_4=\infty$, if $T^{\max}<\infty$.
\end{itemize}

Let
\begin{equation}
\nonumber
T^{*}:=\sup\left\{ \tau \in I_{T^{\max}}~\left|~ \|(\eta,\partial_1\eta,u)(t)\|_4\leqslant \sqrt{{c}_3}c_4^2\delta  \mbox{ for any }t\leqslant \tau\right.\right\}.
\end{equation}
We easily see that the definition of $T^*>0$ makes sense and $T^*>0$. Thus, to obtain the existence of global solutions, next it suffices to verify $T^*=\infty$. Next we shall prove this fact by contradiction.

We assume that $T^*<\infty$. Then, for any given $T^{**}\in I_{T^*}$,
\begin{equation}
\label{201911262sadf202}
\sup\nolimits_{\overline{I_{T^{**}}}}\|(\eta,\partial_1\eta, u)(t) \|_4   \leqslant \sqrt{{c}_3}c_4^2\delta \leqslant \delta_1.
\end{equation}
Thanks to \eqref{201911262sadf202},
  we can follow the argument of \eqref{201910050940}  by further using a  regularity method as in the derivation of \eqref{202104112141} to  verify that
 \begin{align}
 &\label{201910051asdfa012}
\frac{\mm{d}}{\mm{d}t}\tilde{\mathcal{E}}+c\mathcal{D}\leqslant 0 \mbox{ for a.e. }t\in I_{T^*},\mbox{ where } \tilde{\mathcal{E}} \in W^{1,\infty}(I_{T^*}).
\end{align}
Referring to \eqref{estemalas}, we can further derive from \eqref{202012212151} and  \eqref{201910051asdfa012} that \begin{align}
{\|\mathcal{E}(t)\|_{L^\infty(I_{T^*})} + \int_0^t\mathcal{D}(\tau)\mm{d}\tau \leqslant {{c}_3}\lim_{t\to 0}\sup\nolimits_{\tau\in I_t}\|(\eta,\partial_1\eta,u)(\tau)\|_4^2 \leqslant {{c}_3}c_4^2\delta^2.}  \label{202010326asdfas1534}
\end{align}
By the continuity of $(\eta,u)$ and the fact
  \begin{align}\sup\nolimits_{t\in \overline{I_\tau}} \|f\|_0=   \|f\|_{L^\infty_\tau L^2}\mbox{ for any }f\in C^0_{B,\mm{weak}}(\overline{I_\tau}, L^2)\mbox{ with }\tau>0,\nonumber
\end{align}
we further derive from \eqref{202010326asdfas1534} that
\begin{align}
{ \sup\nolimits_{\overline{I_{T^*}}}\|(\eta,\partial_1\eta,u)(t)\|_4  \leqslant \sqrt{{c}_3}c_4\delta.}  \label{2020103261534}
\end{align}

We take $(\eta(T^{**}),u(T^{**}))$ as a initial data. Noting that, by \eqref{2020103261534} and the definition of $\delta$,
 $$ \|( \eta, \partial_1\eta,u)(T^{**})\|_4  \leqslant  {B} := \sqrt{{c}_3}c_4\delta\mbox{ and }
\| \eta(T^{**})\|_4\leqslant \delta_2
 ,$$ then,
by Proposition   \ref{202102182115}, there exists a unique local-in-time classical solution, denoted by $( \eta^*, u^*$, $q^*)$, to the transformed MRT  problem \eqref{01dsaf16asdfasf} with $(\eta(T^{**}),u(T^{**}))$ in place of $(\eta^0,u^0)$; moreover
\begin{align}
 \sup\nolimits_{t\in [T^{**},T]}\| (\eta^*, \partial_1 \eta^*, u^*   )\|_4 \leqslant  c_4B\leqslant  \sqrt{{c}_3}c_4^2\delta \mbox{ and } \sup\nolimits_{t\in [T^{**},T]}\| \eta^*\|_4\leqslant 2\delta_2  , \nonumber
\end{align}
where the local existence time $T>0$ depends on $ b$, $g$, $a$,  $\lambda$, $m$, $\bar{\rho}$ and $\Omega$.

In view of the existence result of $( \eta^*, u^*, q^*)$, the uniqueness conclusion in Proposition  \ref{202102182115}  and the fact that $T^{\max}$ denotes the maximal existence time, we immediately see that $T^{\max}>T^*+T/2$ and
$ \sup_{t\in  {[0,T^*+T/2]}}\| (\eta, \partial_1 \eta, u)\|_4 \leqslant  \sqrt{{c}_3}c_4^2\delta$.
This contradicts with the definition of $T^*$. Hence $T^*=\infty$ and thus $T^{\max}=\infty$.  This completes the proof of the existence of global solutions. The uniqueness of global solution is obvious due to the uniqueness result of local solutions in Proposition \ref{202102182115} and the fact $\sup_{t\geqslant 0}\| \eta\|_4\leqslant 2 \delta_2$. To complete the proof of Theorem \ref{thm2}, we shall verify that the solution $(\eta,u,q)$ satisfies the properties \eqref{1.200xyx}--\eqref{1.200xx}.

Referring to \eqref{201910051asdfa012} and \eqref{202010326asdfas1534}, we see that
 the global solution $(\eta,u)$ enjoys \eqref{1.200xyx} and \eqref{1.200}. Similarly,   we also verify that the global solution also satisfies  \eqref{1.200n0} by  referring to the derivation of \eqref{estemalasn0} and \eqref{202010326asdfas1534}.

Finally, we shall derive \eqref{1.200xx}.
By \eqref{1.200n0}, we easily see that
\begin{align}
\label{202101242112}
&\partial_1\eta(t)\to 0\mbox{ in }H^2\mbox{ as }t\to \infty
\end{align}
and
\begin{align}
\label{20safda2101242112}
\left\|\int_0^t u\mm{d}\tau\right\|_3\lesssim\int_0^t \|u  \|_3\mm{d}\tau\lesssim \sqrt{I_0}  \mbox{ for any }t>0\end{align}
Thanks to \eqref{20safda2101242112}, there exist a subsequence $\{t_n\}_{n=1}^{\infty}$
and some function $\eta^\infty_1\in H^3$ such that
$$ \int_0^{t_n} u_1\mm{d}\tau \to\eta^\infty_1-\eta^0_1\mbox{ weakly in }H^3.$$
Exploiting $\eqref{01dsaf16asdfasf}_1$, \eqref{1.200n0}  and weakly lower semi-continuity, we have
$$
\begin{aligned}
\|\eta_1(t)-\eta^\infty_1\|_3\leqslant &\liminf_{t_n\to \infty} \left\|\int_t^{t_n } u_1\mm{d}\tau\right\|_3
\lesssim \sqrt{I^0} \liminf_{t_n \to \infty}\int_t^{t_n} e^{-c_1 \tau}\mm{d}\tau
 \lesssim \sqrt{I^0}e^{-c_1 t} ,
\end{aligned}$$
which, together with \eqref{202101242112}, yields that \eqref{1.200xx} holds and $\eta_1^\infty $ only depends on ${y_2}$. This completes the proof of Theorem  \ref{thm2}.

\section{Proof of Corollary \ref{cor1}} \label{202103171455}
This section is devoted to the proof of Corollary \ref{cor1}.
Let $(\eta,u,q)$ be the classical solution  constructed by
Theorem  \ref{thm2} with initial data $(\eta^0,u^0)$ further satisfying  the odevity conditions \eqref{202005011004} and \eqref{202005011005}.
Next we first verify the solution preserves the odevity conditions in \eqref{202005012204}.

Let $\psi=(-\eta_1,\eta_2)(-y_1,y_2,t )$, $w=( -u_1,u_2)(-y_1,y_2,t)$ and $p=( q_1, q_2)(-y_1,y_2,t)$.
Since the classical solution of the transformed MRT problem \eqref{01dsaf16asdfasf}  is unique, thus, to get the  relation \eqref{202005012204}, it obviously suffices to verify that $(\psi,w,p)$ is also the classical solution of \eqref{01dsaf16asdfasf}. It is obvious that $(\psi,w)$ satisfies \eqref{01dsaf16asdfasf}$_1$. Next we shall verify that $(\psi,w)$ also satisfies \eqref{01dsaf16asdfasf}$_2$ and \eqref{01dsaf16asdfasf}$_3$.

Defining \begin{align} {\mathcal{A}}=\left(\begin{array}{cc}
            \partial_2\eta_2+1 &  - \partial_1\eta_2\\
                - \partial_2\eta_1&   \partial_1\eta_1+1
                 \end{array}\right)\mbox{ and } \mathcal{B}=\left(\begin{array}{cc}
            \partial_2\psi_2+1 &  - \partial_1\psi_2\\
                - \partial_2\psi_1&   \partial_1\psi_1+1
                 \end{array}\right), \nonumber
\end{align}
then
\begin{align}
& (\mathcal{B}_{11} ,\mathcal{B}_{22})= ({\mathcal{A}}_{11} ,{\mathcal{A}}_{22})|_{y_1=-y_1 }\mbox{ and } (\mathcal{B}_{12}, \mathcal{B}_{21})= -({\mathcal{A}}_{12} ,{\mathcal{A}}_{21})|_{y_1=-y_1}.
\label{202002122032}
\end{align}

By \eqref{202002122032},
\begin{align}\nabla_{{\mathcal{B}}}p   = (-{\mathcal{A}}_{1i} , {\mathcal{A}}_{2i})^{\mm{T}}  \partial_i q |_{y_1=-y_1}
\label{20200212203safdsaf2}
\end{align}
 and
\begin{align}\mm{div}_{\mathcal{B}}
w = \mm{div}_{ {\mathcal{A}}}u |_{y_1=-y_1 }=0.
\label{20204261426}
\end{align}
By  \eqref{20204261426}, we see that $(\psi,w)$ satisfies  \eqref{01dsaf16asdfasf}$_3$ as $(\eta,u)$.

Thanks to  \eqref{20200212203safdsaf2}, we have
$$
\begin{aligned}
 & \bar{\rho}\partial_tw_i +\nabla_{ {\mathcal{B}}}p + a\bar{\rho} w_i =
      \begin{cases}
  - (\bar{\rho}\partial_t u_1+ {\mathcal{A}}_{1i}\partial_i q+\bar{\rho} u_1 )|_{y_1=-y_1 }  & \hbox{for }i=1; \\
   (\bar{\rho}\partial_t u_2+{\mathcal{A}}_{2i}\partial_i q+ a\bar{\rho} u_2
 )|_{y_1=-y_1 }  & \hbox{for }i=2,
      \end{cases}
 \\
& \partial_1^2 (\psi_1,\psi_2)=  \partial_1^2 (-\eta_1,\eta_2)|_{y_1=-y_1 }\mbox{ and }
g{G}_\psi\mathbf{e}_2=gG_{\eta}\mathbf{e}_2\big|_{y_1=-y_1} ,
\end{aligned}$$
where ${G}_\psi$ is defined in \eqref{202009130836} with $\phi_2$ in place of $\eta_2$.
Hence we see that $(\psi,w)$ also satisfies \eqref{01dsaf16asdfasf}$_2$ by the three identity above.  This completes the proof of the verification of preserving odevity of solutions.

Thanks to \eqref{202005012204} and \eqref{poinsafdcaasdsfaressadf1}, we easily get, for $0\leqslant i\leqslant 3$,
$$ \|\partial_2^i\eta_1\|_0\lesssim \|\partial_1\partial_2^i\eta_1 \|_0,$$
which, together the estimate \eqref{omessetsim122n} satisfied by $\eta$,  implies that
$$\mathcal{E}\mbox{ is equivalent to } {\mathcal{D}}. $$
Thus we can further derive \eqref{1.200n} from \eqref{1.200xyx} and \eqref{1.200}. This completes the proof of Corollary \ref{cor1}.

\section{ Proof of Theorem \ref{thm1}}\label{sec:instable}
This section is devoted to the proof of instability of transformed MRT problem in Theorem \ref{thm1}.
We will complete the proof by five subsections. In what follows, the fixed positive constant $c_i^I$ for $i\geqslant 1$ may depend on $g$, $a$,  $\lambda$, $m$,  $\bar{\rho}$  and $\Omega$.
\subsection{Linear instability}\label{sunsec:growing mode}
To begin with, we exploit modified variational method of ODE as in \cite{GYTI2,JFJWGCOSdd}  to prove the existence of unstable solutions of the following linearized MRT problem under the instability condition $|m|\in(0,m_{\mm{C}})$:
\begin{equation}\label{01dsaf16asdfasf0101}
                              \begin{cases}
\eta_t=u   ,\\[1mm]
\bar{\rho}u_t+\nabla q+a\bar{\rho} u=\lambda m^2\partial_1^2\eta+g\bar{\rho}'\eta_2\mathbf{e}_2   ,\\[1mm]
\div u=0 , \\[1mm]
(\eta,u)|_{\partial\Omega}\cdot\vec{\mathbf{n}}=0.
\end{cases}
\end{equation}
\begin{pro}\label{pro:08252100}
Let $a\geqslant 0$ and $\bar{\rho}$ satisfy \eqref{0102} and \eqref{0102n}. If  $|m|\in[0,m_{\mm{C}})$,
 then the zero solution
 is unstable to the linearized MRT problem. That is, there is an unstable solution
$(\eta, u,  q):=e^{\Upsilon t}(w/\Upsilon,w, \beta )$
 to  the above problem, where
 \begin{equation}\nonumber
 (w, \beta )\in H^5_\sigma\times\underline{H}^4
 \end{equation}
 solves  the boundary-value problem  \begin{equation}\label{01dsaf16asdfasf0202}
                              \begin{cases}
\Upsilon^2\bar{\rho}w+ \Upsilon\nabla  \beta
+a\Upsilon \bar{\rho} w=  m^2\partial_1^2w+g\bar{\rho}'w_2\mathbf{e}_2
&\mbox{in } \Omega,\\[1mm]
\div w=0  &\mbox{in } \Omega, \\[1mm]
w\cdot\vec{\mathbf{n}}=0  &\mbox{on } \partial\Omega.
\end{cases}
\end{equation}
  with some
growth rate $\Upsilon\in  (2\Lambda/3,\Lambda]$, where $\Lambda$ satisfies
\begin{equation} \label{201907092040}
E(v)\leqslant ({\Lambda^2}+{a\Lambda })\|\sqrt{\bar{\rho}}v\|_0^2\mbox{ for any }v\in H_{\sigma}^{1} . \end{equation}
In addition,
\begin{align}
& \|  w_i\|_{0}\|\partial_1w_i\|_{0}\|\partial_2w_i\|_0\neq 0
\label{201602081445MH}\mbox{ for }i=1,\ 2. \end{align}
\end{pro}
\begin{pf}
Similarly to \cite{GYTI2,JFJWGCOSdd}, next we divide the proof of Proposition \ref{pro:08252100} into five steps.

(1) Let
$$\tilde{E}(\psi, \xi)=\int_0^h\left(g\bar{\rho}'|\psi|^2-\lambda m^2
\left(\xi^{2}|\psi|^{2}+\left|\psi'\right|^{2}\right)\right)\mm{d} y_{2}$$
and
\begin{equation}\nonumber \mathbb{F}:=\{\xi\in\mathbb{Z}/\{0\}~|~
\tilde{E}(\psi, \xi)>0\mbox{ for some }\psi\in H_0^1(0,h)\}.
\end{equation}
\emph{Next we prove the first assertion that the set of instability frequencies  $\mathbb{F}$ is not empty. }

Recalling the definition of $m_{\mm{C}}$ and   the condition $|m|\in[ 0,m_{\mm{C}})$, we see that there exists a function
$\omega:=(\omega_1,\omega_2)^{\mm{T}}\in H_{\sigma}^1$, such that
\begin{align}\label{202009010930}
\int g\bar{\rho}'\omega_2^2\mm{d}y-\lambda m^2\|\partial_1\omega\|_0^2>0 .
\end{align}
Let $\hat{\omega}(\xi,y_2)$ be the Fourier coefficient of $\omega(y_1,y_2)$ for fixed $y_2$, i.e.,
$$\hat{\omega}(\xi,y_2)=\int_0^{2\pi}\omega(y_1,y_2)\mathrm{e}^{-i\xi y_1}\mm{d}y_1.$$
We define the functions $\varphi$ and $\psi$ by the following relations
$$\hat{\omega}_1(\xi,y_2)= i\varphi (\xi,y_2)\mbox{ and }\hat{\omega}_2(\xi,y_2)=-\psi (\xi,y_2),$$
where $(\hat{\omega}_1, \hat{\omega}_2)^{\mm{T}}=\hat{\omega}$ and $\xi\in\mathbb{Z}$. Obviously,
\begin{align}
\label{2020101291058}
\widehat{\partial_1\omega_1}=\xi\varphi\mbox{ and }\widehat{\partial_1\omega}_2=-i\xi\psi.
\end{align}
Since $\omega_2|_{\partial\Omega}=0$, then $\psi\in H^1_0(0,h)$. By $\mm{div}\omega=0$, we have
$$\xi\varphi+\psi'=0,$$
which, together with the fact $\psi|_{y_2=0,\ h}=0$, implies that
\begin{align}
\label{2021012910578}
\psi(0,y_2)=0\mbox{ for }\xi=0.
\end{align}

Exploiting \eqref{2020101291058}, \eqref{2021012910578}, and Fubini and Parseval Theorems, we have
\begin{align*}
E(\omega )=\frac{1}{2\pi}\int_0^h
\sum_{\xi\in\mathbb{Z}/\{0\}}\left(g\bar{\rho}'|\psi|^2-  \lambda m^2\left(|\xi\psi|^{2}+\left|\psi'\right|^{2}\right)\right)
\mm{d} y_{2}.
\end{align*}
The above   identity together with \eqref{202009010930} imply that there is a $\xi\in\mathbb{Z}/\{0\}$ and $\psi\in H^1_0(0,h)$, such that \begin{align}
\nonumber
\tilde{E}(\psi, \xi)>0.
\end{align}
Hence the set of instability frequencies $\mathbb{F}$ is not empty.

(2)
We define that
 \begin{align}\nonumber
 &\mathfrak{J}(\psi, \xi):=\int_{0}^{h}\bar{\rho}
\left(|\xi \psi|^{2}+\left|\psi'\right|^{2}\right)\mm{d} y_{2},\ \mathcal{H}:=\left\{\psi \in H_{0}^{1}(0,h)~|~ \mathfrak{J}(\psi,\xi)=1 \right\},\nonumber \\
& \Xi(\psi,\xi,s) := \xi^{2}\tilde{E}(\psi, \xi)-as\mathfrak{J}(\psi, \xi)  .\nonumber
\end{align} For any given $\xi\in \mathbb{F}$, we define that
$$ \mathfrak{C}_\xi:= \sup\{s\in \mathbb{R}~|~\Xi(\psi,\xi,s)> 0\mbox{ for some }\psi\in H_0^1(0,h) \}.$$
Recalling the definition of $\mathbb{F}$, we easily see that
\begin{align}
 0<\mathfrak{C}_\xi,\  \mathfrak{C}_\xi=\infty \mbox{ for }a=0 \mbox{ and } \mathfrak{C}_\xi<\infty \mbox{ for }a> 0.  \nonumber
\end{align}

\emph{
 Next we prove the  second assertion that, for any given $(\xi,s)\in \mathbb{F}\times [0,\mathfrak{C}_\xi)$, there exist  $\alpha>0$ and a classical solution $\psi_0\in H_0^1(0,h)\cap H^4(0,h)$, which satisfies the following  modified boundary-value problem:
\begin{equation}\label{202008asdfsa251540}
\begin{cases}
(\alpha+as)\left(\bar{\rho}\xi^2\psi_0-(\bar{\rho}\psi_0')'\right)=
g\xi^2\bar{\rho}'\psi_0-\lambda m^2\xi^2\left(\xi^2\psi_0- \psi_0''\right)\mbox{ in }(0,h),&
\\
\psi_0(0)=\psi_0(h)=0.
\end{cases}
\end{equation}
Moreover, $\sup_{\psi\in\mathcal{H}}\Xi(\psi,\xi,s)=\Xi(\psi_0,\xi,s) $.}

To being with, we shall consider the following variational problem: \begin{align}
\alpha:=\sup_{\psi\in\mathcal{H}}\Xi(\psi,\xi,s),
\label{2021012917205}
\end{align}
where $\xi\in \mathbb{F}$ and $s\in \mathbb{R}_0^+$.
Since $\alpha$ depends on $s$, sometimes we denote $\alpha$ by $\alpha(s)$.
Obviously,
\begin{align}
\label{2019070safa12110}
& \alpha(s_2) -\alpha(s_1)=a(s_1-s_2)\mbox{ for any } s_1,\ s_2\geqslant 0,\\ \label{201907saf012110}
& g\|\bar{\rho}'/\bar{\rho}\|_{L^\infty}\geqslant \alpha >0\mbox{ for any }0\leqslant s<  \mathfrak{C}_\xi ,\\
\label{201907012110}
&\alpha(\mathfrak{C}_\xi)=0\mbox{ if } \mathfrak{C}_\xi\in \mathbb{R}^+, \mbox{i.e., the case $a>0$}.
\end{align}

Since $0\leqslant \alpha<\infty$, there exists a maximizing sequence $\{\psi_n\}_{n=1}^\infty\subset\mathcal{H}$
 such that $\psi_n\to \psi_0$ weakly in $H^1_0(0,h)$ and strongly in $L^2(0,h)$. By the convergence results of $\psi_n$ and the weakly lower semi-continuity, we easily get
  \begin{align}
\mathfrak{J}(\psi_0, \xi) \leqslant  \lim_{n\to\infty}\|\sqrt{\bar{\rho}}\psi_n\|^2_{L^2(0,h)}+ \liminf_{n\to\infty}\|\psi'_n\|^{2}_{L^2(0,h)} \leqslant\liminf_{n\to\infty}\mathfrak{J}(\psi_n, \xi)=  1
\label{202101291659}
\end{align}
and
 \begin{align}
 \alpha= \lim_{n\to\infty}\Xi(\psi_n,\xi,s)  \leqslant & g\lim_{n\to\infty}\int_0^h\bar{\rho}'|\psi_n|^2\mm{d} y_{2} - \lambda m^2 \lim_{n\to\infty}\int_0^h
 \xi^{2}|\psi_n|^{2}\mm{d} y_{2}\nonumber \\
&- \lambda m^2 \liminf_{n\to\infty}\int_0^h  |\psi'_n |^{2} y_{2}-as\leqslant \Xi(\psi_0 , \xi).
\label{202101291659z}
\end{align}

\emph{From now on, we consider the case $s\in [0,\mathfrak{C}_\xi)  $.}
Now we prove $\psi_0\neq0$ by contradiction. We assume that $\psi_0 =0$, then there exists a subsequence $\{\psi_{n}\}_{n=1}^\infty$ (still denoted by $\psi_n$ for simplicity), such that
$ \|\psi_{n}'\|_{L^2(0,h)}\to b >0$ and $\|\sqrt{\bar{\rho}}\psi_{n}'\|_{L^2(0,h)}\to 1$.
Thus we immediately see that
\begin{align}
\Xi(\psi_n,\xi,s)\to \alpha=-(\lambda bm^2\xi^2+as)
                                                       \leqslant  0  \hbox{ for } s\in [0,\mathfrak{C}_\xi) \mbox{ and } a\geqslant 0  ,
 \nonumber
\end{align}
which contradicts with \eqref{201907saf012110}. Thus we get $\psi_0\neq0$.

Since $\psi_0\neq0$, we see from \eqref{202101291659} that
 \begin{align}
 0<\mathfrak{J}(\psi_0, \xi)\leqslant 1. \nonumber
 \end{align}
Thus we derive from \eqref{202101291659z} and the relation above that
$$\alpha \leqslant  \alpha \mathfrak{J}^{-1}(\psi_0, \xi)\leqslant \Xi(\psi_0, \xi)\mathfrak{J}^{-1}(\psi_0, \xi)  \leqslant \alpha,$$
which yields $\mathfrak{J}^{-1} (\psi_0, \xi) =1$.  This means that $\psi_0$
is an achievable point of the variational problem \eqref{2021012917205}. In particular, we have
\begin{align}
\alpha=\frac{\Xi(\psi_0,\xi,s)}{ \mathfrak{J}(\psi_0, \xi)}\geqslant  \frac{\Xi(\psi,\xi,s)}{ \mathfrak{J}(\psi, \xi)} \mbox{ for any }\psi\in H^1_0(0,h).
\label{2022101291720}
\end{align}

For any given $\phi\in H^1_0(0,h)$, let
$$\mathfrak{F}(\tau) =  \Xi((\psi_0+\tau \phi),\xi,s)-(\alpha+as)\mathfrak{J}((\psi_0+\tau\phi), \xi)  \mbox{ for any }\tau\in \mathbb{R}.$$
 Then $\mathfrak{F}(\tau)\in C^\infty(\mathbb{R})$. By \eqref{2022101291720}, $\mathfrak{F}(0)=0$ and $\mathfrak{F}(\tau)\leqslant 0 $. Thus we have $\mathfrak{F}'(0)=0$. In particular, we have
\begin{equation}\nonumber
(\alpha+as)\int_0^h\bar{\rho}\left(\xi^2 \psi_0 \phi+ \psi_0' \phi'\right)\mm{d}y_2=
\int_0^h \xi^2(g\bar{\rho}'\psi_0 \phi -\lambda m^2 (\xi^2\psi_0 \phi+ \psi'_0 \phi'))\mm{d}y_2.
\end{equation}
Noting that $ \bar{\rho} (\alpha+as)+\lambda m^2\xi^2\neq 0  $ for any $y_2\in (0,h)$, thus
we equivalently rewrite the above identity as follows:
\begin{align}\nonumber
&\int_0^h(\bar{\rho} (\alpha+as)+\lambda m^2\xi^2  )^{-1}( \xi^2( g\bar{\rho}'-\lambda m^2\xi^2 -(\alpha+as)\bar{\rho})\psi_0\\
&\qquad + (\alpha+as  )\bar{\rho}' \psi'_0 )\chi\mm{d}y_2  =
\int_0^h \psi'_0 \chi' \mm{d}y_2,  \label{20211918542}
\end{align}
where $\chi:=(\bar{\rho} (\alpha+as)+\lambda m^2\xi^2  )  \phi$.
It is easy to see from \eqref{20211918542} that  $\psi_0\in H^4(0,h)$ and  $\psi_0$ is just the classical solution of the modified boundary-value problem \eqref{202008asdfsa251540}. This completes the proof of the second assertion.

(3) \emph{Next we further prove the third assertion that, for any given $\xi\in \mathbb{F}$, there exist $\gamma>0$  and $\psi\in  H_0^1(0,h)\cap H^4(0,h)$, such that
\begin{equation}\label{202008251540}
\begin{cases}
\gamma^2\left(\bar{\rho}\xi^2\psi-(\bar{\rho}\psi')'\right)=
g\xi^2\bar{\rho}'\psi-\lambda m^2\xi^2\left(\xi^2\psi-\psi''\right)
-a\gamma\left(\bar{\rho}\xi^2\psi- (\bar{\rho}\psi')'\right)\mbox{ in }(0,h),&\\
\psi(0)=\psi(h)=0.
\end{cases}
\end{equation}
Moreover,} \begin{align}
\label{20asfas2101292338}
\sup_{\chi\in\mathcal{H}}\Xi(\chi,\xi,\gamma)= \Xi(\psi ,\xi,\gamma)=\gamma>0.
\end{align}
The above assertion obviously holds for $a=0$ by the second assertion, thus it suffices to consider the case $a> 0$.

Thanks to \eqref{2019070safa12110}--\eqref{201907012110}, it is easy to see  that there exists a fixed point $\gamma$ satisfying \eqref{20asfas2101292338}.
Moreover, by the second assertion, there exists $\psi\in  H_0^1(0,h)\cap H^5(0,h)$ satisfying \eqref{202008251540} with $\gamma$. This completes the proof of the third assertion.

(4) \emph{Now we are in the position to the proof of existence of a solution $(w, \beta)$ to the boundary-value problem \eqref{01dsaf16asdfasf0202} with some $\Upsilon>0$}

Let $\psi$ and $\gamma$ be constructed in the third assertion for any given $\xi\in \mathbb{F}$. From now on, we denote $(\psi,\gamma)$ by $(\psi_\xi,\gamma_\xi)$ to emphasize the dependence of $\xi$. We define $\Lambda :=\sup _{\xi\in\mathbb{F}}\gamma_\xi$. It is easy to see that
\begin{align}
-\xi\in \mathbb{F},\ (\psi_\xi,\gamma_\xi)=(\psi_{-\xi},\gamma_{-\xi})\mbox{ for any }\xi\in \mathbb{F} \nonumber
\mbox{ and }
0<\Lambda <\infty.
\end{align}

Denoting that
$$ \varphi_\xi:=-\xi^{-1}\psi'_\xi\mbox{ and }  \vartheta_\xi:=(\gamma^2\bar{\rho}\varphi_\xi +a\gamma\bar{\rho}\varphi_\xi+\lambda m^2\xi^2\varphi_\xi)  /\gamma_\xi\xi,$$
then it is easy to check that
 \begin{equation}\label{201907011445}
\begin{cases}
\gamma^2_\xi\bar{\rho}\varphi_\xi-\gamma_\xi\xi\vartheta_\xi+a\gamma_\xi\bar{\rho}\varphi_\xi=
-\lambda m^2\xi^2\varphi_\xi &\mbox{in } (0,h),\\[1mm]
\gamma^2_\xi\bar{\rho} \psi_\xi+\gamma_\xi\vartheta'_\xi+a\gamma_\xi\bar{\rho}\psi_\xi=
-\lambda m^2\xi^2\psi_\xi+g\bar{\rho}'\psi_\xi &\mbox{in } (0,h),\\[1mm]
\xi\varphi_\xi+\psi'_\xi=0 &\mbox{in } (0,h),\\[1mm]
\psi_\xi(0)=\psi_\xi(h)=0. &
\end{cases}
\end{equation}

Let
\begin{equation}\nonumber
w_{1}(y) =-i  \varphi_\xi\left(y_{2}\right) ({e}^{i y_1 \xi} -{e}^{-i y_1 \xi}) \mbox{ and }
(w_{2} , {\theta} )(y) =( \psi_\xi\left(y_{2}\right),\vartheta_\xi\left(y_{2}\right) )({e}^{i y_1 \xi}+{e}^{-i y_1 \xi}),
\end{equation}
then $w_1$, $w_2$ and $ {\theta}$ are real-value functions. Recalling \eqref{20asfas2101292338} and \eqref{201907011445}$_4$, we easily see that $\psi_\xi$, $\varphi_\xi$, $\psi_\xi'$ and  $\varphi_\xi'$ are non-zero function. Thus we have $\|w_i\|_0\|\partial_1 w_i\|_0\|\partial_2 w_i\|_0\neq0$ for $i=1$ and $2$.
 Let $w=(w_1,w_2)^{\mm{T}}$ and $ \beta := {\theta}-(\theta)_{\Omega}$,
 then  $(w, \beta )\in H^5_\sigma\times\underline{H}^4$. By  \eqref{201907011445}, we easily see that $(w, \beta )$ satisfies \eqref{01dsaf16asdfasf0202} with $\gamma_\xi$ in place of $\Upsilon$. By the definition of $\Lambda$, there exists an $\xi_0\in \mathbb{F}$ such that $\gamma_{\xi_0} \in(2\Lambda/3,\Lambda]$. Now we take $\xi= \xi_0 $ and $\Upsilon=\gamma_{\xi_0}$.  Consequently, we see that the solution $(w, \beta ,\Upsilon)$ satisfies \eqref{201602081445MH} and the boundary-value problem \eqref{01dsaf16asdfasf0202} with $\Upsilon\in (2\Lambda/3,\Lambda]$.
 In addition, it is easy to see that
 $(\eta, u,  q):=e^{\Upsilon t}(w/\Upsilon,w, \beta )$
 is a solution of \eqref{01dsaf16asdfasf0101}.

 (5) \emph{To complete the proof of Proposition \ref{01dsaf16asdfasf0101}, we shall verifies \eqref{201907092040}.}

Let $\xi\in \mathbb{Z}$, and the two real-valued functions $\tilde{\varphi}$, $\tilde{\psi}\in H_0^1(0,h) $ satisfy
\begin{equation}\label{201907092118}
  \xi\tilde{\varphi}+\tilde{\psi}'=0.
\end{equation}
In particular, we  have
\begin{equation}\label{201907asdfa092118}
\tilde{\psi}=  \tilde{\psi}'=0\mbox{ if }\xi=0.
\end{equation}
Thanks to \eqref{20asfas2101292338} and the fact $\gamma\leqslant \Lambda$, we see that
\begin{equation}\label{201907012126}
\tilde{E}(\tilde{\psi}, \xi)\leqslant {\xi^{-2}}
\left(\Lambda^2 +a\Lambda \right)\mathfrak{J}(\tilde{\psi}, \xi), \mbox{ if }\xi\in\mathbb{F}.
\end{equation}
In addition, using the definition of $\mathbb{F}$ and \eqref{201907asdfa092118},
we have
\begin{equation}\label{201907012126n}
\tilde{E}(\tilde{\psi}, \xi)\leqslant0\mbox{ if }\xi\in\mathbb{Z}/\mathbb{F}.
\end{equation}
Exploiting  \eqref{201907012126} and \eqref{201907012126n}, we obtain that, for any $\xi \in \mathbb{Z}$, and for any $\tilde{\varphi}$, $\tilde{\psi}\in H_0^1(0,h) $ satisfying \eqref{201907092118},
\begin{align}\label{201907092140}
\tilde{E}(\tilde{\psi}, \xi) \leqslant
(\Lambda^2+a\Lambda)\int_0^h\bar{\rho}\left(|\tilde{\varphi}|^{2}+|\tilde{\psi}|^{2}\right) \mathrm{d} y_{2}.
\end{align}

Let $v\in H_\sigma^1$, $\hat{v}(\xi,y_2)$ be the  Fourier coefficient of $v(y_1,y_2)$ for fixed $y_2$,  $\varphi\left(\xi, y_{2}\right)=\mathrm{i} \hat{v}_{1}\left(\xi, y_{2}\right)$ and $\psi\left(\xi, y_{2}\right)=\hat{v}_{2}\left(\xi, y_{2}\right)$. Then $(\varphi,\psi)$ satisfies $\xi{\varphi}+{\psi}'=0$.
Consequently, making use of \eqref{201907092140}, and the Fubini's and Parseval's theorems, we have
\begin{align}
E(v)=&\frac{1}{2 \pi}\sum_{\xi \in \mathbb{Z}}
\tilde{E}({\psi}, \xi)\leqslant
\frac{1}{2 \pi}(\Lambda^2+a\Lambda)
\sum_{\xi \in  \mathbb{Z}}
\int_0^h\bar{\rho}\left(|{\varphi}|^{2}+|{\psi}|^{2}\right) \mathrm{d} y_{2}  \nonumber \\[1mm]
\leqslant&(\Lambda^2+a\Lambda)\|\sqrt{\bar{\rho}}v\|_0^2, \nonumber
\end{align}
which yields \eqref{201907092040}. This completes the proof.
\hfill$\Box$
\end{pf}

\subsection{Nonlinear energy estimates}
This section is devoted to establishing the following Gronwall-type energy inequality for the solutions of the transformed MRT problem.
\begin{pro}  \label{pro:0301n0845}
Let $\Upsilon>0$ be provided by Proposition \ref{01dsaf16asdfasf0101} and $(\eta,u,q)$ be the local  solution constructed by Proposition  \ref{202102182115}. There exist $\delta_1^I$ and $ {c}_1^I>0$ such that, if $\|(\eta,\partial_1\eta,u)\|_4\leqslant \delta_1^I $ in some time interval $ I_{\tilde{T}} \subset I_T$, where $I_T$ is the existence time interval of $(\eta,u,q)$, then there exists a functional ${\mathfrak{E}}(t)$ of
$(\eta,u,q)$ satisfies the  Gronwall-type energy inequality
\begin{align}\label{2019120521430850}
{\mathfrak{E}}(t)\leqslant    {c}_1^I \left(I^0
+\int_0^t\|\eta_2(\tau)\|_{0}^2 \mm{d}\tau\right)+\Upsilon\int_0^t{\mathfrak{E}}(\tau)\mm{d}\tau
\end{align}
for a.e. $t\in I_{\tilde{T}}$,
where the constants $\delta_1^I$ may depend on $g$, $a$, $\lambda$, $m$, $\bar{\rho}$ and $\Omega$, and ${\mathfrak{E}}(t)\in W^{1,\infty}(I_{\tilde{T}})$ satisfies, for a.e. $t\in I_{\tilde{T}}$,
\begin{align}
\label{202103280958}
\mathfrak{E}, \ {\mathcal{E}}\mbox{ and }\|(\eta, \partial_1 \eta,u)\|^2_4\mbox{ are equivalent}.
\end{align}
\end{pro}
\begin{pf}
Let  $(\eta,u,q)$ be the local solution constructed by Proposition  \ref{202102182115}. We further assume that
\begin{align}
\label{xxxx}
\mathcal{E}(t)\leqslant \delta\in (0,1]  \mbox{ for any }t\in I_{\tilde{T}}\subset I_T.
\end{align}

Similarly to Lemma \ref{lem:08241445}, for sufficiently small $\delta$,  we can easily derive that, for a.e. $t\in I_T$
\begin{align}
&
\frac{\mm{d}}{\mm{d}t}\left(\int\bar{\rho}\partial_1^{k}\eta\cdot \partial_1^{k}u\mm{d}y
+\frac{a}{2}\|\sqrt{\bar{\rho}} \eta\|_{k,0}^2\right)
+c \| \eta\|_{k+1,0}^2\lesssim\| (\eta_2, u)\|_{k,0}^2
+\sqrt{\mathcal{E}}\mathcal{D}, \label{202008241446q}\\[1mm]
&\label{202008241448q}
 \frac{\mm{d}}{\mm{d}t} \|(\sqrt{ \bar{\rho} } u,
\sqrt{\lambda}m\partial_1 \eta)\|_{k,0}^2+ c \|  u\|_{k,0}^2
\lesssim \| \eta_2\|_{k,0}^2+ \sqrt{\mathcal{E}}\mathcal{D},
\end{align}
where $0\leqslant k\leqslant 4$.

We also verify that $(\eta,u)$  satisfies \eqref{202005021632} for a.e. $t\in I_{\tilde{T}}$.
Thus we derive from \eqref{202005021632} satisfied by $(\eta,u)$, \eqref{202008241446q} and \eqref{202008241448q}  that,
for some sufficiently  small $\varsigma\in (0,1]$,
\begin{align}\label{201910050940nxq}
\frac{\mm{d}}{\mm{d}t}\mathfrak{E}+c \varsigma
 \|( u,\partial_1 \eta)  \|_{3}^2  \lesssim  (1+\varsigma^{-1})\| \eta_2 \|_{\underline{3},0}^2+(1+\varsigma+\varsigma^{-1})
\sqrt{\mathcal{E}}\mathcal{D},
\end{align}
where
$$\begin{aligned}
&{\mathfrak{E}}:=
\varsigma\mathcal{E}^{\mm{cul}}+  {\varsigma}^{-1}  \|(\sqrt{ \bar{\rho} } u,
\sqrt{\lambda}m\partial_1 \eta)\|_{\underline{3},0}^2
+\frac{a}{2}\|\sqrt{\bar{\rho}}  \eta\|_{\underline{3},0}^2 + \sum_{k=0}^3 \int\bar{\rho}\partial_1^{k}\eta\cdot \partial_1^{k}u\mm{d}y\in  W^{1,\infty}(I_{\tilde{T}}).
\end{aligned}
$$

Similarly to \eqref{202012252005},  \eqref{2017020614181721nm}   and \eqref{201908081232},  for sufficiently small $\delta$, we have,  for a.e. $t\in I_{\tilde{T}}$,
\begin{align}
&\label{11111648}
\mathfrak{E}, \ {\mathcal{E}}\mbox{ and }\|(\eta, \partial_1 \eta,u)\|^2_4\mbox{ are equivalent}, \\
&{\mathcal{D}}\mbox{ is equivalent to }\|(u, \partial_1 \eta)\|^2_4,  \label{202104161417}
\end{align}
where the  equivalent coefficients in \eqref{11111648}  are independent of $\delta$.

Exploiting the interpolation inequality  \eqref{201807291850} yields, for any $\varepsilon\in (0,1]$,
\begin{align}
 \| \eta_2\|_{k,0}  \lesssim
                                         \begin{cases}
                       \varepsilon^{-1}  \| \eta_2\|_0 + \varepsilon \|  \eta_2\|_2   & \hbox{for }k=1; \\
     \varepsilon^{-({k-1})/(4-k)}\| \eta_2\|_{1,0} + \varepsilon \| \partial_1\eta_2\|_{3} & \hbox{for }2\leqslant k\leqslant 4.
                                         \end{cases}
  \label{202102181813}
\end{align}
Consequently, making use  of \eqref{11111648}--\eqref{202102181813}, we easily derive \eqref{2019120521430850} from \eqref{202012212151}  and  \eqref{201910050940nxq} for sufficiently small $\delta$.
\hfill$\Box$
\end{pf}

\subsection{Construction of nonlinear solutions}\label{202102180524}
For any given $\delta>0$, let
\begin{equation}\label{0501}
\left(\eta^\mm{a}, u^\mm{a}, q^\mm{a}\right)=\delta e^{\Upsilon t } (\tilde{\eta}^0, \tilde{u}^0, \tilde{q}^0),
\end{equation}
where $ (\tilde{\eta}^0, \tilde{u}^0, \tilde{q}^0)
:=(w/\Upsilon , w, \beta )$ and $(w, \beta ,\Upsilon )$ is provided by Proposition \ref{pro:08252100}.
Then $\left(\eta^\mm{a}, u^\mm{a}, q^\mm{a}\right)$
is also a  solution to the linearized MRT problem \eqref{01dsaf16asdfasf0101}, and enjoys the estimate, for any $i\geqslant0$,
\begin{equation}
\label{appesimtsofu1857}
\|\partial_{t}^i(\eta^\mm{a}, u^\mm{a})\|_5+\|\partial_{t}^i q^\mm{a} \|_4=\Upsilon^i \delta e^{\Upsilon  t}(\|(\tilde{\eta}^0,\tilde{u}^0)\|_5+\|\tilde{q}^0\|_4)\lesssim \Upsilon^i \delta e^{\Upsilon t}.
\end{equation}
In addition,
we have by \eqref{201602081445MH} that
\begin{eqnarray}\label{n05022052}
\| \chi_j\|_{0}\|\partial_1\chi_j\|_{0}\|\partial_2\chi_j\|_{0} >0,
\end{eqnarray}
where $\chi =\tilde{\eta}^0 $ or $\tilde{u}^0 $, and $j=1$, $2$.

Since the initial data of linear solution
$( {\eta}^{\mm{a}}, {u}^{\mm{a}}, {q}^{\mm{a}})$
does not satisfy the necessary compatibility conditions in the transformed MRT   problem in general.
Therefore, we shall modify the initial data of the linear solution.
\begin{pro}\label{pro:0101}
Let $(\tilde{\eta}^0, \tilde{u}^0 ):=(w/\Upsilon , w)$ be provided by \eqref{0501},
then there exists a constant $\delta_2^I \in (0,1]  $, such that for any $\delta\in(0, \delta_2^I ]$, there exists
$(\eta^{\mm{r}}, u^{\mm{r}})$ enjoying the following properties:
\begin{enumerate}
\item[(1)] The modified initial data
\begin{align}\nonumber
({\eta}^{\delta}_0,{u}^{\delta}_0 )
:=\delta(\tilde{\eta}^0,\tilde{u}^0 )+\delta^2(\eta^{\mm{r}},u^{\mm{r}} )
\end{align}
belongs to
$(H^{5}_1\cap H^{5}_{\mm{s}})\times H^5_{\mm{s}}$ and satisfies
the compatibility condition
\begin{align}\nonumber
&\mm{div}_{\mathcal{A}_0^{\delta}}u_0^{\delta}=0  \mbox{ in } \Omega,
\end{align}  where  $\mathcal{A}^{\delta}_0$ is defined as $\mathcal{A}$ with $\eta^{\delta}_0$ in place of $\eta$.
\item[(2)]
Uniform estimate:
\begin{align}
\label{202103281107}
\|(\eta^{\mm{r}},  u^{\mm{r}}  )\|_5\leqslant {c}_2^I,
\end{align}
where the positive constant $ {c}_2^I $ is independent of $\delta$.
\end{enumerate}
\end{pro}
\begin{pf}
Please refer to \cite[Lemma 4.2]{JFJSARMA2019}  or \cite[Proposition 5.1]{JFJSZWC}.
\hfill$\Box$
\end{pf}

Now we define that
\begin{align}
&\label{201912041727}
 {c}_3^I:= \|(\tilde{\eta}^0,\partial_1\tilde{\eta}^0, \tilde{u}^0)\|_4  + {c}_2^I >0 ,\\[1mm]
&\label{201912041705}
 {\delta}_0:= \min\left\{ \frac{{\delta}_2}{2c_3^Ic_4}, {\delta}_2^I, \frac{\delta_1^I}{2c_3^Ic_4^2} \right\}\leqslant 1  ,
\end{align}
where $c_4\geqslant 1$ is the constant in \eqref{202012212151}.

Let $\delta\leqslant  {\delta}_0$. Since $\delta\leqslant \delta_2^I$, we can use Proposition \ref{pro:0101} to construct $(\eta_0^\delta, u_0^\delta)$, which satisfies
$$ \|(\eta_0^\delta,\partial_1 \eta_0^\delta,u_0^\delta)\|_4
\leqslant c_3^I\delta\leqslant  {\delta}_2 .$$
By Proposition  \ref{202102182115}, there exists a (nonlinear) unique solution $(\eta, u, q)$ of the transformed MRT problem \eqref{01dsaf16asdfasf}
with initial value $({\eta}_0^\delta, {u}_0^\delta)$ in place of $(\eta^0, u^0)$, where $(\eta, u, q)\in \mathfrak{H} ^{1,4}_{1,*,\tau}\times   \mathfrak{U}_{\tau}^4\times  {\mathfrak{Q}} ^4_\tau  $ for any $\tau\in I_{T^{\mm{\max}}}$ and $T^{\mm{max}}$ denotes the maximal time of existence.

Let $\epsilon_0\in (0,1]$ be a constant, which will be given in \eqref{201907111842}.
We define
\begin{align}\label{times}
& T^\delta:={\Upsilon }^{-1}\mm{ln}({\epsilon_0}/{\delta})>0,\quad\mbox{i.e.,}\;
 \delta e^{\Upsilon  T^\delta }=\epsilon_0,\\[1mm]
&T^*:=\sup\left\{t\in I_{T^{\max}}\left|~
\sup\nolimits_{\tau\in [0,t)}\|(\eta,\partial_1 \eta,u )(\tau)\|_4\right.\leqslant 2   {c}_3^Ic_4\delta_0 \right\}, \label{xfdssdafatimes}\\[1mm]
&  T^{**}:=\sup\left\{t\in I_{T^{\max}} \left|~\sup\nolimits_{\tau\in [0,t)}\left\|\eta(\tau)\right\|_{0}\leqslant 2 {c}_3^I\delta  e^{\Upsilon  \tau}
 \right\}.\right.  \label{xfdstimes}
\end{align}
Since $(\eta,u)$ satisfies \eqref{202012212151} with $(\eta_0^\delta,u_0^\delta)$ in place of $(\eta^0,u^0 )$ for some $T\in I_{T^{\max}}$,
\begin{align}
c_4\left.\|(\eta,\partial_1\eta,u)(t)\|_4^2 \right|_{t=0}
=  c_4\|(\eta_0^\delta,\partial_1 \eta_0^\delta,u_0^\delta)\|_4
\leqslant c_3^Ic_4\delta  <  2c_3^Ic_4\delta ,
\label{201809121553}
\end{align} and
\begin{align}
\left.\| \eta (t)\|_3 \right|_{t=0}
=  \| \eta_0^\delta \|_3
\leqslant c_3^I\delta  <2 c_3^I\delta , \nonumber
\end{align}
thus  $T^{*}>0$, $T^{**}>0$ and
\begin{align}
\label{0502n111}  & \left\|\eta (T^{**}) \right\|_0
=2 {c}_3^I\delta e^{\Upsilon  T^{**}},\mbox{ if }T^{**}<T^{\max}.
\end{align}

Noting that $2c_3^Ic_4\delta \leqslant \delta_2 $, thus by Proposition  \ref{202102182115},
 we can check that
 \emph{if $T^*<\infty$, then
there exists a $T$ (independent of $\delta$) such that
\begin{align}
T^{*}_T:=T^*+T/2<T^{\mm{max}},\ \sup\nolimits_{\tau\in [0,T^{*}_T)} \|(\eta,\partial_1\eta ,u) (\tau)\|_3 \leqslant 2  {c}_3^I  {c}_4^2 \delta_0,  \label{0502n1}
\end{align}
and,  for any $\varsigma\in (T^{*},T^*_T)$, there always exists a $t_0 \in [ T^{*},\varsigma)$ such that}
\begin{align}
  \|(\eta,\partial_1\eta ,u) (t_0)\|_4 \geqslant 2  {c}_3^I  {c}_4 \delta_0 .\label{0502safdsan1}
\end{align}
From now on, we denote
$${T}^{*}_{**}:=
                  \begin{cases}
                T^{**} & \hbox{for } T^{**} \leqslant T^{*}; \\
                   T^*+ \min\{T^{**}-T^{*},T\}/2 & \hbox{for } T^{*}<T^{**},
                  \end{cases} $$ where $T$ comes from \eqref{0502n1}.

By \eqref{xfdssdafatimes} and \eqref{0502n1}, $(\eta,u)$ satisfies
$$
 \sup\nolimits_{ t\in {T}^{*}_{**}}\|(\eta,\partial_1 \eta,u)(t)\|_4\leqslant 2{c}_3^I  {c}_4^2 \delta_0 \leqslant \delta_1^I .$$
Thus, by Proposition \ref{pro:0301n0845}, $(\eta,u)$ enjoys  the  Gronwall-type energy inequality \eqref{2019120521430850}  for a.e. $t\in I_{{T}^{*}_{**}}$.
Using this fact, \eqref{xfdstimes} and   \eqref{201809121553},  we further have, for a.e. $t\in I_{{T}^{*}_{**}}$,
\begin{align}
\mathfrak{E} (t)  \leqslant c\delta^2e^{2\Upsilon  t}+\Upsilon \int_0^t \mathfrak{E} (\tau)  \mm{d}\tau
 \nonumber \end{align}
for some positive constant $c$.
Making use of the fact $\mathfrak{E} (t)\in W^{1,\infty}(I_{{T}^{*}_{**}})$, Gronwall's lemma, and the equivalence \eqref{202103280958}, we further deduce from the above estimate that, for a.e. $t\in I_{{T}^{*}_{**}}$,
\begin{align}
 \mathcal{E} (t) \lesssim \delta^2 e^{2\Upsilon  t},
\nonumber
\end{align}
which implies that,  for sufficiently small $\delta$,
\begin{align} \sup\nolimits_{ t\in [0 , T^{*}_{**})} ( \|(\eta,\partial_1\eta,u)\|_4+\|q\|_3+\|u_t\|_2 )\leqslant c_4^I \delta  e^{ \Upsilon  t}.
\label{20191204safda2114}
\end{align}

\subsection{Error estimates}\label{201907111049}
This subsection is devoted to the derivation of error estimates
between  $(\eta,u)$  and $(\eta^{\mm{a}}, u^{\mm{a}})$, where the (nonlinear) solution $(\eta,u,q)$   and the linear solution $(\eta^{\mm{a}}, u^{\mm{a}},q  ^{\mm{a}})$ are constructed in Section \ref{202102180524}.

\begin{pro}Let
$(\eta^{\mathrm{d}}, u^{\mathrm{d}},q^{\mathrm{d}})=(\eta, u,q)-(\eta^{\mm{a}}, u^{\mm{a}}, q^{\mm{a}})$, then there exists a constant $ {c}_5^I$ such that, for any $\delta\in(0,\delta_0]$,
\begin{align}\label{201907121534}
\sup\nolimits_{t\in I_{T^*_{**}}}\|\left(\eta^{\mm{d}}, u^{\mm{d}}\right)(t)\|_{W^{1,i}}\leqslant
 {c}_5^I\sqrt{\delta^3e^{3\Upsilon  t}}\mbox{ for }i=1\mbox{ and }2,
\end{align}
where $\delta_0$ is defined by \eqref{201912041705},  $ {c}_5^I$ only depends on   $g$, $a$,  $\lambda$, $m$, $\bar{\rho}$ and $\Omega$.
\end{pro}
\begin{pf}
Subtracting the both transformed MRT problem \eqref{01dsaf16asdfasf} and the linearized problem \eqref{01dsaf16asdfasf0101} with $\left(\eta^\mm{a}, u^\mm{a}, q^\mm{a}\right)$ in place of $(\eta, u, q)$, we get
\begin{equation}\label{201907121600}
 \begin{cases}
\eta_t^{\mm{d}}=u^{\mm{d}} ,\\[1mm]
\bar{\rho}u_t^{\mm{d}}+\nabla_{\mathcal{A}} q^{\mm{d}}+a\bar{\rho}u^{\mm{d}}-\Upsilon  m^2\partial_1^2\eta^{\mm{d}}
-g \bar{\rho}'\eta_2^{\mm{d}}\mathbf{e}_2
=g \mathcal{R}\mathbf{e}_2  -\nabla_{\tilde{\mathcal{A}}}q^{\mm{a}} ,\\[1mm]
\mm{div}_{{\mathcal{A}}}u^{\mm{d}}=-\mm{div}_{\tilde{\mathcal{A}}}u^{\mm{a}}  , \\[1mm]
\left(u^{\mm{d}},\eta^{\mm{d}}\right)|_{\partial\Omega}\cdot\vec{\mathbf{n}}=0  ,\\[1mm]
(\eta^{\mm{d}},u^{\mm{d}})|_{t=0}=\delta^2(\eta^{\mm{r}},u^{\mm{r}}),
\end{cases}
\end{equation}
where $\mathcal{R}:={ \int_{0}^{\eta_2}\left(\eta_2(y,t) -
z\right)\bar{\rho}''(y_2+z)\mm{d}z}$, $\tilde{\mathcal{A}}$ is given by \eqref{202103181309} and $\mathcal{A}=\tilde{\mathcal{A}}+I$.

Differentiating \eqref{201907121600}$_2$--\eqref{201907121600}$_4$ with respect to time $t$ and then using \eqref{201907121600}$_1$, we get
\begin{equation}\label{201907121602}
 \begin{cases}
\bar{\rho}u_{tt}^{\mm{d}}+\nabla_{\mathcal{A}} q_t^{\mm{d}}
+a\bar{\rho}u_t^{\mm{d}}-\Upsilon  m^2\partial_1^2u^{\mm{d}}-g\bar{\rho}'u_2^{\mm{d}}\mathbf{e}_2
\\
=g \mathcal{R}_t\mathbf{e}_2-\nabla_{\mathcal{A}_t}q-\nabla_{\tilde{\mathcal{A}} }q^{\mm{a}}_t=:\mathcal{N},\\[1mm]
\mm{div}_{{\mathcal{A}}}u_t^{\mm{d}}=-\mm{div}_{\mathcal{A}_t}u
-\mm{div}_{\tilde{\mathcal{A}}}u^{\mm{a}}_t , \\[1mm]
 u_t^{\mm{d}} |_{\partial\Omega}\cdot\vec{\mathbf{n}}=0 .
\end{cases}\end{equation}
Taking the inner product of \eqref{201907121602}$_1$ with $u_t^{\mm{d}}$ in  $L^2$, we obtain, for a.e. $t\in I_{{T}^{*}_{**}}$,
\begin{align}\label{201907092017}
\begin{aligned}
\frac{1}{2} \frac{\mm{d}}{\mm{d} t}\left(\|\sqrt{\bar{\rho}}u_{t}^{\mm{d}}\|_0^{2}-E(u^{\mm{d}})\right)
+a\left\|\sqrt{\bar{\rho}}u_{t}^{\mm{d}}\right\|_0^{2}
= \int \mathcal{N} \cdot  u_{t}^{\mm{d}}\mm{d}y-\int \nabla _{{\mathcal{A}}} q_{t}^{\mm{d}}\cdot u_{t}^{\mm{d}}\mm{d}y.
\end{aligned}
\end{align}

Exploiting the integral by parts, the boundary-value conditions of $(\eta,u,u^{\mm{a}})$, \eqref{201907121602}$_2$ and
 the relation
$$\partial_j\mathcal{A}_{ij}=0\mbox{ for }j=1,\ 2,$$
it is easy to compute out that
\begin{align}\nonumber
 -\int \nabla_{\mathcal{A}}q_{t}^{\mm{d}}\cdot u_{t}^{\mm{d}}\mm{d}y
&=\frac{\mm{d}}{\mm{d} t}
\int \nabla q^{\mm{d}}\cdot({\mathcal{A}}_t^{\mm{T}}  u + \tilde{\mathcal{A}}^{\mm{T}}  u^{\mm{a}}_t) \mm{d}y
 -\int  \nabla q^{\mm{d}}\cdot\partial_t({\mathcal{A}}_t^{\mm{T}}  u + \tilde{\mathcal{A}}^{\mm{T}}  u^{\mm{a}}_t)    \mm{d}y. \nonumber
\end{align}
Putting the above estimate  into \eqref{201907092017}, we conclude that
\begin{align}\label{201907092017nxq}
\begin{aligned}
&\frac{1}{2} \frac{\mm{d}}{\mm{d} t}\left(\|\sqrt{\bar{\rho}}u_{t}^{\mm{d}}\|_0^{2}-E(u^{\mm{d}})
-{I}_4\right)
+a\left\|\sqrt{\bar{\rho}}u_{t}^{\mm{d}}\right\|_0^{2}=I_5,
\end{aligned}
\end{align}
where we have defined that
\begin{align}
&{I}_4:=2\int\nabla q^{\mm{d}}\cdot\left({\mathcal{A}}_t^{\mm{T}}  u + \tilde{\mathcal{A}}^{\mm{T}}  u^{\mm{a}}_t  \right)\mm{d}y,\nonumber\\[1mm]
&I_5:=\int \mathcal{N}\cdot  u_{t}^{\mm{d}}\mm{d}y-\int   \nabla q^{\mm{d}}\cdot \partial_t({\mathcal{A}}_t^{\mm{T}}  u + \tilde{\mathcal{A}}^{\mm{T}}  u^{\mm{a}}_t) \mm{d}y\nonumber.
\end{align}

Thanks to  \eqref{appesimtsofu1857}, \eqref{202103281107}, \eqref{20191204safda2114} and \eqref{201907121600}$_5$,
it is easy to estimate that,   for all $t\in [0,{{T}^*_{**}})$,
\begin{align}\label{202008291723}
I_5(t) + {I}_4(t) - {I}_4(0) - E(u^{\mm{d}}|_{t=0}) \lesssim\delta^3e^{3\Upsilon  t}.
\end{align}
Taking the inner product of \eqref{201907121600}$_2$ with $t=0$ and $u_t^{\mm{d}}|_{t=0}$ in $L^2$,
and using the integration by parts and \eqref{201907121602}$_2$, we have
\begin{align}
 \int\bar{\rho}|u_t^{\mm{d}}|_{t=0}|^2\mm{d}y
 =&\int\left(g\mathcal{R}\mathbf{e}_2
-\nabla_{\tilde{\mathcal{A}}}q^{\mm{a}}+\Upsilon  m^2\partial_1^2\eta^{\mm{d}}
+g\bar{\rho}'\eta_2^{\mm{d}}\mathbf{e}_2-a\bar{\rho}u^{\mm{d}}\right)\cdot u_t^{\mm{d}}\mm{d}y \bigg|_{t=0}\nonumber \\
&+\int \nabla q^{\mm{d}}\cdot \mathcal{A}^{\mm{T}}_t(u+u^{\mm{a}})\mm{d}y \bigg|_{t=0}.\nonumber
\end{align}
Similarly to \eqref{202008291723}, we easily further deduce from the above identity that
\begin{align}\label{202008291630}
\|\sqrt{\bar{\rho}}u_t^{\mm{d}}|_{t=0}\|_0^2
\lesssim \delta^3e^{3\Upsilon  t}.
\end{align}
Using  \eqref{202008291723} and \eqref{202008291630}, we easily deduce from \eqref{201907092017nxq} that, for all $t\in [0,T^*_{**})$,
\begin{equation}\label{201907100815}
\begin{aligned}
\left\| \sqrt{\bar{\rho}}u_{t}^{\mm{d}}\right\|_0^{2}
+2a\int_0^t\left\| \sqrt{\bar{\rho}}u_{\tau}^{\mm{d}}(\tau)\right\|_0^{2}\mm{d}\tau
&\leqslant E(u^{\mm{d}}) +c\delta^3e^{3\Upsilon  t}.
\end{aligned}
\end{equation}
Next we shall deal with the term $E(u^{\mm{d}})$ in right hand side of \eqref{201907100815}.

By the existence theory of Stokes problem, for any given $t\in [0,T^*_{**})$, there exists a unique $(\tilde{u},\varpi)\in H^3\times \underline{H}^2$ such that
\begin{equation} \nonumber
\begin{cases}
\nabla \varpi- \Delta \tilde{u} =0,\;\;
\mm{div}\tilde{u} =-\mm{div}_{\tilde{\mathcal{A}}}u &\mbox{in }\Omega,\\[1mm]
\tilde{u} =0&\mbox{on } \partial\Omega.
\end{cases}
\end{equation}
Moreover, the solution $(\tilde{u},\varpi)$ enjoys
\begin{align}\label{201907121642}
\|\tilde{u}\|_3^2\lesssim
\|-\mm{div}_{\tilde{\mathcal{A}}}u\|_2^2\lesssim{\delta^4 e^{4\Upsilon  t}} .
\end{align}
It is easy to see that $v^{\mm{d}}:=u^{\mm{d}}-\tilde{u}\in H_{\sigma}^1$. Then we derive from \eqref{201907092040}  that
\begin{align}\nonumber
E(v^{\mm{d}})\leqslant ({\Lambda  ^2}+  {a\Lambda   }) \left\|\sqrt{\bar{\rho}}v^{\mm{d}}\right\|_0^2 ,
\end{align}
which, together with \eqref{201907121642}, implies that
\begin{align}\nonumber
E(u^{\mm{d}})\leqslant  ({\Lambda  ^2}+  {a\Lambda   }) \left\|\sqrt{\bar{\rho}}u^{\mm{d}}\right\|_0^2 +c\delta^3 e^{3\Upsilon t}.
\end{align}

Putting the above estimate  into \eqref{201907100815} yields
\begin{equation}\label{201907121715}
\begin{aligned}
\left\| \sqrt{\bar{\rho}}u_{t}^{\mm{d}}\right\|_0^{2}
+2a\int_0^t\left\|\sqrt{\bar{\rho}} u_{\tau}^{\mm{d}}(\tau)\right\|_0^{2}\mm{d}\tau
&\leqslant ({\Lambda  ^2}+  {a\Lambda   }) \left\|\sqrt{\bar{\rho}}u^{\mm{d}}(t)\right\|_0^2
+c\delta^3e^{3\Upsilon   t}.
\end{aligned}
\end{equation}
Using Newton--Leibniz's formula and Young's inequality, we have
\begin{equation}\label{201907100900}
\begin{aligned}
{a\Lambda   }\left\|\sqrt{\bar{\rho}}u^{\mm{d}}(t)\right\|_0^2
&\leqslant {a\Lambda   }\left\|\sqrt{\bar{\rho}}u^{\mm{d}}(0)\right\|_0^2
+\int_0^{t}\left\|(\sqrt{a\bar{\rho}}u_{\tau}^{\mm{d}}, \Lambda   \sqrt{a\bar{\rho}}u^{\mm{d}})(\tau)\right\|_0^2\mm{d}\tau.
\end{aligned}
\end{equation}
Thus we deduce from   \eqref{201907121715}--\eqref{201907100900}  that
\begin{equation}\label{201907100940}
\begin{aligned}
&\left\|\left (\sqrt{ {\bar{\rho}}{\Lambda  }^{-1}}u_{t}^{\mm{d}}, \sqrt{a\bar{\rho}}u^{\mm{d}}(t)\right)\right\|_0^2 \leqslant {\Lambda   }\left\|\sqrt{\bar{\rho}}u^{\mm{d}}(t)\right\|_0^2
+2a\Lambda  \int_0^{t}\left\|\sqrt{\bar{\rho}} u^{\mm{d}}(\tau)\right\|_0^2\mm{d}\tau
+c\delta^3e^{3\Upsilon   t}.
\end{aligned}
\end{equation}
In addition,
\begin{equation} \label{201907121720}
\frac{\mm{d}}{\mm{d}t}\| \sqrt{\bar{\rho}}u^\mm{d} \|^2_{0}
\leqslant \Lambda  ^{-1} \|\sqrt{\bar{\rho}}u_t^\mm{d} \|^2_{0}
+\Lambda  \|\sqrt{\bar{\rho}}u^\mm{d} \|^2_{0} .
\end{equation}
Combining \eqref{201907121720} with \eqref{201907100940}, we obtain the following differential inequality
\begin{equation}\label{201907101012}
\begin{aligned}
&\frac{\mathrm{d}}{\mathrm{d} t}\|\sqrt{\bar{\rho}} u^\mm{d}\|_{0}^{2}
+{a}\left\|\sqrt{\bar{\rho}}u^\mm{d}(t)\right\|_0^2 \leqslant {2\Lambda   }\left(\left\|\sqrt{\bar{\rho}}u^\mm{d}(t)\right\|_0^2
+a\int_0^{t}\left\|\sqrt{\bar{\rho}}u^\mm{d} (\tau)\right\|_0^2\mm{d}\tau\right)
+c\delta^3e^{3\Upsilon  t}.
\end{aligned}
\end{equation}

Applying Gronwall's inequality to \eqref{201907101012} then yields
\begin{equation}\nonumber
\begin{aligned}
\|u^\mm{d}(t)\|_{0}^{2}+  \int_{0}^{t}\|u^\mm{d}(\tau)\|_{0}^{2} \mathrm{d} \tau
&\lesssim\delta^3e^{3\Upsilon  t}.
\end{aligned}
\end{equation}
We further derive from \eqref{201907121600}$_1$ and the above estimate  that
\begin{equation}
\left\|\eta^\mm{d}(t)\right\|_{0}
\lesssim \sqrt{\delta^3e^{3\Upsilon  t}}. \label{2021003281325}
\end{equation}
Putting the above two estimates and \eqref{201907100815} together yields that
\begin{equation}
\left\|(\eta^\mm{d},\partial_1u^\mm{d},u^\mm{d},u^\mm{d}_t)(t)\right\|_{0}
\lesssim \sqrt{\delta^3e^{3\Upsilon  t}}.   \label{2021023281319}
\end{equation}
Similarly to  \eqref{2021003281325}, we also derive  from \eqref{201907121600}$_1$ and the above estimate
\begin{equation}
\left\| \partial_1\eta^\mm{d}(t)\right\|_{0}
\lesssim \sqrt{\delta^3e^{3\Upsilon  t}}.   \label{2021dfsafs023281319}
\end{equation}

Applying $\mm{curl}$ to \eqref{201907121600}$_2$, and then multiplying the resulting identity by $\mm{curl}u^{\mm{d}} $ in $L^2$, we get
 \begin{align}
 &\frac{1}{2}\frac{\mm{d}}{\mm{d}t}(\|\sqrt{\bar{\rho}} \mm{curl}u ^{\mm{d}}\|_0^2
+\Upsilon  \|m \mm{curl}\eta^{\mm{d}}\|_{1,0}^2)+a\|\sqrt{\bar{\rho}}\mm{curl}u^{\mm{d}}\|_0 \nonumber
\\
&=\int (\mm{cur}(g \mathcal{R}\mathbf{e}_2 -
\nabla_{\tilde{\mathcal{A}}} q^{\mm{d}})+\bar{\rho}' \partial_t u^{\mm{d}}_1 + a\bar{\rho}'u^{\mm{d}}_1+ g \bar{\rho}'\partial_1 \eta_2^{\mm{d}})\mm{curl}u^{\mm{d}} \mm{d}y,\nonumber
\end{align}
Thanks to  \eqref{2021023281319} and \eqref{2021dfsafs023281319}, we easily derive from the above identity that
 \begin{align}
 & \| \mm{curl}u ^{\mm{d}}\|_0^2
+ \| \mm{curl}\eta^{\mm{d}}\|_{1,0}^2 \lesssim {\delta^3e^{3\Upsilon  t}}. \nonumber
\end{align}
In addition, by \eqref{201907121600}$_3$, we have
 \begin{align}
 & \| \mm{\mm{div}}u ^{\mm{d}}\|_0^2  \lesssim {\delta^3e^{3\Upsilon  t}}.\nonumber
\end{align}
Thus we further derive from  the Hodge-type elliptic estimate \eqref{202005021302} and the above two estimates
that
 \begin{align} \nonumber
 & \| \nabla u ^{\mm{d}}\|_0^2  \lesssim {\delta^3e^{3\Upsilon  t}}.
\end{align}

Similarly to  \eqref{2021003281325}, we also derive  from \eqref{201907121600}$_1$ and the above estimate
 \begin{align}\| \nabla \eta ^{\mm{d}}\|_0^2  \lesssim {\delta^3e^{3\Upsilon  t}}. \nonumber
\end{align}
Putting \eqref{2021023281319} and the above two estimates together, and then using H\"older's inequality, we get \eqref{201907121534}. This completes the proof. \hfill $\Box$
\end{pf}
\subsection{Existence of  escape times}
Let
$$
\begin{aligned}
m_0 =& \min_{  \stackrel{\chi=\tilde{\eta}^0, \tilde{u}^0 }{1\leqslant  j\leqslant 2} } \{ \| \chi_j\|_{L^1},\|\partial_1\chi_j\|_{L^1},\|\partial_2\chi_j\|_{L^1} \}.
\end{aligned}
$$
Then $m_0>0$ by \eqref{n05022052}.

Now we define that
\begin{equation}\label{201907111842}
\epsilon_{0}=\min \left\{ \left(\frac{ {c}_3^I}{2 c_{5}^I }\right)^2, \frac{c_3^I c_4\delta_0}{c_4^I} , \frac{m_0^2}{4 |{c}_{5}^I|^{2}} \right\}>0.
\end{equation}
We claim that
\begin{equation}\label{201907111840}
T^{\delta}=T^{\mm{min}}=
\min \left\{T^{\delta}, T^{*}, T^{* *}\right\}
\neq T^{*}\;\;\mbox{or}\;\;T^{**},
\end{equation}
which can be showed by contradiction as follows:
\begin{enumerate}[(1)]
  \item If $T^{\mm{min}}=T^{**} $, then $T^{**}< T^{*}\leqslant T^{\max} $ or $T^{**}= T^{*} <+\infty$.
Noting that $\sqrt{\epsilon_0}\leqslant  {c}_3^I/2{c}_5^I$, then by
\eqref{0501}, \eqref{201912041727}, \eqref{times} and \eqref{201907121534}  that
\begin{equation*}\begin{aligned}
\|  \eta  (T^{**})\|_{0}&\leqslant \| \eta^\mm{a} (T^{**})\|_{0}+\| \eta^\mm{d} (T^{**})\|_{0}\\
&\leqslant  \delta e^{{\Upsilon  T^{**}}}(c_3^I+ c_5^I\sqrt{\delta e^{\Upsilon  T^{**}}})
\leqslant  \delta e^{{\Upsilon  T^{**}}}(c_3^I+ c_5^I\sqrt{\epsilon_0})\\
&\leqslant 3c_3^I \delta e^{\Upsilon  T^{**}}/2< 2c_3^I\delta e^{\Upsilon  T^{**}},
\end{aligned} \end{equation*}
which contradicts to \eqref{0502n111}.
  \item If $T^{\mm{min}}=T^{*}$, then $T^{*}<T^{**}  $. Recalling
$  \epsilon_0\leqslant {c}_3^I c_4\delta_0/c_4^{I}$, then we deduce from \eqref{20191204safda2114}
that, for any $t\in I_{T^*_{**}}$,
\begin{equation}\nonumber
 {\|(\eta,\partial_1 \eta , u )(t)\|_3  }  \leqslant c_4^I \delta  e^{ \Upsilon  T^\delta} \leqslant  {c}_3^Ic_4 \delta_0<2 {c}_3^Ic_4\delta_0,
\end{equation}
which contradicts \eqref{0502safdsan1}. Hence $T^{\mm{min}}\neq T^{*}$.
\end{enumerate}

Since $T^{\delta} $ satisfies \eqref{201907111840}, then \eqref{201907121534} holds to $t=T^\delta$.
Using this fact, \eqref{0501} and the condition $\epsilon_0\leqslant    {m_0^2}/ 4 |{c}_{5}^I|^{2}$, we  find that
$$
\begin{aligned}
\|\partial_j^k\chi_i  (T^\delta)\|_{L^1} \geqslant & \|\partial_j^k\chi_i^{\mm{a}}(T^\delta)\|_{L^1}-\|\partial_j^k\chi_i^{\mm{d}}(T^\delta)\|_{L^1} \\
\geqslant  &  \delta e^{\Upsilon  T^\delta }( \|\partial_j^k\tilde{\chi}_i^{0}\|_{L^1}
- c_5^I\sqrt{\delta e^{\Upsilon  T^\delta }})\geqslant  (m_0 -c_5^I\sqrt{\epsilon_0})\epsilon_0 \geqslant m_0 \epsilon_0 /2,
\end{aligned}
$$
where $\chi$ represents $\eta$ or $u$, $1\leqslant i$, $j\leqslant 2$ and $k=0$, $1$.
This completes the proof of  Theorem \ref{thm1} by taking $\epsilon:= m_0\epsilon_0 /2$.

\section{Local well-posedness}\label{202102241211}
 This section is devoted to the proof of local well-posedness of the transformed MRT problem \eqref{01dsaf16asdfasf}.  Recalling \eqref{202009130836}, the transformed MRT problem \eqref{01dsaf16asdfasf} is equivalent to the following initial-boundary value problem
\begin{equation}\label{01dsaf16asdfasf0000}
                              \begin{cases}
\eta_t=u ,\\[1mm]
\bar{\rho}u_t+\nabla_{\mathcal{A}} Q+a\bar{\rho} u=\lambda m^2\partial_1^2\eta-g\bar{\rho}\mathbf{e}_2 ,\\[1mm]
\div_{\mathcal{A}} u=0  , \\[1mm]
(\eta,u)|_{t=0}=(\eta^0,u^0) , \\[1mm]
(\eta,u)|_{\partial\Omega}\cdot\vec{\mathbf{n}}=0   ,
\end{cases}
\end{equation}
where $Q=q+\bar{P}(\zeta_2)+\lambda|\bar{M}|^2/2$.
Hence it suffices to establish the  local well-posedness result of the  initial-boundary value problem above.  Next we roughly sketch the proof frame to establish the  local well-posedness result.

Similarly to \cite{GXMWYJJMPA2019}, in which Gu--Wang investigated the well-posdeness problem of  incompressible inviscid MHD fluid equations with free boundary, we first alternatively  use an iteration method  to establish the existence result of  unique local solutions of the linear $\kappa$-approximate problem
\begin{equation} \label{01dsaf16asdfasf0000xx}
                              \begin{cases}
\eta_t-\kappa\partial_1^2\eta=u  ,\\
\bar{\rho}u_t+\nabla_{{\mathcal{B}}} Q+a\bar{\rho} u=\lambda m^2\partial_1^2\eta-g\bar{\rho}\mathbf{e}_2 ,\\
\div_{{\mathcal{B}}} u=0  , \\[1mm]
(\eta,u)|_{t=0}=(\eta^0,u^0) , \\[1mm]
(\eta,u)|_{\partial\Omega}\cdot\vec{\mathbf{n}}=0 ,\\ \nabla_{{\mathcal{B}}}Q|_{\partial\Omega}\cdot\vec{\mathbf{n}}=-g \bar{\rho}\mathbf{e}_2\cdot\vec{\mathbf{n}} ,
\end{cases}
\end{equation}
where $\mathcal{B}:=(\nabla \varsigma +I)^{-\mm{T}}$, and $\varsigma$ is given and belongs to some proper function space,
see Proposition \ref{202102151544} in Section \ref{subsec:local01} for details. It should be noted that the Neumann boundary-value condition \eqref{01dsaf16asdfasf0000xx}$_6$ makes sure the solvability of $Q$.

Then we derive the $\kappa$-independent estimates of the solutions, denoted by $(\eta^\kappa,u^\kappa, Q^\kappa)$, of the $\kappa$-approximate problem above,
and then  take the limit of $(\eta^\kappa,u^\kappa, Q^\kappa)$ with respect to $\kappa\to 0$ in some common definition domain. In particular, the obtained limit function, denoted by $(\eta,u,Q)$,
is the unique local  solution to the linearized problem
\begin{equation}
\label{01dsaf16asdfsafasf0000}
\begin{cases}
\eta_t =u  , \\[1mm]
\bar{\rho}u_t+\nabla_{{\mathcal{B}}} Q+a\bar{\rho} u=\lambda m^2\partial_1^2\eta-g\bar{\rho}\mathbf{e}_2  ,\\[1mm]
\div_{{\mathcal{B}}} u=0  , \\[1mm]
(\eta,u)|_{t=0}=(\eta^0,u^0)   , \\[1mm]
(\eta,u)|_{\partial\Omega}\cdot\vec{\mathbf{n}}=0  ,
\end{cases}
\end{equation}
see Proposition \ref{thm08} in Section \ref{subsubsce:03} for details. It should be noted that the solution $Q$ of \eqref{01dsaf16asdfsafasf0000} automatically satisfies \eqref{01dsaf16asdfasf0000xx}$_6$ due to $u_t |_{\partial\Omega}\cdot\vec{\mathbf{n}}=0$.

Finally, since the linearized problem \eqref{01dsaf16asdfsafasf0000} admits a unique local solution for any given $\varsigma$, thus we easily arrive at Proposition \ref{202102182115}  by a standard iteration method as in \cite{Coutand1,Coutand2}, in which  Coutand--Shkoller investigated the well-posdeness problem of  incompressible Euler equations with free boundary.

Now we turn to introducing some new notations appearing in this section.
 We denote $\Omega\times I_d$ by $\Omega_d$ for some constant $d>0$.
The notations $P(x_1,\ldots, x_n )$  and $\dot{P}(x_1,\ldots, x_n )$ represent the generic polynomials with respect to the parameters $x_1$, $\ldots$, $x_n$, where all the coefficients in $P$ and $\dot{P}$ are equal one, and $\dot{P}$ further satisfies $\dot{P}(0,\ldots,0)=0$.  It should be noted that $P(x_1,\ldots, x_n )$, $\dot{P}(x_1,\ldots, x_n )$ and $c^\kappa$ may vary from line to line. $a\lesssim_\kappa b$ means that $a\leqslant c^\kappa b$, where $c^\kappa$ denotes a generic positive constant,  which may depend on $\kappa$, $a$, $g$, $\lambda$, $m$, $\bar{\rho}$ and $\Omega$.

We always use the notations $\mathcal{B}$, resp. $\mathcal{J}$ to represent $(\nabla \varsigma +I)^{-\mm{T}}$, resp. $\det(\nabla \varsigma+I)$, where $\varsigma$ at least satisfies
\begin{align}
\varsigma \in {C}^0(\overline{I_T},{H}^4)\mbox{ and }\inf\nolimits_{(y,t)\in {\Omega_T}}\det(\nabla \varsigma +I)\geqslant 1/4 \mbox{ for some }T>0 .  \label{202safas2103250842}
\end{align}
In addition, we define that  \begin{align}
\mathbb{A}^{4,1/4}_{T,\iota }:=\{\psi\in  {C}^0(\overline{I_T},{H}^4_{\mm{s}})~|~  \|\psi\|_3\leqslant \iota\}  \label{2022104161633}
  \end{align}
and
  \begin{align}
 \mathbb{S}_T:=\{&(\xi,w,\beta) \in C(\overline{I_{T}} ,{H}^4_{ \mathrm{s}}) \times C^0( \overline{I_{T}},{H}^4_{\mm{s}})\times (C ^0(\overline{I_{T}},\underline{H}^3 ) \cap  L^\infty_TH^4)~|~\nonumber \\[1mm]
& \partial_1^2\xi\in  L^2_{T_1}H^4,\ \nabla^4\partial_1\xi \in C^0_{ B, \mm{weak}}( \overline{I_T},L^2), \ \nabla_{{\mathcal{B}}} \beta\in {L_{T}^2H^4},\nonumber  \\
& (\nabla_{{\mathcal{B}}} \beta/\bar{\rho})
|_{\partial\Omega}\cdot\vec{\mathbf{n}}=-g\mathbf{e}_2\cdot\vec{\mathbf{n}}  \} \label{20210413asfda21251},
\end{align}
for some constant $T>0$, $\iota$ is the positive constant provided in Lemma \ref{pro:1221},  $\bar{\rho}$  satisfies \eqref{0102} and  $g$ is the gravity constant. \emph{It should be noted that the function $\varsigma$, which belongs to $\mathbb{A}^{4,1/4}_{T,\iota }$, automatically satisfies \eqref{202safas2103250842} by Lemma \ref{pro:1221}.}

Finally, some preliminary estimates for $\mathcal{B}$ and $\mathcal{J}$ are collected as follows.
\begin{lem}\label{lem:0817}
Let $\varsigma $ satisfy \eqref{202safas2103250842} with $T>0$. Then
  \begin{align}
&\label{202008121535nnn}
\|{\mathcal{B}}-I\|_3 \lesssim \dot{P}( \|{\varsigma}\|_4),\\[1mm]
&\label{202008121505}
\|{\mathcal{B}}\|_{C^0(\overline{\Omega})}
+\|\mathcal{J}\|_{C^0(\overline{\Omega})}+\|(\mathcal{J},\mathcal{J}^{-1})\|_3
 \lesssim  P(\|{\varsigma}\|_4)\mbox{ for any }t\in \overline{I_T}.
\end{align}
\begin{enumerate}
   \item[(1)]  If additionally  $\varsigma$ further satisfies   $(\partial_1\varsigma, \varsigma_t)\in  L^\infty_TH^4$ , then
 \begin{align}
&\label{202008121535}
\|\partial_t(\mathcal{B}, \mathcal{J})\|_3\lesssim \dot{P}( \|(\varsigma,\varsigma_t)\|_4)  ,\\&\label{202008121535n}
 \|\partial_1 \mathcal{B} \|_3
 \lesssim \dot{P}( \| \varsigma \|_{\underline{1},4})\mbox{ for a.e. }t\in {I_T}  .\end{align}
   \item[(2)] If additionally  $\varsigma$ further satisfies
 $\varsigma_t\in  L^\infty_TH^3$ and $\varsigma_{tt}\in  L^\infty_TH^2$, then
 \begin{align}
&\label{202008121535xx}
\|\partial_{tt}(\mathcal{B}, \mathcal{J})\|_{1}\lesssim \dot{P}( \| \varsigma\|_4,\|\varsigma_t\|_3,\|\varsigma_{tt}\|_2)  \mbox{ for a.e. }t\in {I_T}  .\end{align}
 \end{enumerate}
\end{lem}
\begin{pf} Since the derivations of the estimates \eqref{202008121535nnn}--\eqref{202008121535xx} are very elementary, we omit it. \hfill$\Box$
\end{pf}

\subsection{Solvability of the linear $\kappa$-approximate problem \eqref{01dsaf16asdfasf0000xx}}\label{subsec:local01}

In this section, we construct the unique local solution of the the linear $\kappa$-approximate problem \eqref{01dsaf16asdfasf0000xx}. To this purpose, we shall rewrite  the linear $\kappa$-approximate problem \eqref{01dsaf16asdfasf0000xx} as the following equivalent problem (in the sense of classical solutions):
\begin{equation} \label{01dsaf16asdfasf0safd000}
 \begin{cases}
\eta_t-\kappa\partial_1^2\eta=u ,\\[1mm]
\bar{\rho}u_t+\nabla_{{\mathcal{B}}} Q+a\bar{\rho} u=\lambda m^2\partial_1^2\eta-g\bar{\rho}\mathbf{e}_2  ,\\[1mm]
-\mm{div}_{{\mathcal{B}}}\left(\nabla_{{\mathcal{B}}}Q/\bar{\rho}\right)=K^1  , \\[1mm]
(\eta,u)|_{t=0}=(\eta^0,u^0)  , \\[1mm]
(\eta,u)|_{\partial\Omega}\cdot\vec{\mathbf{n}}=0 ,\ \nabla_{{\mathcal{B}}}Q|_{\partial\Omega}\cdot\vec{\mathbf{n}}=-g \bar{\rho}\mathbf{e}_2\cdot\vec{\mathbf{n}} ,
\end{cases}
\end{equation}
where $\varsigma|_{t=0}=\eta^0$, $\eta^0$ satisfies
\begin{align}
\mm{div}_{\mathcal{A}^0}u^0=0,\ \mathcal{A}^0=(\nabla \eta^0+I)^{-\mm{T}},\  \inf\nolimits_{y \in \overline{\Omega}}(\det{\mathcal{A}^0})>0 , \nonumber
\end{align} and
\begin{align}
& K^1:=a\mm{div}_{{\mathcal{B}}}{u} -\mm{div}_{\mathcal{B}_t}u-\mathcal{J}_t\mm{div}_{{\mathcal{B}}}u/\mathcal{J}
  -\lambda m^2\mm{div}_{{\mathcal{B}}}\left(\partial_1^2{\eta}/\bar{\rho}\right) .  \label{202102141605}
\end{align}
Then the solvability of the linear $\kappa$-approximate problem  reduces to the solvability of the initial-boundary value problem \eqref{01dsaf16asdfasf0safd000}.

We want to establish the local well-posdeness result for \eqref{01dsaf16asdfasf0safd000} by an iteration method. To begin with, we shall investigate the solvability of the linear problem
\begin{equation} \label{01dsa000}
                              \begin{cases}
\eta_t-\kappa\partial_1^2\eta=w , \\[1mm]
u_t+a u= K^2 , \\[1mm]
-\mm{div}_{{\mathcal{B}}}\left(\nabla_{{\mathcal{B}}}Q/\bar{\rho}\right)=K^1 ,  \\[1mm]
(\eta,u)|_{t=0}=(\eta^0,u^0) ,   \\[1mm]
(\eta,u)|_{\partial\Omega}\cdot\vec{\mathbf{n}}=0 ,
\ \nabla_{{\mathcal{B}}}Q\cdot\vec{\mathbf{n}}|_{\partial\Omega}=-g \bar{\rho}\mathbf{e}_2\cdot\vec{\mathbf{n}}  ,
\end{cases}
\end{equation}
where $(\varsigma,w,\theta)$ is given, and
 \begin{align}
& K^2:=(\lambda m^2 \partial_1^2\xi-\nabla_{{\mathcal{B}}}\theta)/\bar{\rho}  -g \mathbf{e}_2 .  \nonumber
\end{align}

It is easy see that the above linear problem is equivalent to the following three sub-problems: the  initial-boundary value problem of partly dissipative equation for $\eta$
\begin{equation}\label{parabolic}
                              \begin{cases}
\eta_t-\kappa\partial_1^2\eta=w ,\\[1mm]
\eta|_{t=0}=\eta^0  ,\\[1mm]
\eta|_{\partial\Omega}\cdot\vec{\mathbf{n}}=0  ,
\end{cases}
\end{equation}
the initial-value problem of ODE  for $u$
\begin{equation}\label{transport}
                              \begin{cases}
 u_t+a u=K^2 ,\\[1mm]
u|_{t=0}=u^0  , \\[1mm]
u|_{\partial\Omega}\cdot\vec{\mathbf{n}}=0
\end{cases}
\end{equation}
and the Neumann boundary-value problem of elliptic equation for $q$
\begin{equation}\label{elliptic}
                              \begin{cases}
-\mm{div}_{{\mathcal{B}}}\left(\nabla_{{\mathcal{B}}}Q/\bar{\rho}\right)=K^1 &\mbox{in } \Omega,\\[1mm]
\nabla_{{\mathcal{B}}}Q/\bar{\rho}\cdot\vec{\mathbf{n}}=-g\mathbf{e}_2\cdot \vec{\mathbf{n}} &\mbox{on } \partial\Omega.
\end{cases}
\end{equation}
Thus the solvability of the linear problem \eqref{01dsa000} reduces to the solvability of the  three sub-problems above. Next we establish the global well-posedness results for the above three sub-problems in sequence.
 \begin{pro}\label{pro:parabolic}
Let $T>0$, $w\in L^2_{T}H^4_{\mm{s}}$
and $\eta^0\in H^{1,4}_{\mm{s}}$.
Then the initial-boundary value problem \eqref{parabolic} admits a unique solution $\eta\in C(\overline{I_T} ,{H}^4_{ \mathrm{s}}) $, which satisfies  $\nabla^4\partial_1\eta \in C^0_{B,\mm{weak}}( \overline{I_T}, L^2)$  and  the estimates
\begin{align}
&
\sup\nolimits_{t\in \overline{I_T}}(\| \eta \|_4+ \sqrt{\kappa}\|\partial_1\eta \|_4) +\kappa \|\partial_1^2\eta\|_{L^2_TH^4}
\nonumber \\
&\lesssim \| \eta^0 \|_4+ \sqrt{\kappa}\| \eta^0 \|_{1,4}+(1+ \sqrt{T})\|w\|_{L^2_TH^4} .\label{20200805100sadfsa5}
\end{align}
\end{pro}
\begin{pf}
We define the difference quotient with respect variable $y_1$ as follows:
$$D_1^{\tau}f (y_1,y_2)=(f(y_1+\tau,y_2 )-f(y_1,y_2 ))/\tau\mbox{ for }
\tau \mbox{ satisfyin }|\tau|\in (0,1).$$
Let $\varepsilon\in (0,1)$,
$\chi$ be a 1D standard mollifier (see \cite[pp. 38]{NASII04} for the definition), and $\chi^{\varepsilon}(s):=\chi(s/\varepsilon)/\varepsilon$.
Let $\tilde{w}=w$ in $\Omega_T $ and $\tilde{w}=0$ outside $\Omega\times (\mathbb{R}\backslash I_T)$. We define the mollification of $\tilde{w}$ with respect to $t$:
\begin{align}
\label{202210205162222}
S^{t}_\varepsilon (\tilde{w}): =\chi^{\varepsilon} * \tilde{w}.
\end{align}
Then
$S^{t}_\varepsilon  (\tilde{w}) \in C^\infty(\mathbb{R},H^4_{\mm{s}})$. We can check that
\begin{align}
&
\| S^{t}_\varepsilon (\tilde{w} ) \|_{L^2_TH^4}  \lesssim \| w\|_{L^2_TH^4}
\label{202210208asdfsa2248},\\
& S^{t}_\varepsilon  (\tilde{w}) \to  w \mbox{ strongly in }{L^2_TH^4}.
\nonumber
\end{align}

Now we consider the  $\tau$-approximate problem for \eqref{parabolic}:
\begin{equation}\label{parabsalic}
                              \begin{cases}
 \eta^\tau_t= {L}(\eta^\tau)+ S^{t}_\varepsilon  (\tilde{w} )
 &\mbox{in } \Omega_{T} ,\\[1mm]
\eta^\tau|_{t=0}=\eta^0 &\mbox{in } \Omega,\\
\eta ^\tau\cdot\vec{\mathbf{n}}=0  &\mbox{on }\partial\Omega,
\end{cases}
\end{equation}
where $L:$ $H^4_{\mm{s}}\to H^4_{\mm{s}}$ by the rule $L(f)= \kappa D^{-\tau}_1 D^\tau_1f$ for $f\in  H^4_{\mm{s}}$.

It is easy to see  that  $L(\varpi)\in  C^0(\overline{I_T},H^4_{\mm{s}})$ for $\varpi\in  C^0(\overline{I_T},H^4_{\mm{s}})$ and $L\in \mathcal{L}(H^{4}_{\mm{s}})$, where $\mathcal{L}(H^4_{\mm{s}})$ is a set of all linear bounded operators of $H^4_{\mm{s}}$. In particular, $ L\in  C^0(\overline{I_T},\mathcal{L}(H^4_{\mm{s}}))$. By existence theory of the initial-value problem of a abstract ODE equation (see \cite[Proposition 2.17]{NASII04}), there exists a unique solution $\eta^\tau\in C^0(\overline{I_T},H^4_{\mm{s}})\cap C^1( \overline{I_T},H^4_{\mm{s}})$ to \eqref{parabsalic}. Obviously
$\eta^\tau$,  $\eta^\tau_t$ automatically belong to $L^2_TH^4_{\mm{s}}$
by the second conclusion in Lemma \ref{20021032019018} and the separability of $H^{ 4}_{\mm{s}}$.

Let $\alpha$ satisfy $0\leqslant |\alpha|\leqslant  4$.
Applying $\partial^\alpha$ to \eqref{parabsalic}$_1$, and then multiplying the resulting identity $ \partial^\alpha \eta^\varepsilon$  in $L^2 $, we have, for a.e. $t\in I_T$,
\begin{align}
 \frac{\mm{d}}{\mm{d}t} \int |\partial^{\alpha}\eta^\tau|^2\mm{d}y
+ \kappa\int |D_1^h\partial^{\alpha}\eta^\tau |^2\mm{d}y
  =\int S^{t}_\varepsilon ( \partial^{\alpha} \tilde{w}) \cdot  \partial^{\alpha}\eta^\tau \mm{d}y . \label{202008051043}
\end{align}
Making use of  \eqref{202210208asdfsa2248}, we easily deduce from \eqref{202008051043} that,
\begin{align}
\|\eta^\tau\|_{C^0(\overline{I_T},H^4)} \lesssim\| \eta^0 \|_4+  \sqrt{T}\|w\|_{L^2_TH^4}.  \label{202102082127}
\end{align}

In addition, we easily deduce from  \eqref{parabsalic}$_1$ that
\begin{align}\label{2020safda08051043}
& \frac{\kappa}{2}\frac{\mm{d}}{\mm{d}t}\| D_1^{\tau}\partial^{\alpha}\eta^\tau\|^2_0
+  \int| \partial^{\alpha} \eta^\tau_t |^2\mm{d}y =\int \partial^{\alpha} (S^{t}_\varepsilon (w))\cdot  \partial^{\alpha} \eta^\tau_t \mm{d}y.
\end{align}
Noting that
\begin{align}
\label{2021032101049}
\| D_1^{\tau}\partial^{\alpha}\eta^\tau|_{t=0}\| _0\lesssim \| \partial^{\alpha}\eta^0\| _{1,0},
\end{align}
thus, making use of  Young's inequality, \eqref{202210208asdfsa2248} and \eqref{parabsalic}$_1$, we deduce from \eqref{2020safda08051043}  and \eqref{2021032101049} that
\begin{align}
 & {\kappa} \| D_1^{\tau}\eta^\tau \|_{C^0(\overline{I_T},H^4)}^2
+ \|({L}( \eta^\tau),  \eta^\tau_t)\|_{L^2_TH^4}^2
\lesssim  \kappa\| \eta^0 \|_{1,4}^2+  \|w\|_{{L^2_TH^4}}^2,  \label{202008051056}
\end{align}

Thanks to the regularity of $\eta^\tau$, we can we easily derive \eqref{parabsalic}$_1$ that, for any $\varphi\in H^1$,
\begin{align}
& \int D_1^{\tau}\partial^{\alpha} (\eta^\tau(y,t) -  \eta^\tau(y,s))\cdot\varphi\mm{d}y \nonumber \\
&=-
 \int_s^t\int ( \partial^{\alpha}S^{t}_\varepsilon ( \tilde{w}) +
L(\partial^{\alpha}\eta^\tau)) \cdot D_1^{-\tau}\varphi \mm{d}y \mm{d}\tau .
\label{202210041111912}
\end{align}
Exploiting  \eqref{202210208asdfsa2248}, the uniform estimate of $L( \eta^\tau)$ in \eqref{202008051056} and the fact
$$\|D_1^{-\tau}\varphi\|_0\lesssim \| \varphi\|_{1,0},$$  we easily derive from the identity \eqref{202210041111912} that
\begin{align}
&\label{20202108021asdf27}
 D_1^{\tau}\nabla^4 \eta^\tau\mbox{ is uniformly continuous in }H^{-1}.
\end{align}
Here and in what follows, $H^{-1}$ denotes the dual space of $H^1_0:= \{w\in {H}^1~|~w|_{\partial\Omega} =0\} $.

Making use of \eqref{202102082127}, \eqref{202008051056}, \eqref{20202108021asdf27} and \eqref{202104141653asfda}--\eqref{202104adsa141653}, there exists a subsequence (still denoted by $\eta^\tau$)   of $\{\eta^\tau\}_{|\tau|\in (0,1)}$ such that, for $\tau\to 0$,
\begin{align}
& \eta_t^\tau \rightharpoonup \eta_t^\varepsilon \mbox{ weakly in }L^2_TH^4_{\mm{s}}, \ \eta^\tau \to  \eta^\varepsilon  \mbox{ strongly in }C^0(\overline{I_T},H^3), \nonumber\\
& \eta^\tau \rightharpoonup  \eta^\varepsilon  \mbox{ weakly-* in }L^\infty_TH^4_{\mm{s}} ,\ D_1^\tau\nabla^4\eta^\tau \to\nabla^4\partial_1\eta^\varepsilon \mbox{ in }C^0(\overline{I_T},L^2_{\mm{weak}}),\nonumber\\
& D_1^\tau\eta^\tau \rightharpoonup \partial_1\eta^\varepsilon \mbox{ weakly-* in }L^\infty_TH^4\mbox{ and } D_1^{-\tau}D_1^\tau(\eta^\tau) \rightharpoonup \partial_1^2 \eta^\varepsilon \mbox{ weakly in }L^2_TH^4. \nonumber
\end{align}
Moreover, the limit function $\eta^\varepsilon$ is just the unique solution to the problem \begin{equation}\nonumber
                              \begin{cases}
 \eta^\varepsilon_t=\kappa\partial_1^2\eta^\varepsilon + S^{t}_\varepsilon  (\tilde{w} )
 &\mbox{in } \Omega_{T} ,\\[1mm]
\eta^\varepsilon|_{t=0}=\eta^0 &\mbox{in } \Omega,\\
\eta ^\varepsilon\cdot\vec{\mathbf{n}}=0  &\mbox{on }\partial\Omega
\end{cases}
\end{equation} and satisfies   \begin{align}
\nonumber
&\sup\nolimits_{t\in \overline{I_T}}(\| \eta^\varepsilon \|_4+\sqrt{\kappa}\|\partial_1\eta^\varepsilon \|_4) +\kappa \|\partial_1^2\eta^\varepsilon\|_{L^2_TH^4}\\
&\leqslant \| \eta^0 \|_4+ \sqrt{\kappa}\| \eta^0 \|_{1,4}+(1+ \sqrt{T})\|w\|_{L^2_TH^4}. \nonumber
\end{align}
Thus we further have $\eta^\varepsilon  \in C^0(\overline{I_T}, H^4_{\mm{s}})$.

Noting that $\eta^\varepsilon$ enjoys the uniform estimate above, thus we again get a limit function $\eta$, which satisfies the desired conclusion in Proposition \ref{pro:parabolic}, by taking limit of some sequence of $\{\eta^\varepsilon\}_{\varepsilon>0}$. We omit such limit process, since it is very similar to the argument of obtaining  $\eta^\varepsilon$.
\hfill $\Box$
\end{pf}
\begin{pro}\label{pro:transport}
Let $a\in \mathbb{R}$ and $K^2\in L^1_{T}H^4_{\mm{s}}$ and $u^0\in H^4_{\mm{s}} $, then the initial-boundary value problem \eqref{transport} admits a unique solution $u\in C^0(\overline{I_T},H^4_{\mm{s}})$, which satisfies $u_t\in L^1_TH^4_{\mm{s}}$ and
\begin{align}\label{202008051405}
\|u(t)\|_4\lesssim\|u^0\|_4+\|K^2\|_{L^1_TH^4}.
\end{align}
\end{pro}
\begin{pf}
Let
\begin{align}\label{202008051415}
 u
 = u^0e^{-at}+\int_0^tK^2(\tau)e^{-a(t-\tau)}\mm{d}\tau.
\end{align}
Since $K^2\in L^1_{T}H^4_{\mm{s}}$ and $u^0\in H^4_{\mm{s}} $, it is easy to check that $u$ given by \eqref{202008051415}  is the unique solution of  \eqref{transport}$_1$--\eqref{transport}$_2$; moreover, $u$ belongs to $ C^0(\overline{I_T},H^4_{\mm{s}})$ and satisfies \eqref{202008051405}.
\hfill $\Box$
\end{pf}

\begin{pro}\label{pro:elliptic}
Let $T>0$, $\iota$ be the positive constant provided in Lemma \ref{pro:1221},  $\bar{\rho}$ satisfy \eqref{0102} and $\varsigma\in \mathbb{A}^{4,1/4}_{T,\iota }$ defined by \eqref{2022104161633}. If
$K^1\in  C^0(\overline{I_T},H^1)\cap L^\infty_TH^2$ and satisfies $\int  K^1\mathcal{J}\mm{d}y=0$ for each $t\in \overline{I_T}$, then there exists a unique solution $Q\in C^0(\overline{I_T},\underline{H}^3)$, which satisfies the boundary-value problem \eqref{elliptic} for each $t\in I_T$ and enjoys the estimate
\begin{align}
\label{202008051safa442}
&\| Q\|_{L^\infty_TH^4}\lesssim P( \|{\varsigma}\|_{L^\infty_TH^4})\left(\|g\|_3+\|K^1\|_{L^\infty_TH^2}\right).
\end{align}
\begin{enumerate}
  \item[(1)] If we additionally assume that $K^1\in   L^2_TH^3$, then
\begin{align}\label{202008051442}
&\|\nabla_{{\mathcal{B}}}Q\|_{L^2_TH^4}\lesssim P( \|{\varsigma}\|_{L^\infty_TH^4})\left(\sqrt{T}\|g\|_4+\|K^1\|_{L^2_TH^3}\right).
\end{align}
  \item[(2)] If we additionally assume that $\varsigma_t \in L^\infty_TH^3$ and  $K_t^1\in L^\infty_TH^1$, then
\begin{align}
& \| Q_t\|_{L^\infty_TH^3}\lesssim P( \| {\varsigma}\|_{L^\infty_TH^4},\|{\varsigma}_t \|_{L^\infty_TH^3}) (\|g\|_3+\|K^1\|_{L^\infty_TH^2}+\|K^1_t\|_{L^\infty_TH^1}).  \label{2020104092126}
\end{align}
\end{enumerate}
\end{pro}
\begin{pf}
 Let $\varphi=\varsigma+y $, then $\varphi$ is a $C^1$-diffeomorphism mapping by Lemma \ref{pro:1221}.
 We denote $\varphi^{-1}$ the inverse mapping of $\varphi$ with respect to variable $y$, and then
define $\tilde{K}^1:=K^1 (\varphi^{-1},t)$.
By \eqref{202104031901}, $\tilde{K}^1\in C^0(\overline{I_T},H^1)\cap L^\infty_TH^2  $.

Now we consider the following Neumann boundary-value problem of elliptic equation
\begin{align}\label{elliptic1}
\begin{cases}
-\mm{div}\left(\nabla p/\tilde{\rho}\right)=\tilde{K}^1 &\mbox{in } \Omega,\\[1mm]
\nabla p/\tilde{\rho}\cdot\vec{\mathbf{n}}=-g\mathbf{e}_2\cdot \vec{\mathbf{n}}  &\mbox{on } \partial\Omega,
\end{cases}
\end{align}
where $\tilde{\rho}=\bar{\rho}(\varphi^{-1}_2)$. It is easy to see that
\begin{align}
\int \tilde{K}^1\mm{d}x=\int K^1\mathcal{J}\mm{d}y=0=-\int_{\partial\Omega}g\mathbf{e}_2\cdot \vec{\mathbf{n}}\mm{d}y_1. \nonumber
\end{align}
By Lemma \ref{lem:08181945}, the second assertion in Lemma \ref{20021032019018}, \eqref{202104032134}, \eqref{2022104101908}  and the separability of $\underline{H}^3$,  there exists
   a unique solution $p\in C^0(\overline{I_T},\underline{H}^3)\cap L^\infty_TH^4$ of the boundary-value problem \eqref{elliptic1} such that
\begin{align}\label{2020080516sasafddfa42}
\|p \|_{L^\infty_TH^4}\lesssim P(\|{\varsigma}\|_{L^\infty_TH^4})(\|g\|_3+\| {K}^1\|_{L^\infty_TH^2}).
\end{align}

Let $\tilde{Q}=p(\varphi)$, then $\tilde{Q}\in C^0(\overline{I_T}, {H}^3)\cap L^\infty_TH^4$ by \eqref{2020103250855}. We further define that $ {Q}:=\tilde{Q}- (\tilde{Q})_{\Omega} $, then $ {Q} \in C^0(\overline{I_T}, \underline{H}^3)\cap L^\infty_TH^4$ and satisfies \eqref{202008051safa442} by \eqref{2021sfa04031901}.
 Since $p$ satisfies \eqref{elliptic1}, we have
\begin{equation}\nonumber
                              \begin{cases}
-\mm{div}_{{\mathcal{B}}}\left(\nabla_{{\mathcal{B}}}Q/\bar{\rho}\right)=K^1&\mbox{in } \Omega,\\[1mm]
\nabla_{{\mathcal{B}}}Q/\bar{\rho}\cdot\vec{\mathbf{n}}=-g\mathbf{e}_2\cdot \vec{\mathbf{n}}  &\mbox{on } \partial\Omega.
\end{cases}
\end{equation}
In addition, the uniqueness of $Q$ is obvious in the class $   C^0(\overline{I_T}, \underline{H}^3)\cap L^\infty_TH^4$.

(1)  If we additionally assume that $K^1\in   L^2_TH^2$, then  $\tilde{K}^1\in  L^2_TH^2$. By Lemma \ref{20021032019018}, we have $p\in L^2_TH^4$ and
\begin{align}\label{2020080516sadfa42}
\|p \|_{L^2_TH^4}\lesssim P(\|{\varsigma}\|_{L^\infty_TH^3})(\sqrt{T}\|g\|_3+\|\tilde{K}^1\|_{L^2_TH^2}).
\end{align}
By \eqref{2020080516sadfa42}, \eqref{2021sfa04031901} and \eqref{2022104101908},
\begin{align}
&\|\nabla_{\mathcal{B}}Q\|_{L^2_TH^4}=\|\nabla p(x,t)|_{x=\varphi(y,t)}\|_{L^2_TH^4}\nonumber \\
&\lesssim P(\|\varsigma\|_{L^\infty_TH^4}) \|\nabla p\|_{L^2_TH^4}\lesssim P(\|{\varsigma}\|_{L^\infty_TH^4}) \left(\sqrt{T}\|g\|_4+\|{K}^1\|_{L^2_TH^3}\right), \nonumber
\end{align}
which yields \eqref{202008051442}.

(2) By \eqref{2021040sdaf32132}, it is easy to see that
\begin{align}
\label{202104131946}
\tilde{K}^1_t \in   L^\infty_TH^1\mbox{ and }\partial_t( \tilde{\rho},(1/\tilde{\rho}) )\in     L^\infty_TH^3.
\end{align}
Thanks to Lemma \ref{lem:08181945}, it is easy to verify that there exists a unique solution $\chi\in    L^\infty_TH^3 $ to
\begin{align}\label{ellsadfiptic1}
\begin{cases}
-\mm{div}\left(\nabla \chi/\tilde{\rho}\right)=\tilde{K}^1_t+\mm{div}\left(\partial_t(1/\tilde{\rho})\nabla p\right)  &\mbox{in } \Omega,\\[1mm]
\nabla \chi/\tilde{\rho}\cdot\vec{\mathbf{n}}=\partial_t(1/\tilde{\rho})  \nabla p\cdot\vec{\mathbf{n}}&\mbox{on } \partial\Omega.
\end{cases}
\end{align}
Moreover, by \eqref{2020080516sasafddfa42}, \eqref{neumaasdfann1n}, \eqref{2022104101908} and \eqref{2021040sdaf32132},
\begin{align}
& \|p\|_{L^\infty_TH^4} +\|\chi\|_{L^\infty_TH^3}\nonumber \\
&\lesssim P(\|{\varsigma}\|_{L^\infty_TH^4},\|{\varsigma}_t\|_{L^\infty_TH^3})(\|g\|_3+\| {K}^1
\|_{L^\infty_TH^2} + \|K_t\|_{L^\infty_TH^1} ).\label{ellsadfipsafatic1}\end{align}

Let $t\in I_T$ and $D_s\vartheta=\big(\vartheta(y,t+s)-\vartheta(y,t)\big)/s$ where $t+s\in I_T$. Since $p$ satisfies \eqref{elliptic1}, thus
\begin{align}\label{ellsadfiptic1sdfsf}
\begin{cases}
-\mm{div}\left(\nabla D_s p/\tilde{\rho}(x,t+s)\right)=D_s\tilde{K}^1+ \mm{div}\left(D_s(1/\tilde{\rho})\nabla p \right)&\mbox{in } \Omega,\\[1mm]
\nabla D_s p/\tilde{\rho}(x,t+s)\cdot\vec{\mathbf{n}}= -D_s(1/\tilde{\rho})\nabla  p \cdot\vec{\mathbf{n}}&\mbox{on } \partial\Omega.
\end{cases}
\end{align}
Subtracting \eqref{ellsadfiptic1} from \eqref{ellsadfiptic1sdfsf} yields
that
\begin{align}\label{ellsadfasdfasiptic1}
\begin{cases}
-\mm{div}\left(\nabla (\chi -D_s p)/\tilde{\rho}(x,t+s)\right)&\\
= \tilde{K}^1_t - D_s\tilde{K}^1 +\mm{div} ( ( \partial_t(1/\tilde{\rho})-D_s(1/\tilde{\rho})) \nabla p &\\
\quad +(1/\tilde{\rho}(x,t) -1/\tilde{\rho}(x,t+s))\nabla \chi  ) &\mbox{in } \Omega,\\[1mm]
(\nabla (\chi- D_s p)/\tilde{\rho}(x,t+s))\cdot\vec{\mathbf{n}} &\\
= (( \partial_t(1/\tilde{\rho} )-D_s(1/\tilde{\rho})) \nabla  p(x,t)&\\ \quad +(1 /\tilde{\rho} (x,t+s))-1/\tilde{\rho}(x,s)) \nabla \chi)\cdot\vec{\mathbf{n}}&\mbox{on } \partial\Omega.
\end{cases}
\end{align}
Applying the estimate \eqref{neumaasdfann1n} to \eqref{ellsadfasdfasiptic1},  we have,  for a.e. $t\in I_T$,
\begin{align}
 \| \chi- D_s p \|_3\lesssim  &\|\tilde{K}^1_t -D_s\tilde{K}^1  \|_1\nonumber
\\
& +\| ((\partial_t(1/\tilde{\rho}) -   D_s(1/\tilde{\rho}))  \nabla p  ,
 (1 /\tilde{\rho} (x,t+s) -1/\tilde{\rho}(x,s)) \nabla \chi  )\|_2.
\label{2022104092056}
\end{align}

Noting that the generalized derivative with respect to $t$ is automatically strong derivative, we easily see that the two terms on the right hand of the inequality \eqref{2022104092056} converge to $0$ for a.e. $t\in I_T$ by \eqref{202104131946}. So, $\|(D_sp-\chi)\|_3^2\to 0$ as $s\to 0$ for a.e. $t\in I_T$.
This means that the strong derivative of $ p$ with respect to $t$ is equal to that of $\chi$. In addition, it is easy to check that
$p \in AC^0(\overline{I_T},H^3)$, thus $p_t=\chi$, where $p_t$ denotes the generalized derivative of $p$.
Hence, $p_t\in  L^\infty(I_TH^3)$ satisfies \eqref{ellsadfipsafatic1} with $p_t$ in place of $\chi$. Thanks to \eqref{2021sfa04031901} and \eqref{202104032132}, we immediately get \eqref{2020104092126} from \eqref{ellsadfipsafatic1}.  This completes the proof of Proposition \ref{pro:elliptic}. \hfill $\Box$
\end{pf}

With Propositions \ref{pro:parabolic}--\ref{pro:elliptic} in hand, next we will use   an iteration method  to establish an existence result of a unique local solution to the linearized $\kappa$-problem \eqref{01dsaf16asdfasf0000xx}.

\begin{pro}
\label{202102151544}
Let $\mathbb{A}^{4,1/4}_{\alpha ,\iota }$ and $\mathbb{S}_{ \alpha  }$ are defined by \eqref{2022104161633} and \eqref{20210413asfda21251}  with some positive constant $\alpha  $ in place of $T$, resp..
 We assume that  $a\geqslant 0$,  $\bar{\rho}$  satisfy \eqref{0102},   $(\eta^0, u^0)\in {H}^{1,4}_{\mm{s}}\times {H}^{4}_{\mm{s}}$  and $\varsigma\in\mathbb{A}^{4,1/4}_{\alpha ,\iota }$  satisfies  $ \varsigma_t \in C^0(\overline{I_{\alpha }},H^2_{\mm{s}})\cap {L^{\infty}_{\alpha }H^4}$, then  the $\kappa$-approximate problem \eqref{01dsaf16asdfasf0000xx} defined on $\Omega_\alpha$  admits a unique solution, denoted by $(\eta^\kappa,u^\kappa,Q^\kappa)$, which belongs to $\mathbb{S}_{ \alpha  }$.
 \end{pro}
\begin{pf}
Let $T\leqslant \alpha$.  Thanks to the regularity of $(\varsigma,\varsigma_t)$ and the relation
\begin{align}
\partial_{j}(\mathcal{J} \mathcal{B}_{ij})=0\mbox{ for }i=1,\ 2,
\label{2020221042141629}
 \end{align}it is easy to check that, for any $ (\xi,w,\beta)\in \mathbb{S}_T$,
 $$K^1( \xi,w ) \in   C^0(\overline{I_T  },H^1)\cap L^\infty_T  H^1\cap L^2_T  {H}^3 \mbox{ and }\int  K^{1}( \xi,w )\mathcal{J}\mm{d}y=0\mbox{ for each }t\in \overline{I_T},$$
where $K^1( \xi,w)$ is defined by \eqref{202102141605} with $( \xi,w)$ in place of $(\eta,u)$.

By Proposition  \ref{pro:elliptic}, there exists a function $Q^1\in C^0(\overline{I_T},\underline{H}^3)\cap L^\infty_TH^4$ such that
$\nabla_{{\mathcal{B}}}Q^1\in {L^2_TH^4}$ and
\begin{align}
\nonumber
\begin{cases}
-\mm{div}_{{\mathcal{B}}}\left(\nabla_{{\mathcal{B}}}Q^1/\bar{\rho}\right)=0 &\mbox{in } \Omega,\\[1mm]
\nabla_{{\mathcal{B}}}Q^1/\bar{\rho}\cdot\vec{\mathbf{n}}=-g\mathbf{e}_2\cdot \vec{\mathbf{n}} &\mbox{on } \partial\Omega.
\end{cases}
\end{align}

 In view of Propositions \ref{pro:parabolic}--\ref{pro:elliptic} and the facts above, we easily see that there exist  a solution sequence $\{ (\eta^n,u^n,Q^n)\}_{n=1}^\infty$ defined on $I_T  $ and the solutions $(\eta^n,u^n,Q^n)$ enjoy the following properties:
\begin{enumerate}
  \item[(1)] $(\eta^1,u^1)=0$, and $Q^1$ satisfies $\nabla_{{\mathcal{B}}}{Q}^{1} |_{\partial\Omega}\cdot\vec{\mathbf{n}}=-g\mathbf{e}_2\cdot\vec{\mathbf{n}}$.
 \item[(2)]
$(\eta^n,u^n,Q^n)\in \mathbb{S}_{T  }$ for $n\geqslant 1$.
\item[(3)]
 for $n\geqslant2$, $(\eta^n,u^n,Q^n)$ satisfies the following relations:
\begin{equation}\nonumber
                              \begin{cases}
{\eta}_t^{n}-\kappa\partial_1^2{\eta}^{n}={u}^{n-1},\\[1mm]
{u}^{n}_t+a {u}^{n}= (\lambda m^2 \partial_1^2\bar{\eta}^{n-1}
-\nabla_{{\mathcal{B}}}{Q}^{n-1})/\bar{\rho} -g \mathbf{e}_2=:K^{2,n}  ,\\[1mm]
-\mm{div}_{{\mathcal{B}}}\left(\nabla_{{\mathcal{B}}}{Q}^{n}/\bar{\rho}\right) = K^{1}(\eta^n,u^n), \\[1mm]
({\eta}^{n},{u}^{n})|_{t=0}=(\eta^0,u^0) ,\\
({\eta}^{n}, {u}^{n})|_{\partial\Omega}\cdot\vec{\mathbf{n}}=0 ,\  \nabla_{{\mathcal{B}}}{Q}^{n} |_{\partial\Omega}\cdot\vec{\mathbf{n}}=-g\bar{\rho}\mathbf{e}_2\cdot\vec{\mathbf{n}},
\end{cases}
\end{equation}
where  $K^1( \eta^n,u^n)$ is defined by \eqref{202102141605} with $( \eta^n,u^n)$ in place of $(\eta,u)$, $K^{1,n} \in  C^0(\overline{I_T  },H^1)\cap L^\infty_T  H^2\cap L^2_T  {H}^3 $, $K^{2,n} \in  L^2_{T  }H^4_{\mm{s}}$,  and
\begin{align}
 \int  K^{1,n}\mathcal{J}\mm{d}y=0\mbox{ for each }t\in \overline{I_T} . \nonumber
\end{align}
\end{enumerate}

 We further define
$(\bar{\eta}^{n}, \bar{u}^{n},\bar{Q}^{n}):=(\eta^{n}-\eta^{n-1},u^{n}-u^{n-1},Q^{n}-Q^{n-1})$ for $n\geqslant 3$.
Then
\begin{equation}\nonumber
                              \begin{cases}
\bar{\eta}_t^{n}-\kappa\partial_1^2\bar{\eta}^{n}=\bar{u}^{n-1},\\[1mm]
\bar{u}^{n}_t+a \bar{u}^{n}= (\lambda m^2 \partial_1^2\bar{\eta}^{n-1}
-\nabla_{{\mathcal{B}}}\bar{Q}^{n-1})/\bar{\rho},\\[1mm]
-\mm{div}_{{\mathcal{B}}}\left(\nabla_{{\mathcal{B}}}\bar{Q}^{n}/\bar{\rho}\right)
&\\
=a\mm{div}_{{\mathcal{B}}}\bar{u}^{n}-\mm{div}_{\mathcal{B}_t}\bar{u}^{n}
-\mathcal{J}_t\mm{div}_{{\mathcal{B}}}\bar{u}^{n}/\mathcal{J}
-\lambda m^2\mm{div}_{{\mathcal{B}}}(\partial_1^2\bar{\eta}^{n}/\bar{\rho})  , \\[1mm]
(\bar{\eta}^{n},\bar{u}^{n})|_{t=0}=(0,0) ,\\
(\bar{\eta}^{n}, \bar{u}^{n})|_{\partial\Omega}\cdot\vec{\mathbf{n}}=0 ,\  \nabla_{{\mathcal{B}}}\bar{Q}^{n} |_{\partial\Omega}\cdot\vec{\mathbf{n}}=0.
\end{cases}
\end{equation}
 Thanks to  Propositions \ref{pro:parabolic}--\ref{pro:elliptic}, \eqref{202008121505} and \eqref{202008121535}, we easily estimate that, for $n\geqslant 3$,
\begin{align}
&
\| \bar{\eta}^{n}\|_{C^0(\overline{I_T  }, H^4_{\mm{s}})}+\|\partial_1\bar{\eta}^{n} \|_{L_{T  }^{\infty}H^4}
+\|\partial_1^2\bar{\eta}^{n}\|_{L_{T  }^2H^4}\nonumber \\
&\lesssim_{  {\kappa}}(1+\sqrt{T  })\|\bar{u}^{n-1}\|_{L_{T  }^2H^4}
\lesssim_{\kappa}T  (1+\sqrt{T  })\|\bar{u}^{n-1}\|_{C^0(\overline{I_T  },H^4_{\mm{s}})},\nonumber\\[2mm]
&\nonumber
\|\bar{u}^{n}\|_{C^0(\overline{I_T  },H^4_{\mm{s}})}
\lesssim \sqrt{T   } \|(\nabla_{{\mathcal{B}}}\bar{Q}^{n-1},\partial_1^2\bar{\eta}^{n-1})\|_{L_{T  }^2H^4}, \\
&\|\bar{Q}^{n}\|_{L_{T  }^\infty H^4} \lesssim
P ( \| ({\varsigma},{\varsigma}_t)\|_{L^\infty_\alpha  H^4} )
\left(  \|\bar{u}^{n}\|_{C^0(\overline{I_T  }, H^4_{\mm{s}})} +\|\partial_1 \bar{\eta}^{n}\|_{L_{T  }^\infty H^4} \right),\nonumber\\
&\|\nabla_{{\mathcal{B}}}\bar{Q}^{n}\|_{L_{T  }^2H^4} \lesssim
P ( \|({\varsigma}, {\varsigma}_t)\|_{L^{\infty}_{\alpha }H^4})
\left( \sqrt{T  }\|\bar{u}^{n}\|_{C^0(\overline{I_T  }, H^4_{\mm{s}})} +\|\partial_1^2\bar{\eta}^{n}\|_{L_{T  }^2H^4} \right).\nonumber
\end{align}

We immediately see from the above four estimates that  there exists a  sufficiently small $T  _1$  (the smallness  depends on $\kappa$,  $g$, $\lambda$, $m$,  $\bar{\rho}$, $\Omega$ and the norm $\|({\varsigma}, {\varsigma}_t)\|_{L^{\infty}_{\alpha }H^4}$) such that, for $T  = \min\{T  _1,\alpha \}$,
$$\{(\eta^n,\partial_1\eta^n,\partial_1^2\eta^n,u^n,Q^n,\nabla_{\mathcal{B}}Q^n)\}_{n=1}^\infty$$ is a Cauchy sequence in  $C^0(\overline{I_T  },H^4_{\mm{s}})\times L_{T  }^{\infty}H^4\times{L_{T  }^2H^4}\times C^0(\overline{I_T  },H^4_{\mm{s}})\times {L_{T  }^{\infty}\underline{H}^4}\times {L_{T  }^2H^4}$.  Thus we can get one limit function $(\eta,u,Q) $, which is the unique local solution of  the initial-boundary value problem  \eqref{01dsaf16asdfasf0safd000} and  also the unique local solution of  the linear  $\kappa$-approximate problem \eqref{01dsaf16asdfasf0000xx}.
In addition, it is easy to see that   $(\eta,u,Q)\in \mathbb{S}_{T  }$ by further using  Propositions \ref{pro:parabolic}, \ref{pro:elliptic} and trace theorem.

Noting that the local time $T  _1$ is independent of the initial data and  the local solution constructed above satisfies $(\eta,u)|_{t=T  }\in {H}^{1,4}_{\mm{s}}\times {H}^{4}_{\mm{s}} $, thus, if $T <\alpha$, we can further extend the local solution to be a global solution defined on $\overline{I_\alpha }$ by finite steps; moreover the obtained global solution is the unique solution of of  the linearized  $\kappa$-approximate problem \eqref{01dsaf16asdfasf0000xx} and  belongs to $\mathbb{S}_{\alpha}$.
This completes the proof  of Proposition \ref{202102151544}. \hfill$\Box$
\end{pf}

\subsection{Solvability of linearized  problem \eqref{01dsaf16asdfsafasf0000}}\label{subsubsce:03}

To investigate the solvability of linearized problem \eqref{01dsaf16asdfsafasf0000},
we shall first derive $\kappa$-independent estimates of the solutions of the linear problem \eqref{01dsaf16asdfasf0000xx}.

\begin{lem}\label{pro:0812}
Under the assumptions of Proposition \ref{202102151544},
we further assume that $\kappa \in (0,1]$, $\|\eta^0\|_4\leqslant \delta \in \mathbb{R}^+$ and $\varsigma$  satisfies
\begin{align}
&
 \partial_1\varsigma\in L^\infty_\alpha H^4,\   \varsigma|_{t=0}=\eta^0 .  \nonumber
\end{align}
 Then there exist  polynomials $\dot{P }( \|{\varsigma}\|_{L^{\infty}_{\alpha}H^4})$, $P  ( \|({\varsigma},\partial_1{\varsigma},  {\varsigma}_t)\|_{L^{\infty}_{\alpha}H^4} )$,  $P  ( I^0,\|({\varsigma},\partial_1{\varsigma},  {\varsigma}_t)\|_{L^{\infty}_{\alpha}H^4} )$ and positive constants $c$, $c_5\geqslant 1$,   $\delta_3$ (may depending on $g$, $a$,  $\lambda$, $m$, $\bar{\rho}$ and $\Omega$), such that, for any $\delta\leqslant \delta_3 $, the local solution  $(\eta^\kappa ,u^\kappa ,Q^\kappa )  $ provided by Proposition \ref{202102151544} enjoys the following estimates:
\begin{align}
&\label{2020081safdasaf21000}
\sup\nolimits_{t\in \overline{I_{T_1}}}\mathcal{E} ^\kappa(t) \leqslant 2 c_5(  I^0 +T\dot{P}( \|{\varsigma}\|_{L^{\infty}_{T_1}H^4}) )   ,\\[1mm]
&\label{202008121000}
\sup\nolimits_{t\in \overline{I_{T_1}}}\mathcal{E} ^\kappa(t)+    \kappa \|\partial_1^2\eta^\kappa \|_{L^2_{T_1}H^4}^2\leqslant 2 c_5I^0 +1 ,\\[1mm]
&\label{202008160930}
\|Q^\kappa  \|_{L^\infty_{T_1}H^4}+\|  u^\kappa_t \|_{L^\infty_{T_1}H^3}\lesssim P\left( I^0,\|({\varsigma},  {\varsigma}_t)\|_{L^{\infty}_{T_1}H^4}\right),
\end{align}
where $\mathcal{E}^\kappa (t):=\|(\eta^\kappa,\partial_1\eta^\kappa ,u^\kappa )(t)\|_4^2$,  $I^0$ is defined by \eqref{202104121056}, and
\begin{align}
T_1:=\min\{ 1/3c_5 P  ( \|({\varsigma},\partial_1{\varsigma},  {\varsigma}_t)\|_{L^{\infty}_{\alpha}H^4}  ) ,\alpha, 1\}   .
\label{20220103021516}
\end{align}
Moreover, for $(\eta^\kappa ,u^\kappa )$   restricted in $\overline{I_{T_1}}$,
\begin{align}
&\nabla^4\partial_1\eta^\kappa   \mbox{ and }  \nabla^4u^\kappa   \mbox{are uniformly continuous in }H^{-1}.
\label{202103021540}
\end{align}

  If additionally  $\varsigma_{tt}\in L^\infty_{T_1}H^2$, then
 \begin{align}
&\label{20200asfda8121535}
\|  Q^\kappa_t\|_{L^\infty_{T_1}H^3}\lesssim {P}( I^0,\|(\varsigma,\varsigma_t)\|_{L^\infty_{T_1}H^4},\|\varsigma_{tt}\|_{L^\infty_{T_1}H^2}) .\end{align}
\end{lem}
\begin{pf}
Frow now on we denote $(\eta^{\kappa} ,u^{\kappa} ,Q^{\kappa} )$   by $(\eta,u,Q)$, and let $T\leqslant \min\{\alpha,1\}$ and
 $$\|\eta^0\|_4\leqslant \delta \in (0,1].$$
Next we establish the desired uniform estimates for $(\eta,u,Q)$  by eight steps.

(1) \emph{Estimate of $\eta$.}

Recalling \eqref{20200805100sadfsa5}, we immediately get
\begin{align}\label{202008121045}
\sup\nolimits_{t\in\overline{I_T}}\|\eta (t)\|_4^2+\kappa^2\|\partial_1^2\eta \|^2_{L^2_{T}H^4}
 \lesssim \|\eta^0\|_{\underline{1},4}^2+ \|u \|_{L^2_{T}H^4}^2
\lesssim I^0+T\sup\nolimits_{t\in\overline{I_T}}\mathcal{E}^\kappa(t) .
\end{align}

(2) \emph{Estimate of $Q$.}

Noting that $Q$ satisfies
\begin{equation}\label{elliptasfasic}
                              \begin{cases}
-\mm{div}_{{\mathcal{B}}}\left(\nabla_{{\mathcal{B}}}Q/\bar{\rho}\right)=K^1 &\mbox{in } \Omega,\\[1mm]
\nabla_{{\mathcal{B}}}Q/\bar{\rho}\cdot\vec{\mathbf{n}}=-g\mathbf{e}_2\cdot \vec{\mathbf{n}} &\mbox{on } \partial\Omega,
\end{cases}
\end{equation}
where $K^1$ is defined by \eqref{202102141605}. By \eqref{202008121505}, \eqref{202008121535}, \eqref{202008051safa442}  and \eqref{202008121045}, we have
\begin{align}
\| Q\|_{L^\infty_TH^4}\lesssim & P( \|{\varsigma}\|_{L^\infty_TH^4})\left(\|g\|_3+\|K^1\|_{L^\infty_TH^2}\right)\nonumber \\
\lesssim & P\left( \|({\varsigma},   {\varsigma}_t)\|_{{L^\infty_TH^4}}\right)\left(1+
\sqrt{\mathcal{E}^\kappa(t)}\right).  \label{202008121055}
\end{align}

(3) \emph{$L^2$-norm energy estimate of $u$}.

 Let $q_\varsigma=Q-\bar{P}(\varsigma_2+y_2)$, and $G_\varsigma:=\bar{\rho}(\varsigma_2(y,t)+y_2)-\bar{\rho}(y_2)$. Then \eqref{01dsaf16asdfasf0000xx}$_2$ can be rewritten as follows:
\begin{equation}
\bar{\rho}u_t+\nabla_{{\mathcal{B}}}q_\varsigma +a\bar{\rho} u=\lambda m^2\partial_1^2\eta+g G_\varsigma\mathbf{e}_2   .\label{01dsasadff16asdfasf0000xx}
\end{equation}

Multiplying  \eqref{01dsasadff16asdfasf0000xx} with $\mathcal{J}u$ in $L^2$ then yields that
\begin{align}
&\frac{1}{2}\frac{\mm{d}}{\mm{d}t}\int\bar{\rho}|u|^2\mathcal{J}\mm{d}y
+a\int\bar{\rho}|u|^2 \mathcal{J}\mm{d}y\nonumber \\
& =\frac{1}{2} \int\bar{\rho}|u|^2\mathcal{J}_t\mm{d}y
  +\lambda m^2\int \partial_1^2\eta\cdot u  \mathcal{J}\mm{d}y +\int(g G_\varsigma-\nabla_{{\mathcal{B}}}q_\varsigma) \cdot u\mathcal{J}\mm{d}y.
\label{202008121145}
\end{align}

Using the boundary-value condition of $(\varsigma,u)$, the integration by parts, \eqref{01dsaf16asdfasf0000xx}$_3$, \eqref{202008121505} and \eqref{2020221042141629}, we have
\begin{align}
\label{202008121207}
 \int(g G_\varsigma-\nabla_{{\mathcal{B}}}q_\varsigma) \cdot u\mathcal{J}\mm{d}y =g \int G_\varsigma   u_2\mathcal{J}\mm{d}y\lesssim \dot{P}(\|\varsigma \|_3) \|u_2\|_0.
\end{align}
In addition, making use of \eqref{01dsaf16asdfasf0000xx}$_1$,
\eqref{202008121505}, \eqref{202008121535} and the integration by parts, we get
\begin{align}
& \lambda m^2\int\partial_1^2\eta\cdot u  \mathcal{J} \mm{d}y
 \nonumber \\
& =-\frac{\lambda m^2}{2}\frac{\mm{d}}{\mm{d}t}\int|\partial_1\eta|^2 \mathcal{J}\mm{d}y
+\lambda m^2 \int (|\partial_1\eta|^2 \mathcal{J}_t/2 -  \partial_1\eta\cdot \eta \partial_1\mathcal{J}- \kappa|\partial_1^2\eta|^2 \mathcal{J})\mm{d}y\nonumber \\
&\leqslant  -\frac{\lambda m^2}{2}\frac{\mm{d}}{\mm{d}t}\int|\partial_1\eta|^2 \mathcal{J}\mm{d}y
 - \kappa\lambda m^2 \int|\partial_1^2\eta|^2 \mathcal{J} \mm{d}y+ cP ( \|({\varsigma},{\varsigma}_t)\|_{ 3} ) \|\eta\|_2^2.\label{202008121210}
\end{align}

Noting that $1/4\leqslant \mathcal{J}\lesssim 1$,
thus, plugging \eqref{202008121207}--\eqref{202008121210} into \eqref{202008121145} and then integrating the resulting inequality over $(0,t)$, we immediately get, for any $t\in\overline{I_T}$,
\begin{align}
& \| (u, \partial_1\eta )(t) \|_0^2 +    \kappa \int_0^t\|\partial_1^2\eta (\tau)\|_{0}^2\mm{d}\tau\nonumber \\
&\lesssim I^0 +T\dot{P}(\|\varsigma \|_{L^{\infty}_{T}H^3 })+ TP\left( \|({\varsigma},{\varsigma}_t)\|_{L^{\infty}_{T}H^3}\right)  \sup\nolimits_{t\in\overline{I_T}}\|(\eta,u)(t)\|_2^2.\label{202008121220}
\end{align}

(4) \emph{Curl estimate of $(\eta,u)$.}

Applying $\mm{curl}_{{\mathcal{B}}}$ to \eqref{01dsaf16asdfasf0000xx}$_2$  yields that
\begin{align}\label{201910072117nbnm}
\partial_t\mm{curl}_{{\mathcal{B}}}\left(\bar{\rho}u\right)+a\mm{curl}_{{\mathcal{B}}}\left(\bar{\rho}u\right)
=\lambda m^2\mm{curl}_{{\mathcal{B}}}\left({\bar{\rho}}^{-1}\partial_1^2\left(\bar{\rho}\eta\right)\right)
+\mm{curl}_{ {{\mathcal{B}}}_t}\left(\bar{\rho}u\right)-g\bar{\rho}'{\mathcal{B}}_{12}.
\end{align}
Let the multiindex $\alpha$ satisfy $|\alpha|\leqslant2$. Applying $\partial^{\alpha}$
to \eqref{201910072117nbnm} yields
\begin{align}\label{202005021542nnbm}
\begin{aligned}
\partial_t\partial^{\alpha}\mm{curl}_{{\mathcal{B}}}\left(\bar{\rho}u\right)
+a\partial^{\alpha}\mm{curl}_{{\mathcal{B}}}\left(\bar{\rho}u\right)
=&{\lambda m^2}{\bar{\rho}}^{-1}\partial_1\partial^{\alpha}\mm{curl}_{{\mathcal{B}}}\left(\bar{\rho}\partial_1\eta\right)
+K^{\alpha}_3+K^{\alpha}_4,
\end{aligned}
\end{align}
where we have defined that
$$\begin{aligned}
&K^3_{\alpha}:=\lambda m^2\partial_1 \big([\partial^{\alpha}\mm{curl}_{{\mathcal{B}}}, {\bar{\rho}}^{-1}]\left(\bar{\rho}\partial_1 \eta\right)\big)
-\partial^{\alpha}\left(g\bar{\rho}'{\mathcal{B}}_{12}\right),\\[1mm]
&K^4_{\alpha}:= \partial^{\alpha}\mm{curl}_{{{\mathcal{B}}}_{t}}\left(\bar{\rho}u\right)
-\lambda m^2\partial^{\alpha}\mm{curl}_{\partial_1{{\mathcal{B}}}}\left(\partial_1\eta\right).
\end{aligned}$$

Multiplying \eqref{202005021542nnbm} by $\partial^{\alpha}\mm{curl}_{ {\mathcal{B}}}\left(\bar{\rho}u\right)$ in $L^2$, we obtain, for a.e. $t\in I_T$,
\begin{align}
&\frac{1}{2}\frac{\mm{d}}{\mm{d}t}\int|\partial^{\alpha}\mm{curl}_{{\mathcal{B}}}\left(\bar{\rho}u\right)|^2\mm{d}y
+a\int|\partial^{\alpha}\mm{curl}_{{\mathcal{B}}}\left(\bar{\rho}u\right)|^2\mm{d}y
\nonumber
\\&-\lambda m^2\int  {\bar{\rho}}^{-1}\partial_1\partial^{\alpha}\mm{curl}_{{\mathcal{B}}}\left(\bar{\rho}\partial_1\eta\right)
\partial^{\alpha}\mm{curl}_{{\mathcal{B}}}\left(\bar{\rho}u\right)\mm{d}y\nonumber
\\&
=\int K^3_{\alpha} \partial^{\alpha}\mm{curl}_{{\mathcal{B}}} \left(\bar{\rho}u\right)\mm{d}y
+\int K^4_{\alpha} \partial^{\alpha}\mm{curl}_{{\mathcal{B}}} \left(\bar{\rho}u\right)\mm{d}y
=:I_6 +I_7.\label{2020081214246}
\end{align}
Using the integral by parts and \eqref{01dsaf16asdfasf0000xx}$_1$ again, we obtain
\begin{align}
&-\lambda m^2\int  {\bar{\rho}}^{-1}\partial_1\partial^{\alpha}\mm{curl}_{{\mathcal{B}}}\left(\bar{\rho}\partial_1\eta\right)
 \partial^{\alpha}\mm{curl}_{{\mathcal{B}}}\left(\bar{\rho}u\right)\mm{d}y \nonumber \\
&=\frac{{\lambda}}{2}\frac{\mm{d}}{\mm{d}t}\int {\bar{\rho}}^{-1}| m\partial^{\alpha}\mm{curl}_{{\mathcal{B}}}\left(\bar{\rho}\partial_1\eta\right)|^2\mm{d}y
+ \kappa \lambda m^2\int {\bar{\rho}}^{-1}|\partial^{\alpha}\mm{curl}\left(\bar{\rho}\partial_1^2\eta\right)|^2\mm{d}y
+ I_8+I_9, \label{202008121536}
\end{align}
where we have defined that
$$
\begin{aligned}
&I_8:=\lambda m^2\int {\bar{\rho}}^{-1}\partial^{\alpha}\mm{curl}_{{\mathcal{B}}}\left(\bar{\rho}\partial_1\eta\right)
 \partial^{\alpha}(\mm{curl}_{\partial_1{\mathcal{B}}}\left(\bar{\rho}u\right)-
\mm{curl}_{\partial_t{\mathcal{B}}}\left(\bar{\rho}\partial_1\eta\right))\mm{d}y,\\
&I_9:=
\kappa\lambda m^2
\int\bigg( {\bar{\rho}}^{-1}\partial^{\alpha}\mm{curl}_{\partial_1{\mathcal{B}}}\left(\bar{\rho}\partial_1\eta\right)
 \partial^{\alpha}\mm{curl}_{{\mathcal{B}}}(\bar{\rho}\partial_1^2\eta) + |\partial^{\alpha}\mm{curl}_{{\mathcal{B}}-I}\left(\bar{\rho}\partial_1^2\eta\right)|^2  \\
&\qquad
+2 \partial^{\alpha} \mm{curl}_{{\mathcal{B}}-I}\left(\bar{\rho}\partial_1^2\eta\right)
\partial^{\alpha}\mm{curl}\left(\bar{\rho}\partial_1^2\eta\right)\bigg)\mm{d}y.
\end{aligned}
$$

 Thanks to the estimates \eqref{202008121535nnn}--\eqref{202008121535n}, we can estimate that
\begin{align}
&
I_6+I_7+I_8
\lesssim \dot{P}(\| \varsigma  \|_4) +   P(\|({\varsigma},\partial_1{\varsigma}, {\varsigma}_t)\|_4)
 \|(u,\partial_1\eta)\|_4^2,\nonumber \\[1mm]
&\nonumber
I_9\lesssim \kappa ( P(\|({\varsigma},\partial_1{\varsigma})\|_4)\|
\eta\|_{ {1},4}(\|\mm{curl}\left(\bar{\rho} \eta\right)\|_{2,3}+\|\eta\|_{1,4})
\nonumber \\
&\qquad +\|{\mathcal{B}}-I\|_3(\|{\mathcal{B}}-I\|_3\|  \eta \|_{2,4}^2
+  \| \eta \|_{2,4}\|\mm{curl}\left(\bar{\rho} \eta\right)\|_{2,3}).\nonumber
\end{align}
Noting that $\varsigma|_{t=0}=\eta^0$, by Newton--Leibniz formula and \eqref{202008121535nnn} with $t=0$, we have
\begin{align}
\|{\mathcal{B}}(t)-I\|_3\leqslant \left\|\left({\mathcal{B}}^0-I, \int_0^t{\mathcal{B}}_{\tau}(\tau) \mm{d}\tau\right)\right\|_3
\lesssim\|\eta^0\|_4+T P(\|({\varsigma}, {\varsigma}_t)\|_{L^{\infty}_{T}H^5}) ,
\label{202103020944}
\end{align}
where ${\mathcal{B}}^0:={\mathcal{B}}|_{t=0}$.
Making use of the above three estimates,  we deduce from \eqref{2020081214246}--\eqref{202008121536} that
\begin{align}
& \frac{\mm{d}}{\mm{d}t}
 \left\| \left(\mm{curl}_{{\mathcal{B}}}( \bar{\rho}u ), \sqrt{ {\lambda}/{\bar{\rho}}}m\mm{curl}_{{\mathcal{B}}}
\left(\bar{\rho}\partial_1\eta\right)\right)\right\|_3^2 +{c\kappa}\| \mm{curl}\left(\bar{\rho} \eta\right)\|_{2,3}^2\nonumber \\[1mm]
&\lesssim \dot{P}( \|\varsigma\|_4)
 + P( \|({\varsigma},\partial_1{\varsigma}, {\varsigma}_t)\|_4)\|(u,\partial_1\eta)\|_4^2 +   \kappa(\|\eta^0\|_4+T P(\|({\varsigma}, {\varsigma}_t)\|_{ L^{\infty}_TH^4}))\| \eta \|_{2,4}^2 . \label{202008121724}
\end{align}

In addition, similarly to \eqref{202103020944}, we have
\begin{align}
\|\mm{curl}f(t)\|_3&\lesssim\left \|\left(
  \mm{curl}_{\mathcal{B}(0)-I}f(t),
\mm{curl}_{{\mathcal{B}}}f(t),\mm{curl}_{\int_0^t {\mathcal{B}}_{\tau}\mm{d}\tau}f(t)\right)\right\|_3  \label{202104232112} \\[1mm]
&\lesssim \|\eta^0\|_4 \|f(t)\|_4+ T{P}( \|({\varsigma},{\varsigma}_t)\|_4)\|
 f(t)\|_4+\|
\mm{curl}_{{\mathcal{B}}}f(t)\|_3.\nonumber
\end{align}
Integrating \eqref{202008121724} over $(0,t)$, and then using  H\"older's inequality and the above estimate and \eqref{202103020944}, we conclude that
\begin{align}
&\| \mm{curl}(u , \partial_1\eta)(t)\|_3^2
+  \kappa \int_0^t\| \mm{curl}\left(\bar{\rho} \eta\right)(\tau)\|_{2,3}^2\mm{d}\tau\nonumber \\[1mm]
&\lesssim   I^0+T \dot{P}( \|{\varsigma}\|_{L^{\infty}_{T}H^4})+(\|\eta^0\|_4+
T P(  \|({\varsigma},\partial_1{\varsigma},  {\varsigma}_t)\|_{L^{\infty}_{T}H^4})
 \sup\nolimits_{t\in\overline{I_T}}\|(u,\partial_1\eta)\|_4^2  \nonumber \\
&\quad+ \kappa \left(\|\eta^0\|_4+T P(\|({\varsigma}, {\varsigma}_t)\|_{L^{\infty}_{T}H^4})\right)
  \int_0^t\| \eta (\tau)\|_{2,4}^2\mm{d}\tau+\|  (u , \partial_1\eta)(t)\|_3^2.
\label{202008121745}
\end{align}

(5) \emph{Divergence estimate of $(\eta,u)$.}

Similarly to \eqref{202104232112},   we have
\begin{align}
\|\mm{div}f(t)\|_3\leqslant \left\| \left(\mm{div}_{{\mathcal{B}}^0-I}f(t),\mm{div}_{{\mathcal{B}}}f(t),\mm{div}_{\int_0^t  {\mathcal{B}}_{\tau} \mm{d}\tau } f(t)\right) \right\|_3. \nonumber
\end{align}
Noting that $\mm{div}_{{\mathcal{B}}}u=0$ for any $t\in  \overline{I_T} $, taking $f=u$ in the above estimate yields
\begin{align}
\|\mm{div}u\|_3 \lesssim  (\|\eta^0\|_4+TP(\|({\varsigma},{\varsigma}_t)\|_4) ) \|u\|_4.\label{202008121802}
\end{align}

Applying $\mm{div}_{{\mathcal{B}}}\partial_1$ to \eqref{01dsaf16asdfasf0000xx}$_1$ and then using \eqref{01dsaf16asdfasf0000xx}$_3$, we have
\begin{align}\label{202008122040}
\partial_t(\mm{div}_{{\mathcal{B}}}\partial_1\eta)-\kappa\partial_1\mm{div}_{{\mathcal{B}}}\partial_1^2\eta=
\mm{div}_{\mathcal{B}_t}\partial_1\eta-\mm{div}_{\partial_1{\mathcal{B}}}u
-\kappa\mm{div}_{\partial_1{\mathcal{B}}}\partial_1^2\eta=: K _5.
\end{align}
We define the mollification of  $f\in L^2_TL^2$ with respect to  $y_1$ as follows:
 \begin{align}
\label{202121020281617}
 S_{ \varepsilon}^1 (f):=\chi^{\varepsilon} * \eta.
\end{align}
It is easy to check that
 \begin{align}
 S_{ \varepsilon}^1 (f)\to f \mbox{ in }L^2_TL^2.
\label{202121020281617safda}
\end{align}
Then we can derive from \eqref{202008122040} that
\begin{align}
&\partial_t S_{ \varepsilon}^1(\partial^\alpha \mm{div}_{{\mathcal{B}}}\partial_1\eta)-\kappa\partial_1 S_{ \varepsilon}^1 (\partial^\alpha\mm{div}_{{\mathcal{B}}}\partial_1^2\eta) =
S_{ \varepsilon}^1 (\partial^\alpha K _5),
\end{align}
where the multiindex $\alpha$ satisfies $|\alpha|\leqslant 3$.

For $\phi\in C_0^\infty(I_{T})$, we multiply the above identity by $S_\varepsilon^1(\partial^\alpha\mm{div}_{{\mathcal{B}}}\partial_1\eta)\phi$ in $L^2(\Omega_T)$ to get that
\begin{align}
& \int_0^T
\left(\frac{1}{2} \|S_{ \varepsilon}^1(\partial^{\alpha }\mm{div}_{\mathcal{B}}  \partial_1\eta  )\|^2_0  \phi_{\tau}+
\kappa  \|S_{ \varepsilon}^1 (\partial^\alpha\mm{div}_{{\mathcal{B}}}\partial_1^2\eta)\|_2^2\phi\right)\mm{d}\tau\nonumber \\
&=\int_0^T \int (S_{ \varepsilon}^1(\partial^\alpha K _5) S_{ \varepsilon}^1 (\partial^\alpha \mm{div}_{{  \mathcal{B}}}\partial_1\eta)-\kappa S_{ \varepsilon}^1 (\partial^\alpha \mm{div}_{{\mathcal{B}}}\partial_1^2\eta) S_{ \varepsilon}^1 (\partial^\alpha \mm{div}_{{\partial_1 \mathcal{B}}}\partial_1\eta))\mm{d}y
 \phi\mm{d}\tau . \nonumber
\end{align}
Thanks to \eqref{202121020281617safda}, we can take limits by $\varepsilon\to 0$ in the above identity to get that
\begin{align}
& \int_0^T
\left(\frac{1}{2} \| \partial^{\alpha }\mm{div}_{\mathcal{B}}  \partial_1\eta   \|^2_0  \phi_t+
\kappa  \|  \partial^\alpha \mm{div}_{{\mathcal{B}}}\partial_1^2\eta \|_2^2\phi\right)\mm{d}\tau\nonumber \\
&=\int_0^T \int (\partial^\alpha K _5  \partial^\alpha \mm{div}_{{  \mathcal{B}}}\partial_1 \eta -\kappa \partial^\alpha \mm{div}_{{\mathcal{B}}}\partial_1^2\eta \partial^\alpha \mm{div}_{{\partial_1 \mathcal{B}}}\partial_1\eta) \mm{d}y
 \phi\mm{d}\tau . \nonumber
\end{align}
In particular, we have, for a.e. $t\in I_T$,
\begin{align}
& \frac{1}{2} \frac{\mm{d}}{\mm{d}t}\| \partial^{\alpha }\mm{div}_{\mathcal{B}}  \partial_1\eta   \|^2_0 +
\kappa  \|  \partial^\alpha \mm{div}_{{\mathcal{B}}}\partial_1^2\eta  \|_2^2  \nonumber \\
&= \int (\partial^\alpha K _5 \partial^\alpha \mm{div}_{{  \mathcal{B}}}\partial_1\eta -\kappa \partial^\alpha \mm{div}_{{\mathcal{B}}}\partial_1^2\eta \partial^\alpha \mm{div}_{{\partial_1 \mathcal{B}}}\partial_1\eta )\mm{d}y .  \label{202104112141}
\end{align}
Follow the argument of \eqref{202008121745},  we can derive from \eqref{202104112141} that
\begin{align}
&\| \mm{div}\partial_1\eta  (t)\|_3^2
+ {\kappa} \int_0^t\|\mm{div} \eta(\tau)\|_{2,3}^2\mm{d}\tau\nonumber \\[1mm]
&\lesssim I^0+ (\|\eta^0\|_4^2+ T
P(\|({\varsigma},\partial_1{\varsigma}, {\varsigma}_t)\|_{L^{\infty}_{T}H^4})) \sup\nolimits_{t\in\overline{I_T}}\|(u,\partial_1\eta)\|_4^2
\nonumber \\[1mm]
&\quad
+ \kappa  \left(\|\eta^0\|_4+T P(\|({\varsigma},\varsigma_t)\|_{L^{\infty}_{T}H^4})\right)
\int_0^t\|  \eta(\tau)\|_{2,4}^2\mm{d}\tau\nonumber \\
&\quad+  \kappa P(\| \varsigma \|_{L^{\infty}_{T}H^{1,4}})
\int_0^t\| \eta(\tau)\|_{2,4}\| \partial_1\eta(\tau)\|_4\mm{d}\tau .
 \label{202008122100}
\end{align}
Consequently,  we immediately deduce from \eqref{202008121745}, \eqref{202008121802} and \eqref{202008122100} that
\begin{align}
&\|(\mm{div}u,\mm{curl}u,\mm{div}\partial_1\eta, \mm{curl}\partial_1\eta)(t)\|_3^2
+ { \kappa} \int_0^t\|(\mm{div} \eta,
  \mm{curl}(\bar{\rho} \eta))(\tau)\|_{2,3}^2\mm{d}\tau \nonumber \\[1mm]
&\lesssim   I^0+T \dot{P}( \|{\varsigma}\|_{L^{\infty}_{T}H^4})+(\|\eta^0\|_4+
T P(  \|({\varsigma},\partial_1{\varsigma},  {\varsigma}_t)\|_{L^{\infty}_{T}H^4})
 \sup\nolimits_{t\in\overline{I_T}}\|(u,\partial_1\eta)(t)\|_4^2  \nonumber \\
&\quad+ \kappa \left(\|\eta^0\|_4+T P(\|({\varsigma}, {\varsigma}_t)\|_{L^{\infty}_{T}H^4})\right)
  \int_0^t\|  \eta (\tau)\|_{2,4}^2\mm{d}\tau\nonumber \\
&\quad +  \kappa P(\| \varsigma \|_{L^{\infty}_{T}H^{1,4}})
\int_0^t\| \eta(\tau)\|_{2,4}\|\partial_1\eta(\tau)\|_4\mm{d}\tau+\|  (u , \partial_1\eta)(t)\|_3^2  \label{202008122105}
\end{align}

(6) \emph{Summing up the estimates $(\eta,u,Q)$.}

Thanks to the estimates \eqref{202008121045}, \eqref{202008121220}, \eqref{202008122105}, the interpolation inequality \eqref{201807291850}
and the Hodge-type elliptic estimate \eqref{202005021302}, we have, for sufficiently small $\delta$,
\begin{align}
 &\sup\nolimits_{t\in\overline{I_T}}\mathcal{E}^\kappa(t)+   \kappa \|\partial_1^2\eta \|_{L^2_{T}H^4}^2\nonumber \\
&\leqslant c_5( I^0+ T\dot{P}( \|{\varsigma}\|_{L^{\infty}_{T}H^4} ) ) \nonumber \\
&\quad + c_5 TP( \|({\varsigma},\partial_1{\varsigma}, {\varsigma}_t)\|_{L^{\infty}_{T}H^4})
\left(\sup\nolimits_{t\in\overline{I_T}}\mathcal{E}^\kappa(t)+   \kappa  \|\partial_1^2\eta  \|_{L^2_{{T}}H^4}^2\right)
 \label{202008122137}
\end{align}
and
 \begin{align}
 &\sup\nolimits_{t\in\overline{I_T}}\mathcal{E}^\kappa(t)+  \kappa \|\partial_1^2\eta  \|_{L^2_{T}H^4}^2\nonumber \\
&\leqslant c_5  I^0+    c_5 TP( \|({\varsigma},\partial_1{\varsigma}, {\varsigma}_t)\|_{L^{\infty}_{T}H^4})
\left(1+\sup\nolimits_{t\in\overline{I_T}}\mathcal{E}^\kappa(t)+  \kappa  \|\partial_1^2\eta  \|_{L^2_{{T}}H^4}^2\right).
 \label{202008122137xx}
\end{align}
 It should be noted that the two polynomials above are same.

Now we use $c_5$ and $P$  in \eqref{202008122137xx} to define $T_1$ by \eqref{20220103021516}, and thus get \eqref{2020081safdasaf21000}, resp. \eqref{202008121000} from \eqref{202008122137}, resp. \eqref{202008122137xx} by taking $T={T}_1$. Finally, making use of \eqref{01dsaf16asdfasf0000xx}$_2$ satisfied by $(\eta,u,Q)$,  \eqref{202008121000} and \eqref{202008121055}, we easily get \eqref{202008160930}.

(7) Let the multiindex $\beta$ satisfy $ |\beta|=4$.
Obviously, there exists $i$ such that $1\leqslant i\leqslant 4$ and $\beta _i\neq 0$. Let $\beta ^-$ satisfy $\beta ^-_i=\beta _i-1$ and $\beta ^-_j=\beta _j$ for $j\neq i$, and
$\beta ^+$ satisfy $\beta ^-_i=1$ and $\beta ^+_j=0$ for $j\neq i$.
Similarly to \eqref{202210041111912}, we can deduce from \eqref{01dsaf16asdfasf0000xx}$_1$ and \eqref{01dsaf16asdfasf0000xx}$_2$ that, for any $\varphi  \in H^1_0$ and for any $s$, $t\in I_T$,
\begin{align}
&\int\partial^\beta  \partial_1 (\eta(t)- \eta(s))
\cdot\varphi\mm{d}y=-\int_s^t\int\partial^\beta  (  u + \kappa \partial_1^2 \eta) \cdot \partial_1\varphi\mm{d}y\mm{d}\tau  ,  \label{2022104191625} \\[1mm]
&  \int \partial^\beta  (  \bar{\rho}u ( t) -
 \bar{\rho}u( s))  \cdot \varphi\mm{dy}\nonumber \\
& = \int_s^t \int
 ( \partial^{\beta^-} \nabla_{\mathcal{B}}Q \cdot \partial^{\beta^+}\varphi  - \lambda m^2\partial^\beta \partial_1 \eta \cdot \partial_1 \varphi -\partial^\beta  (g\bar{\rho}\mathbf{e}_2+ a\bar{\rho} u ) \cdot  \varphi )
\mm{d}y \mm{d}\tau.\label{2022104191625xxx}
\end{align}
Making use of the uniform estimates \eqref{202008121000} and \eqref{202008160930}, we easily deduce the assertion in \eqref{202103021540} from the two identities above.

(8)  If additionally  $\varsigma_{tt}\in L^\infty_{T_1}H^2$, we can apply the second conclusion in  Proposition \ref{pro:elliptic}  to \eqref{elliptasfasic}. Then, by further using the estimates \eqref{202008121535xx}, \eqref{2020104092126},  \eqref{202008121000} and \eqref{202008160930},
we can easily get \eqref{20200asfda8121535}.
This completes the proof.
\hfill $\Box$
\end{pf}

Thanks to  Lemma \ref{pro:0812}, we establish the unique local solvability of the $\kappa$-approximate problem \eqref{01dsaf16asdfsafasf0000}  by a compactness argument.
\begin{pro}\label{thm08}
Let the assumptions of Lemma \ref{pro:0812} be satisfied,  $ \delta_0$, $T_1$ be provided by Lemma \ref{pro:0812} and $\varsigma_{tt}\in L^\infty_{T_1}H^2$, then, for any $\delta\leqslant \delta_3$,  the linearized problem
  \eqref{01dsaf16asdfsafasf0000}
defined on $\Omega_{T_1}$ admits a unique solution
$(\eta^{{L}},u^{{L}}, {Q}^{L} )\in \mathfrak{C}^0(\overline{I_{T_1}},{H}^{1,4}_{\mm{s}} )\times   \mathfrak{U}_{T_1}^4\times   \mathfrak{Q}_{T_1}^4   $; moreover the solution satisfies
\begin{align}
&\label{2020081asdfa21045}
\sup\nolimits_{t\in {I_{T_1}}}\|\eta (t)\|_4
 \leqslant \| \eta^0 \|_4+ \sqrt{{T}(2c_5 I^0+1)},\\
&\label{2020081safsafadasaf21000}
\sup\nolimits_{t\in {I_{T_1}}}\mathcal{E}^L(t) \leqslant 2 c_5( I^0+ T\dot{P}( \|{\varsigma}\|_{L^{\infty}_{T}H^4}) )   ,\\
&\label{202008131221}
\sup\nolimits_{t\in  {I_{T_1}}}\mathcal{E}^L(t)\leqslant 2 c_5 I^0+1,\\
&\label{202008160asdadd954}
  \|{ {Q}^{L}} \|_{L^\infty_{T_1}H^4}+\| u^{{L}}_t \|_{L^\infty_{T_1}H^3}
\lesssim P\left(  I^0,\|({\varsigma},  \varsigma_t)\|_{L^{\infty}_{T_{\alpha} }H^4} \right),\\
&\label{202008160954}
  \| Q^{{L}}_t\|_{L^\infty_{T_1}H^3}
\lesssim P\left(  I^0,\|({\varsigma},  \varsigma_t)\|_{L^{\infty}_{T_{\alpha} }H^4},\|\varsigma_{tt}\|_{L^\infty_{T_1}H^2}\right),
\end{align}
where $\mathcal{E}^L(t):= \sup\nolimits_{t\in {I_{T_1}}}\|(\eta^{{L}},\partial_1\eta^{{L}} ,u^{{L}})(t)\|_4^2 $.
\end{pro}
\begin{pf}
Let  $(\eta^\kappa,u^\kappa,Q^\kappa)\in\mathbb{S}_{{T} }$ be the  solution of the linearized problem \eqref{01dsaf16asdfasf0000xx} stated as in Lemma \ref{pro:0812}. Thanks to the $\kappa$-independent estimates \eqref{2020081safdasaf21000}--\eqref{202008160930}, \eqref{202103021540} and \eqref{20200asfda8121535}, we can easily follow the compactness argument as in the proof of Proposition \ref{pro:parabolic} to obtain a limit function $(\eta^{{L}},u^{{L}},Q^{{L}})$, which is the solution of the linearized problem
 \eqref{01dsaf16asdfsafasf0000}, satisfies \eqref{2020081safsafadasaf21000}--\eqref{202008160954} and belongs to $ \mathfrak{C}^0(\overline{I_{T_1}},{H}^{1,4}_{\mm{s}} )\times   \mathfrak{U}_{T_1}^4\times  \mathfrak{Q}_{T_1}^4$.

The estimate \eqref{2020081asdfa21045} is obvious by \eqref{01dsaf16asdfsafasf0000}$_1$ satisfied by $(\eta,u)$  and  \eqref{202008131221}.
In addition, the uniqueness can be easily verified by a standard energy method.  The proof of Proposition \ref{thm08} is complete.
\hfill$\Box$
\end{pf}

\subsection{Proof of Proposition \ref{202102182115}}\label{subsce:03}
Let $(\eta^0,u^0)$ satisfy the assumptions in Proposition \ref{202102182115}, $\|\eta^0\|_4\leqslant  \delta\leqslant \max\{\delta_3, \iota\}/2$, and
\begin{align}
\nonumber
{ T_2}:=  \min\{1/3c_5 P  (c_4^2 (b+\delta_3)^2) , \delta^2/ (2c_5(b+\delta_3)+1),1/4c_5  \}< 1,
\end{align}
where the positive constant $c_5\geqslant 1$ and the polynomial are provided by \eqref{20220103021516} and $\delta_3$ is the constant in Lemma  \ref{pro:0812} with $\alpha={ T_2}$. By Lemma \ref{pro:0812}, Proposition \ref{202102151544}  and Lemma \ref{pro:1221}, for any $T\leqslant T_2$, we can construct a solution sequence
$$\{(\eta^{n}, u^{n}, Q^{n})\}_{n=1}^\infty\subset   \mathfrak{C}^0 (\overline{I_{T }},{H}^{1,4}_{\mm{s}} )\times   \mathfrak{U}_{T }^4\times \mathfrak{Q}_T^4$$ such that
\begin{enumerate}[(1)]
\item $( \eta^1, u^{1},Q^{1})=(\eta^0, u^0,0)$.
  \item $(\eta^{n+1}, u^{n+1}, Q^{n+1})$ satisfies
\begin{equation}\label{01dsaf16asdfasf00002}
                              \begin{cases}
 \eta^{n+1}_t =u^{n+1} ,\\[1mm]
\bar{\rho}u^{n+1}_t+\nabla_{\mathcal{A}^{n}}Q^{n+1}+a\bar{\rho}u^{n+1}=\lambda m^2 \partial_1^2\eta^{n+1}-g\bar{\rho}\mathbf{e}_2   ,\\[1mm]
\div_{{\mathcal{A}}^{n}}u^{n+1}=0   , \\[1mm]
(\eta^{n+1},u^{n+1})|_{t=0}=(\eta^0,u^0) , \\[1mm]
(\eta^{n+1},u^{n+1})\cdot\vec{\mathbf{n}}=0 ,\ \nabla_{\mathcal{A}^{n}}Q^{n+1}\cdot\vec{\mathbf{n}}=
-g\mathbf{e}_2 \bar{\rho}\cdot\vec{\mathbf{n}}  \mbox{ on } \partial\Omega
\end{cases}
\end{equation}
for any $n\geqslant 1$, where $\mathcal{A}^{n}=(\nabla(\eta^{n}+I))^{-\mm{T}}$.
  \item for $n\geqslant 2$, $(\eta^{n
}, u^{n}, Q^{n})$ enjoys the following estimates:
\begin{align}
&  \sup\nolimits_{t\in {I_{T }}} \| \eta^{n} (t)\|_4   \leqslant \| \eta^0 \|_4+  \sqrt{T (2c_5 I^0+1) }\leqslant 2\delta\leqslant \min\{\delta_3, \iota\}\leqslant 1, \label{20dsfa2008131747} \\[1mm]&\sup\nolimits_{t\in \overline{I_{T }}}\mathcal{E}^n(t) \leqslant  c_4^2 I^0,  \label{202saf008131740n} \\
&\sup\nolimits_{t\in {I_{T}}} \mathcal{E}^n(t)+
 \|Q^{n}\|_{L^\infty_TH^4}+\| \partial_t(u^{n},Q^n)  \|_{L^\infty_TH^3}
\lesssim P (  I^0 ),
 \label{202008131740n}
\end{align}
where $\mathcal{E}^n(t):= \|(\eta^{n},\partial_1\eta^n,u^n)(t)\|_3^2 $.
\end{enumerate}

Similarly to \eqref{202008121505}--\eqref{202008121535}, by using \eqref{20dsfa2008131747}  and
\eqref{202008131740n}, we have, for any $n\geqslant 1$,
\begin{align}
&\label{202008141600}
\sup\nolimits_{t\in  {I_{T}}} \|( \mathcal{A}^{n},(J ^{n})^{-1},\mathcal{A}^{n}_t, J^{n}_t)\|_3\lesssim P(I^0),
\end{align}
where $J^n:=\det(\nabla \eta^n+I)$.
In addition, similarly to \eqref{2022104191625} and \eqref{2022104191625xxx}, by a regularity method, we can easily derive from \eqref{01dsaf16asdfasf00002}$_1$ and \eqref{01dsaf16asdfasf00002}$_2$ that, for any $\varphi \in H^1_0$, $(\eta^{n}, u^{n})$  satisfies \begin{align}
&\int\partial^\beta \partial_1 (\eta^{n}(t)- \eta^{n}(s))
\cdot\varphi\mm{d}y=-\int_s^t\int  \partial^\beta  u^{n}    \cdot \partial_1\varphi \mm{d}y\mm{d}\tau  ,  \nonumber\\[1mm]
&  \int \partial^\beta (  \bar{\rho}u^{n} (t) -
 \bar{\rho}u^{n}(s))  \cdot \varphi\mm{dy}\nonumber \\
& = \int_s^t\int
 ( \partial^{ \beta^-} \nabla_{\mathcal{A}^{n-1}}Q^{n} \cdot \partial^{\beta^+}\varphi  - \lambda m^2\partial^\beta\partial_1 \eta^{n} \cdot \partial_1 \varphi -\partial^\beta (g\bar{\rho}\mathbf{e}_2+ a\bar{\rho} u^{n} ) \cdot  \varphi )
\mm{d}y\mm{d}\tau. \nonumber
\end{align}
Thus we further derive from the above two identities, \eqref{202008131740n} and \eqref{202008141600}  that
\begin{align}
\nabla^3\partial_1\eta^{n}  \mbox{ and }\nabla^3u^{n} \mbox{ are uniformly continuous in }H^{-1},
\label{20210302154safdsadf0}
\end{align}
where $n\geqslant 2$.
Next we further prove that
$\{(\eta^{n}, u^{n}, Q^{n})\}_{n=1}^{\infty}$
is a Cauchy sequence.

From now on, we always assume $n\geqslant 2$.
We define that
$$\bar{\eta}^{n+1}:={\eta}^{n+1}-{\eta}^{n},\quad
\bar{u}^{n+1}:={u}^{n+1}-{u}^{n},\quad
\bar{Q}^{n+1}:={Q}^{n+1}-{Q}^{n}\mbox{ and }\bar{\mathcal{A}}^{n}:={\mathcal{A}}^{n}-{\mathcal{A}}^{n-1}.$$
Then it follows from \eqref{01dsaf16asdfasf00002} that
\begin{equation}\label{01dsaf16asdfasf00003}
\begin{cases}
 \bar{\eta}^{n+1}_t =\bar{u}^{n+1} ,\\[1mm]
\bar{\rho}\bar{u}^{n+1} _t+\nabla_{\mathcal{A}^{n}}\bar{Q}^{n+1}+a\bar{\rho} \bar{u}^{n+1}=
\lambda m^2 \partial_1^2\bar{\eta}^{n+1}
-\nabla_{\bar{\mathcal{A}}^{n}}{Q}^{n} ,\\[1mm]
\div_{{\mathcal{A}}^{n}}\bar{u}^{n+1}=-\div_{{\bar{\mathcal{A}}}^{n}}{u}^{n}  , \\[1mm]
(\bar{\eta}^{n+1},\bar{u}^{n+1})|_{t=0}=(0,0) , \\[1mm]
(\bar{\eta}^{n+1},\bar{u}^{n+1})|_{\partial\Omega}\cdot\vec{\mathbf{n}}=0   .
\end{cases}
\end{equation}

It follows from \eqref{01dsaf16asdfasf00003}$_1$ and \eqref{01dsaf16asdfasf00003}$_4$ that, for $0\leqslant i\leqslant 4$,
\begin{align}
\label{202008140740}
\begin{aligned}
\sup\nolimits_{t\in  {I_T}}\|\bar{\eta}^{n+1}(t)\|_i+
\sqrt{\kappa}\|\partial_1\bar{\eta}^{n+1}\|_{L^2_TH^i}
\lesssim  {T}\|\bar{u}^{n+1}\|_{L^{\infty}_{T}H^i}.
\end{aligned}
\end{align}
By \eqref{202008131740n}, \eqref{202008141600}, \eqref{01dsaf16asdfasf00003}$_1$ and \eqref{202008140740}, we can estimate that, for $0\leqslant i\leqslant 3$,
\begin{align}
\label{202008140824}
&\sup\nolimits_{t\in  {I_{T }}} \|\bar{\mathcal{A}}^{n} \|_i
\lesssim T P(I^0)
 \|\bar{u}^{n}\|_{{L^{\infty}_{T}H^{i+1}}} ,\\
\label{2020081608}
&
\sup\nolimits_{t\in  {I_{T }}}\| \bar{\mathcal{A}}^{n}_t \|_i
\lesssim  P(I^0) \|\bar{u}^{n}\|_{{L^{\infty}_{T}H^{i+1}}}.
\end{align}

Let $\zeta^n =\eta^n+y$. $(\zeta^n)^{-1}$ denotes the inverse function of $\zeta^n$ with respect to the variable $y$. We define that \begin{align}
&K^5:= ( \lambda m^2 \partial_1^2\bar{\eta}^{n+1}
-\nabla_{\bar{\mathcal{A}}^{n}}{Q}^{n})/\bar{\rho} ,\
K^6:=   K^5- \bar{u}^{n+1}_t- a  \bar{u}^{n+1}  ,\nonumber \\
&(\beta,\varrho,\tilde{K}^5,\tilde{K}^6):=(\bar{Q}^{n+1},\bar{\rho},  {K}^5,K^6)|_{y=(\zeta^n)^{-1}(x,t)},\nonumber
 \end{align}then, by \eqref{01dsaf16asdfasf00003}$_2$,
\begin{equation}\nonumber
  \begin{cases}
\mm{div} \left(\nabla \beta/\varrho\right)
=\mm{div}\tilde{K}^6 &\mbox{in }\Omega,\\
  ( \nabla \beta /\varrho)\cdot\vec{\mathbf{n}}= \tilde{K}^5 \cdot\vec{\mathbf{n}}  &\hbox{on }\partial\Omega.
\end{cases}
\end{equation}
Applying the elliptic estimate \eqref{neumaasdfann1n} to the above boundary-value problem and  then making use of \eqref{2021sfa04031901} and \eqref{2022104101908}, we have
\begin{align}
\|\bar{Q}^{n+1} \|_{2}\lesssim & P(\|\eta^n\|_3) \|\beta \|_{2}
 \lesssim P(\|\eta^n\|_3)( \|   \tilde{K}^5\|_1+\|\mm{div}\tilde{K}^6\|_1)\nonumber \\
\lesssim   &P(\|\eta^n\|_3)(\|{K}^5\|_1+ \| \mm{div}_{\mathcal{A}}{K}^6\|_{0}).
\label{2sfa02008141000}
\end{align}

Noting that
$$
\mm{div}_{\mathcal{A}^n}{K}^6=
\mm{div}_{\mathcal{A}}( (\lambda m^2 \partial_1^2\bar{\eta}^{n+1}
-\nabla_{\bar{\mathcal{A}}^{n}}{Q}^{n})/\bar{\rho} )+\div_{{\bar{\mathcal{A}}}^{n}}{u}^{n} +\div_{{\mathcal{A}}^{n}_t}\bar{u}^{n+1}+
\partial_t\div_{{\bar{\mathcal{A}}}^{n}}{u}^{n},$$
thus, by using  \eqref{202008131740n}, \eqref{202008141600}, \eqref{202008140824} and \eqref{2020081608}, we easily estimate that
$$\|{K}^5\|_1+ \| \mm{div}_{\mathcal{A}}{K}^6\|_0
\lesssim P(I^0)
 \|\left(\bar{u}^{n}, \bar{u}^{n+1}, \partial_1\bar{\eta}^{n+1}\right)\|_{{L^{\infty}_{T}H^2}}  .$$
Putting the above estimate into \eqref{2sfa02008141000} yields
\begin{align}\label{202008141000}
\begin{aligned}
\| \bar{Q}^{n+1}\|_{C^0(\overline{I_T},H^2)}
&\lesssim P(b) \|\left(\bar{u}^{n}, \bar{u}^{n+1}, \partial_1\bar{\eta}^{n+1}\right)\|_{{L^{\infty}_{T}H^2}} .
\end{aligned}
\end{align}

Similarly to the derivation of \eqref{202008121220}, we  have that
\begin{align}\label{202008141635}
\begin{aligned}
 \sup\nolimits_{t\in  {I_{T }}}\|(\bar{u}^{n+1},\partial_1\bar{\eta}^{n+1})(t)\|_0^2
 \lesssim T P(I^0)
 \|\left(\bar{u}^{n}, \bar{u}^{n+1}, \partial_1\bar{\eta}^{n+1}\right)\|_{{L^{\infty}_{T}H^2}}^2.
\end{aligned}
\end{align}

By \eqref{01dsaf16asdfasf00003}$_1$--\eqref{01dsaf16asdfasf00003}$_3$, we have
\begin{align}
&
\partial_t\left(\mm{div}_{{\mathcal{A}}^{n}}\partial_1\bar{\eta}^{n+1}\right)=\mm{div}_{ {\mathcal{A}}^{n}_t}\partial_1\bar{\eta}^{n+1}
-\partial_1\left(\mm{div}_{{\bar{\mathcal{A}}}^{n}} {u}^{n}\right)  -\mm{div}_{\partial_1{{\mathcal{A}}}^{n}}\bar{u}^{n+1}
,\label{202008141704}\\[1mm]
&
\partial_t\left(\mm{curl}_{{\mathcal{A}}^{n}}\left(\bar{\rho}\bar{u}^{n+1}\right)\right)
+a\mm{curl}_{{\mathcal{A}}^{n}}\left(\bar{\rho}\bar{u}^{n+1}\right)
\nonumber\\
& = {\lambda m^2}({\bar{\rho}}^{-1}\partial_1\left(\mm{curl}_{\mathcal{A}^{n}}\left(\bar{\rho}\partial_1\bar{\eta}^{n+1}\right)\right)
 +  \partial_1\left[\mm{curl}_{\mathcal{A}^{n}}, {\bar{\rho}}^{-1}\right](\bar{\rho}\partial_1 \bar{\eta}^{n+1} )
 \nonumber\\[1mm]
&\quad -\mm{curl}_{\partial_1\mathcal{A}^{n}}\partial_1 \bar{\eta}^{n+1}) +\mm{curl}_{ {\mathcal{A}}^{n}_t}\left(\bar{\rho}\bar{u}^{n+1}\right)  +\mm{cur}_{\bar{\mathcal{A}}^n}(\nabla_{\mathcal{A}^{n-1}}Q^n) .\label{202008141724}
\end{align}
Making use of \eqref{202008131740n}, \eqref{202008141600} and \eqref{202008140824}, we can follow the argument  of \eqref{202008122105} to derive from \eqref{01dsaf16asdfasf00003}$_3$, \eqref{202008141704} and \eqref{202008141724} that, for a.e. $t\in I_T$,
\begin{align}
& \|(\mm{div}\bar{u}^{n+1},\mm{curl}\bar{u}^{n+1},\mm{div}\partial_1\bar{\eta}^{n+1},
 \mm{curl}\partial_1\bar{\eta}^{n+1})(t)\|_1^2  \nonumber \\
&\lesssim (\|\eta^0\|_3+ T P(I^0) )  \|\left( \bar{u}^{n}, \partial_1\bar{\eta}^{n}, \bar{u}^{n+1}, \partial_1\bar{\eta}^{n+1}\right)\|_{{L^{\infty}_{T}H^2}}^2
 + \|\left(\bar{u}^{n+1},  \partial_1\bar{\eta}^{n+1}\right)\|_1^2  . \label{202008142035}
\end{align}

 Consequently, summing up the estimates \eqref{202008140740},  \eqref{202008141635} and \eqref{202008142035} and then using  the interpolation inequality \eqref{201807291850}, Hodge-type elliptic estimate \eqref{202005021302} and Young's inequality, we obtain, for sufficiently small $\delta$,
\begin{align}
&\sup\nolimits_{t\in {I_{T}}}  \|\left(\bar{\eta}^{n+1}, \partial_1\bar{\eta}^{n+1},\bar{u}^{n+1} \right)\|_2^2  \lesssim T P(b)  \|\left(\bar{u}^{n} ,\partial_1\bar{\eta}^{n}, \bar{u}^{n+1}, \partial_1\bar{\eta}^{n+1}\right)\|_{L^2_TH^2}
.\label{202008142102}
\end{align}
In addition, by \eqref{202008131740n}, \eqref{202008141600} and \eqref{01dsaf16asdfasf00003}$_2$, we get that
\begin{align}\label{202008151041}
 \|  \bar{u}^{n+1}_t\|_{C^0(\overline{I_T},H^1)} \lesssim P(b)
\left(\|\left(\bar{u}^{n+1}, \partial_1\bar{\eta}^{n+1}\right)\|_{{C^0(\overline{I_T},H^2)}}+
\|  \bar{Q}^{n+1}\|_{{C^0(\overline{I_T},H^2)}}\right).
\end{align}

By \eqref{202008141000}, \eqref{202008142102} and \eqref{202008151041}, we immediately see that the sequence $\{( \eta^{n}, u^{n},{u}^{n}_t,Q^{n})\}_{n=1}^{\infty}$ is a Cauchy sequence in $C^0(\overline{I_T},H^{1,2}\times  H^2\times H^1\times\underline{H}^2)$  for sufficiently small $T\in (0,T_2]$, where the smallness depends  on $b$, $g$, $a$,  $\lambda$, $m$, $\bar{\rho}$ and $\Omega$. Hence
$$( \eta^{n},u^{n},u_t^n, Q^{n})\to (\eta,u,u_t,Q)
\mbox{ strongly in } C^0(\overline{I_T},H^{1,2}\times  H^2\times H^1\times \underline{H}^2).$$
Thanks to the strong convergence of $(\eta^{n},u^{n},u_t^n,Q^{n})$, we can take the limit in \eqref{01dsaf16asdfasf00002}, and then get a local classical solution $(\eta,u,Q)$ to the problem \eqref{01dsaf16asdfasf0000}. Let $q=Q-\bar{P}(\zeta_2)-\lambda|\bar{M}|^2/2$, then  $(\eta,u,q)$ is also a local classical solution to the transformed MRT problem
\eqref{01dsaf16asdfasf} by \eqref{dstist01} and \eqref{202009130836}.
In addition, thanks to \eqref{202saf008131740n}, \eqref{202008131740n} and \eqref{20210302154safdsadf0}, we easily follow the  compactness argument as in the proof of Proposition \ref{thm08} to further derive that
\begin{align}
& (\eta,u, q)
\in    \mathfrak{C}^0( \overline{I_{T}},{H}^{1,4}_{\mm{s}} )\times  \mathfrak{U}_{T}^4\times \mathfrak{Q}_{T}^4,\nonumber  \\
&    \sup\nolimits_{t\in {I_{T }}} \| \eta^{n} (t)\|_4^2 \leqslant 2\delta\leqslant \iota  \label{202201002042414107} ,\\
&\sup\nolimits_{t\in \overline{I_T}}\|  (\eta , \partial_1\eta,u  )\|_4^2 \leqslant c_4 I^0. \nonumber
\end{align}
In addition, it is easy to check that $(\eta,u,q)$ is  the unique   solution of the transformed MRT problem \eqref{01dsaf16asdfasf} in $\mathfrak{C}^0( \overline{I_{T}},{H}^{1,4}_{\mm{s}} )\times  \mathfrak{U}_{T}^4\times \mathfrak{Q}_{T}^4$ due to \eqref{202201002042414107}.  This completes the proof of Proposition \ref{202102182115}.

\appendix
\section{Analytic tools}\label{sec:09}
\renewcommand\thesection{A}
This appendix is devoted to providing some mathematical results, which have been use in previous sections.
 It should be noted that $\Omega$ appearing in what follows is  still defined by \eqref{0101a} and we will also use the simplified notations appearing in Section \ref{subsec:04}. In addition, $\Omega_T:=\Omega\times I_T$ and $a\lesssim b$ still denotes $a\leqslant cb$, however the positive constant $c$ depends on the parameters and domain in the lemma, in which $c$ appears.
\begin{lem}\label{201806171834}
\begin{enumerate}[(1)]
 \item  Embedding inequality (see \cite[4.12 Theorem]{ARAJJFF}): Let $D\subset \mathbb{R}^2$ be a domain satisfying the cone condition, then
\begin{align}
&\label{esmmdforinfty}\|f\|_{C^0(\overline{D})}= \|f\|_{L^\infty(D)}\lesssim\| f\|_{H^2(D)}.
\end{align}
 \item Interpolation inequality in $H^j$
 (see \cite[5.2 Theorem]{ARAJJFF}): Let $D$ be a domain in $\mathbb{R}^2$ satisfying the cone condition, then, for any given $0\leqslant j< i$,
\begin{equation}\label{201807291850}
\|f\|_{H^j(D)}
\lesssim  \|f\|_{L^2(D)}^{(i-j)/i}\|f\|_{H^i(D)}^{j/i} \lesssim  \varepsilon^{-j/(i-j)}\|f\|_{L^2(D)} + \varepsilon \|f\|_{H^i(D)} ,
\end{equation}
where the two estimate constants in \eqref{201807291850} are independent of $\varepsilon$, the positive constant $\varepsilon$ is arbitrary, and we have used Young's inequality in the last inequality above.
\item
Product estimates (see Section 4.1 in \cite{JFJSNS}): Let $D\in \mathbb{R}^2$ be a domain satisfying the cone condition, and the functions $f$, $g$ are defined in $D$. Then
\begin{align}
\label{fgestims}
&
 \|fg\|_{H^i(D)}\lesssim    \begin{cases}
                      \|f\|_{H^1(D)}\|g\|_{H^1(D)} & \hbox{ for }i=0;  \\
  \|f\|_{H^i(D)}\|g\|_{H^2(D)} & \hbox{ for }0\leqslant i\leqslant 2.
                    \end{cases}
\end{align}
                    \end{enumerate}
\end{lem}

\begin{lem}\label{10220830p}
Friedrich's inequality (see \cite[Lemma 1.42]{NASII04}): Let $1\leqslant p<\infty$, $n\geqslant 2$ and
$D\subset \mathbb{R}^n$ be a bounded Lipchitz domain.
Let a set $\Gamma\subset\partial D$ be measurable with respect to the $(n-1)$-dimensional measure $\tilde{\mu}:=\mathrm{meas}_{n-1}$
defined on $\partial D$ and let $\mathrm{meas}_{n-1}(\Gamma)>0$. Then
\begin{equation}
\label{poincare1}
\|w\|_{W^{1,p}(D)}\lesssim \|\nabla w\|_{L^p(D)}
\end{equation}
  for any  $w\in W^{1,p}(D)$ with $u\big|_{\Gamma}=0$ in the sense of trace.
\end{lem}
\begin{rem}\label{10220sasfasaf830p}
By Lemma \ref{10220830p}, we easily deduce that
\begin{equation}
\nonumber
\|w\|_{W^{1,p}(0,a)}\lesssim \| w'\|_{L^p(0,a)}
\end{equation}
 for any $w\in W^{1,p}(0,a)$ with
$w(0)=0$ or $w(a)=0$. Thus we further have
\begin{equation}
\label{poinsafdcasfaressadf1}
\|\varpi\|_{0}\lesssim \| \partial_2 \varpi\|_{0}\mbox{ for any }\varpi\in H_0^1:=\{\upsilon\in
H^1~|~\upsilon|_{\partial\Omega}=0\}.
\end{equation}
\end{rem}
\begin{lem}\label{10220830}
Poincar\'e inequality (see \cite[Lemma 1.43]{NASII04}): Let $1\leqslant p<\infty$, and $D$ be a bounded Lipchitz domain in $\mathbb{R}^n$ for $n\geqslant 2$ or a finite interval in $\mathbb{R}$. Then for any $w\in W^{1,p}(D)$,
\begin{equation}
\label{poincare}
\|w\|_{L^p(D)}\lesssim \|\nabla w\|_{L^p(D)}^p+\left|\int_{D}w\mathrm{d}y\right|^p.
\end{equation}
\end{lem}
\begin{rem}\label{10220saf830p}
By Poincar\'e inequality, we have, for any given $i\geqslant  0$,
\begin{align}
&\label{202012241002}
\| w\|_{1,i}\lesssim \|w\|_{2,i}\mbox{ for any  } w \mbox{ satisfying }\partial_1w,\ \partial_1^2w\in H^i.
\end{align}
 \end{rem}
\begin{rem}\label{10220saasdfsfasaf830p}
By Poincar\'e inequality, we also have
\begin{equation}
\label{poinsafdcaasdsfaressadf1}
\|\varpi\|_{0}\lesssim \| \partial_1 \varpi\|_{0}\mbox{ for any }\varpi\in H^1 \mbox{ satisfying }\varpi(y_1,y_2)=-\varpi(-y_1,y_2).
\end{equation}
\end{rem}
\begin{lem}\label{pro4a}
Hodge-type elliptic estimates:
If $w\in H^i_{\mm{s}}$ with $i\geqslant1$,
then
\begin{align}
&\label{202005021302}
\|\nabla w\|_{i-1}
\lesssim\|(\mm{div}w,\mm{curl}w)\|_{i-1}.
\end{align}
\end{lem}
 \begin{pf}
By a regularity method, we can verify that, for $\leqslant j\leqslant i-1$,
\begin{align}\label{202008152042}
\|\partial_1^j\nabla w\|_{0}^2
=\|\partial_1^j\mm{div}w\|_{0}^2+\|\partial_1^j\mm{curl}w\|_{0}^2,
\end{align}
which yields \eqref{202005021302} for $i=1$.

Next we further consider the case $i\geqslant 2$.
Since
\begin{align}\label{202008152054}
 \Delta w= \nabla\mm{div}w+\nabla^{\perp}\mm{curl}w,
\end{align}
where $\nabla^{\perp}:=(-\partial_2, \partial_1)^{\mm{T}}$,  we derive from \eqref{202008152042} and \eqref{202008152054} that, for any $1\leqslant l+k\leqslant i-1$,
$1\leqslant k$,
\begin{align}
&\|\partial_2^{k+1}\partial_1^{l}w\|_{0}
=\|\partial_2^{k-1}\partial_1^{l}\left(\Delta w-\partial_1^2w\right)\|_{0}\nonumber \\[2mm]
&\leqslant \|\partial_2^{k-1}\partial_1^{l}(\nabla\mm{div}w,\nabla^{\bot}\mm{curl} w)\|_{0}
+\|\partial_2^{k-1}\partial_1^{l+2} w\|_{0}, \nonumber
\end{align}
which further yields that
\begin{align}\label{202008152104}
\begin{aligned}
\|\partial_2^{k+1}\partial_1^lw\|_{0}
\lesssim
\|(\mm{div}w,\mm{curl}w)\|_{i-1}+\|\partial_2^{k-1}\partial_1^{l+2} w\|_{0}.
\end{aligned}
\end{align}

By an induction method, we easily derive from \eqref{202008152042} and \eqref{202008152104} that
\begin{align}\nonumber
 \|\partial_2^{k+1}\partial_1^{l}w\|_{0}
\lesssim\|(\mm{div}w, \mm{curl}w)\|_{i-1},
\end{align}
which, together with \eqref{202008152042}, yields
\begin{align}\nonumber
\|\nabla w\|_{i-1}\lesssim\|(\mm{div}w,\mm{curl}w)\|_{ i-1} .
\end{align}
This completes the proof of  Lemma  \ref{pro4a}.
\hfill$\Box$
\end{pf}

\begin{lem} \label{20201003302206}
Extension theorem: Let $i\geqslant 0$, $h>0$, $\delta=h/(i+1)$, and $f\in H^i$, then there exists a extension operator $\mathbb{E}_{\delta}^i$ such that $\mathbb{E}_{\delta}^i$: $f\in H^i\to  H^i(\mathbb{T}\times \mathbb{R})$ such that
\begin{align}
&\mathbb{E}_{\delta}^i(f)=0\mbox{ for } y_2<-\delta/2,\ h+\delta/2<y_2 , \label{20221042110}  \\
&
\mathbb{E}_{\delta}^i(f)|_{\Omega}=f\mbox{ and }\| \mathbb{E}_{\delta}^i(f)\|_{H^i(\mathbb{T}\times \mathbb{R})}\lesssim  \| {f}\|_i . \label{20221042110x}
\end{align}

\end{lem}
\begin{pf}
Let
\begin{align}
\Omega_-:=\mathbb{T}\times (-\infty, 0), \
\Omega_+:=\mathbb{T}\times(h,+\infty ).
\label{2022104122056}
\end{align}
Let $\chi\in C_0^\infty[0, \delta )$ with $\chi(y_2)=1$ for $y_2\in (0,\delta/2)$ and
\begin{align}\nonumber
\tilde{f} =
\begin{cases}
\chi(-y_2)\sum_{j=1}^{i+1}\lambda_jf(y_1,-jy_2)  &\mbox{ for }y \in \Omega_-,  \\
f &\mbox{ for }y \in \Omega , \\
0 &\mbox{ for }y \in \partial\Omega , \\
\chi(y_2-h)  \sum_{j=1}^{i+1 }\lambda_jf(y_1,h-j(y_2-h)) &  \mbox{ for }y \in \Omega_+,
\end{cases}
\end{align}
where  $\sum_{j=1}^{i+1}(-j)^k\lambda_j =1$  for $0\leqslant k\leqslant i$.

By the definition of $H^i$  and the facts, for $i\geqslant 1$,
\begin{align}
 \nabla^{i-1}\tilde{f}|_{\Omega_-}=\nabla^{i-1}\tilde{f}|_{\Omega}\mbox{ on }\mathbb{T}\times \{0\},\  \nabla^{i-1}\tilde{f}|_{\Omega_+}=  \nabla^{i-1}\tilde{f}|_{\Omega}\mbox{ on }\mathbb{T}\times \{h\}  \nonumber
\end{align}
in the sense of trace by a density argument \cite[5.19 Theorem]{ARAJJFF}, we can easily check that $\tilde{f}$ constructed by above belongs to $H^i(\mathbb{T}\times \mathbb{R})$ and satisfies \eqref{20221042110} and \eqref{20221042110x}. Consequently, Lemma \ref{20201003302206} holds.
This completes the proof.
\end{pf}
\begin{lem}\label{20180812}
Dual estimates: Let $\tau$ satisfying $|\tau|\in (0,1) $, $ {\varphi}$, $ {\psi}\in H^1$, and $D_1^\tau {\varphi}=({\varphi}(y_1+\tau,y_2)- {\varphi}(y_1,y_2))/\tau$. Then
\begin{align}
\label{201808121426}
&\left|\sum_{s=0}^h\int_0^{2\pi}(D_1^\tau {\varphi} {\psi})|_{y_2=s}\mm{d}y_1\right|\lesssim \|\varphi\|_{1} \|\psi\|_{1}.
\end{align}
\end{lem}
\begin{pf}
Let $f$, $\omega \in \mathcal{Y}:=\{w\in H^1 \cap C^1(\overline{\Omega})~|~w(y_1,0)=0\}$. Defining the Fourier coefficient $f$ by $\widehat{f}(\xi,y_2):= \int_0^{2\pi}f(s,y_2)e^{-i\xi s}\mm{d}s$, then, using Parseval's relation on the Torus (see \cite[Proposition 3.1.16]{grafakos2008classical}), Plancherel's identity,  H\"older's inequality in $l^2$, the horizontal periodicity,
and Newton--Leibnitz's formula, we have, for any $y_2\in (0,h)$,
 \begin{align}
&\left|\int_0^{2\pi}D_1^\tau f\omega (y_1,y_2)\mm{d}y_1\right|=\frac{1}{2\pi}\left|  \sum_{\xi\in \mathbb{Z}}\widehat{D_1^\tau f} \overline{\widehat{ { {\omega }}}} (y_2)\right| =\frac{1}{2\pi}\left|\tau^{-1}\sum_{\xi\in \mathbb{Z}}   (e^{i\xi \tau}-1)\widehat{  f} \overline{\widehat{ { {\omega }}}} (y_2)\right|\nonumber \\
&\lesssim  \sum_{\xi\in \mathbb{Z}}  |\xi|\sqrt{\|\widehat{  f}\|_{L^2(0,\tau)}\|\partial_2\widehat{  f}\|_{L^2(0,\tau)}\|
 \overline{\widehat{ { {\omega }}}}\|_{L^2(0,\tau)}\|
 \partial_2\overline{\widehat{ { {\omega }}}}\|_{L^2(0,\tau)}} \nonumber \\
&\lesssim  \left({\sum_{\xi\in \mathbb{Z}}(\xi\|\widehat{  f}\|_{L^2(0,\tau)})^2 \sum_{\xi\in \mathbb{Z}}\|\partial_2\widehat{  f}\|_{L^2(0,\tau)}^2 \sum_{\xi\in \mathbb{Z}} (\xi\|
 \overline{\widehat{ { {\omega }}}}\|_{L^2(0,\tau)})^2 \sum_{\xi\in \mathbb{Z}}\|
 \partial_2\overline{\widehat{ { {\omega }}}}\|_{L^2(0,\tau)}^2}\right)^{1/4}\nonumber \\
&\lesssim   \|f\|_{1} \|\omega \|_{1}.\label{201808121247ss}
\end{align}

Let the extension operator $\mathbb{E} _{\delta}^1 $ be defined in Lemma \ref{20201003302206} with $i=1$ and $\delta=h/2$. For simplicity, we denote  $\mathbb{E} _{\delta}^1 ({\chi})$  by $\tilde{\chi}$, where $\chi$ represents $\varphi$ or $\psi$. 
Let $\varepsilon\in (0,1)$,
then we denote by  $S_\varepsilon(\tilde{\chi})$  the mollification of  $\tilde{f}$ with respect to the 2D  variable $(y_1,y_2)$. Exploiting trace theorem and the properties of mollification, we have, for any given $\tau\in [0,1)$,
\begin{align}
\label{20220103231619}
\|(S_\varepsilon(\tilde{\chi})  -\chi)(y_1+\tau,y_2) \|_{H^{1/2}(0,2\pi)}  \lesssim \|S_\varepsilon(\tilde{\chi}) - {\chi}\|_1\to 0\mbox{ as }\varepsilon\to 0.
\end{align}

Let $\sigma( y_2)\in C_0^\infty(0,h] $  satisfy $\sigma(h)=1$, and $\tilde{\chi}_\sigma^\varepsilon:= S_\varepsilon(\tilde{\chi}) \sigma( y_2)$.
  Then $\chi_\sigma^\varepsilon\in \mathcal{Y}$.
By \eqref{201808121247ss}, we have
\begin{align}
\left|\int_0^{2\pi}(D_1^{\tau}  (\tilde{\varphi}_\sigma^\varepsilon)|_{y_2=h} \tilde{\psi}_\sigma^\varepsilon )|_{y_2=h}\mm{d}y_1\right|  \lesssim \|\tilde{\varphi}_\sigma^\varepsilon\|_{1} \| \tilde{ \psi}_\sigma^\varepsilon\|_{1}\lesssim \|S_\varepsilon(\tilde{\varphi} )\|_{1} \|S_\varepsilon( \tilde{\psi} )\|_{1}. \nonumber
\end{align}
Thanks to  \eqref{20220103231619}, we easily derive from the above estimate that  \begin{align}
\label{201808121426xxx}
&\left| \int_0^{2\pi}(D_1^\tau {\varphi} {\psi})|_{y_2=h}\mm{d}y_1\right|\lesssim \|\varphi\|_1\|\psi\|_{1}.
\end{align}

Similarly, we also have
\begin{align}\nonumber
&\left| \int_0^{2\pi}(D_1^\tau {\varphi} {\psi})|_{y_2=0}\mm{d}y_1\right|\lesssim \|\varphi\|_{1} \|\psi\|_{1},
\end{align}
which, together with \eqref{201808121426xxx}, yields \eqref{201808121426}.
 \hfill $\Box$
\end{pf}

\begin{lem}
\label{lem:08181945}Elliptic estimates:  Let $k\geqslant 2$.
Assume $A$ is a $2\times 2$ symmetry matrix function, each element $A_{ij}$ of which belongs to  $W^{k,\infty}(\Omega)$, and positivity condition, i.e., there exists a positive constant $\theta$ such that $ (A \xi)\cdot \xi\geqslant  \theta |\xi|$ for a.e. $y\in \Omega$ and all $\xi\in \mathbb{R}^2$. If $f^1\in H^{k-2}$ and $f^2\in H^{k-1}$ satisfy the following compatibility condition
\begin{align}\label{202008181950n}
\int f^1\mm{d}y+\int_{\partial\Omega}f^2 \vec{\mathbf{n}}_2\mm{d}y_1=0,
\end{align}
where $\vec{\mathbf{n}}$ denotes the  outward unit normal vector on the boundary $\partial\Omega$,
  there exists a unique strong solution $p\in \underline{H}^k $ to
the Neumann boundary-value problem of elliptic equation:
\begin{align}\label{neumann}
\begin{cases}
 -\mm{div}\left(A\nabla p \right)=f^1 &\mbox{in } \Omega,\\[1mm]
( A\nabla p)\cdot\vec{\mathbf{n}}=f^2\vec{\mathbf{n}}_2  &\mbox{on } \partial\Omega.
\end{cases}
\end{align}
 Moreover $p$ satisfies
\begin{align}\label{neumaasdfann1n}
\|p\|_{k}\lesssim  (1+\|A\|_{W^{k-1,\infty}(\Omega)})(\|f^1\|_{k-2}+\|f^2\|_{k-1}),
\end{align}
\end{lem}

{\begin{pf}
  Noting that  $A$ is a  symmetry matrix function, we define an inner-product of $\underline{H}^{1} $ by
\begin{align*}
(\varphi,\phi)_{\underline{H}^{1}}:=\int  (A\nabla \varphi)\cdot \nabla\phi\mm{d}y\mbox{ for }\varphi,\ \phi
\in \underline{H}^{1} ,
\end{align*}
and then the corresponding norm by $\|\varphi\|_{\mathcal{H}}:=\sqrt{(\varphi,\varphi)_{\underline{H}^{1}}}$. Obviously, by Poincar\'e inequality \eqref{poincare} and the positivity condition, we have
\begin{align*}
\|\varphi\|_{1}\lesssim\|\varphi\|_{\mathcal{H}}\lesssim \|A\|_{L^\infty}\|\varphi\|_{1}.
\end{align*}

We define the functional
\begin{align}
\nonumber
F(\varphi):=\int f^1\varphi\mm{d}y+\int_{\partial\Omega} f^2\vec{\mathbf{n}}_2\varphi\mm{d}y_1\mbox{ for }\varphi\in \underline{H}^{1} .
\end{align}
Then it is easy to check that the functional $F$ is a bounded linear functional on $\underline{H}^1 $.
By Riesz representation theorem, there exists a unique $p\in\underline{H}^1$, such that
\begin{align}\label{202012191937}
(p,\varphi)_{\underline{H}^{1}}=F(\varphi)
\mbox{ for any }\varphi\in \underline{H}^{1}
\end{align}
and
\begin{align}\label{202012191940}
\|p\|_{1}\lesssim\|p\|_{\mathcal{H}} \lesssim \|f^1\|_{k-2}+\|f^2\|_{k-1} .
\end{align}

For any given $\psi\in {H}^{1} $, we denote $\varphi=\psi-(\psi)_{\Omega}$.
Then $\varphi\in \underline{H}^{1} $. Putting it in \eqref{202012191937}, we obtain
\begin{align}
\nonumber
\int (A\nabla   p)\cdot \nabla \psi\mm{d}y
=\int f^1\psi\mm{d}y+\int_{\partial\Omega} f^2 \vec{\mathbf{n}}_2\psi\mm{d}y_1
-(\psi)_{\Omega}\left(\int  f^1\mm{d}y+\int_{\partial\Omega} f^2\vec{\mathbf{n}}_2\mm{d}y_1\right),
\end{align}
which, together with the compatibility condition \eqref{202008181950n}, yields
\begin{align}
 \int (A\nabla  p)\cdot \nabla\psi\mm{d}y= \int  f^1\psi \mm{d}y   +\int_{\partial\Omega} f^2 \vec{\mathbf{n}}_2\psi\mm{d}y_1
 \mbox{ for any } \psi\in {H}^{1}.\label{202012191952}
\end{align}

Next we further improve the regularity of $p$.
\emph{We assert that
\begin{itemize}
  \item for any $0\leqslant l\leqslant k- 1$,
\begin{align}
&p\in \underline{H}^{l+1}, \ \|p\|_{l+1}\lesssim  (1+\|A\|_{W^{k-1,\infty}(\Omega)})(\|f^1\|_{k-2}+\|f^2\|_{k-1}), \label{20201022271158}
\end{align}
  \item for any $0\leqslant l\leqslant k- 2$,
\begin{align}
 \int (A\nabla  \partial_1^{l}p)\cdot \nabla\psi\mm{d}y=&\int  (\partial_1^{l}f^1\psi-
{ [\partial_1^{l},A ]\nabla p} )
\cdot \nabla \psi\mm{d}y \nonumber
\\
& +\int_{\partial\Omega} \partial_1^{l}f^2 \vec{\mathbf{n}}_2\psi\mm{d}y_1
 \mbox{ for any } \psi\in {H}^{1},\label{202sadfa012191952}
\end{align}
\end{itemize}}

Noting that  \eqref{20201022271158}--\eqref{202sadfa012191952} hold for $l=0$ due to  \eqref{202012191940} and \eqref{202012191952}, thus, to obtain the above assertion,  it suffices to prove that
\emph{if \eqref{20201022271158}--\eqref{202sadfa012191952} hold for $0\leqslant l\leqslant k-2$, then }
\begin{itemize}
  \item\emph{ \eqref{202sadfa012191952} holds with $l+1$ in place of $l$, if $l+1\leqslant k-2$.}
  \item\emph{ \eqref{20201022271158}  also holds with $l+1$ in place of $l$.  }
\end{itemize}
Next we  prove these two facts.

(1)
Let $D_1^{\tau}w= (w(y_1+{\tau},y_2)-w(y_1,y_2))/{\tau}$ with $|\tau|\in (0,1)$. Since \eqref{202sadfa012191952} hold for $0\leqslant l\leqslant k-2$, we can take
 $\psi =D_1^{-{\tau}}D_1^{{\tau}}\partial_1^{l}p$ in \eqref{202sadfa012191952} to get
\begin{align}
&\int (A \nabla \partial_1^{l} p) \cdot\nabla \left(D_1^{-{\tau}}D_1^{{\tau}}{ \partial_1^{l}}p\right)\mm{d}y\nonumber \\
&=\int \left(\partial_1^{l} f^1D_1^{-{\tau}}D_1^{{\tau}}{ \partial_1^{l}}p+
{ [\partial_1^{l},A ]\nabla p} )
\cdot \nabla D_1^{-{\tau}}D_1^{{\tau}}{\partial_1^{l}}p\right)\mm{d}y  +\int_{\partial\Omega} \partial_1^{l}f^2 \vec{\mathbf{n}}_2D_1^{-{\tau}}D_1^{{\tau}}{ \partial_1^{l}}p\mm{d}y_1.\label{202012192021}
\end{align}

By the properties of difference quotient, we derive from \eqref{202012192021}
that
\begin{align}\label{202012192014}
&\int (A (y_1+\tau) \nabla \left(D_1^{{\tau}}\partial_1^l p\right))\cdot \nabla D_1^{{\tau}}\partial_1^lp \mm{d}y
\nonumber \\
&=\int ((D_1^{{\tau}} {\left([\partial_1^{l},A]\nabla p\right)} -D_1^{\tau}A   \nabla \partial_1^l p)\cdot \nabla D_1^{{\tau}}\partial_1^lp
-   \partial_1^{l} f^1   D_1^{-{\tau}}D_1^{{\tau}}\partial_1^{l}p
)\mm{d}y\nonumber \\
&\quad +\int_{\partial\Omega} D_1^{{\tau}}\partial_1^{l}f^2 \vec{\mathbf{n}}_2\left(D_1^{{\tau}}\partial_1^{l}p\right)\mm{d}y_1=:I_{10} .
\end{align}

 Thanks to Lemma \ref{20180812}, we can estimate that
$$
\begin{aligned}
&I_{10} \lesssim
( \|A\|_{W^{k-1 ,\infty}(\Omega)}  \|\nabla \partial_1^{l}p \|_0  + \| \partial_1^{l}f^1\|_{0} + \| {  \partial_1^{l}} f^2\|_{1} ) \|  D_1^{\tau}\partial_1^{l}p\|_{1} .
\end{aligned}$$
Making use of the above estimate, Poincar\'e inequality \eqref{poincare}   and the assumption that  \eqref{20201022271158} holds for $0\leqslant l\leqslant k-2$, we can deduce from \eqref{202012192014} that
\begin{align}
\nonumber
\| D_1^{h}\partial_1^{l}p\|_{1} \lesssim  (1+\|A\|_{W^{k-1 ,\infty}(\Omega)})\left(\| f^1\|_{k-2}
+\|  f^2\|_{k-1}\right),
\end{align}
which implies
\begin{align}\label{202012192022}
\|  \partial_1^{l+1}p\|_{1}\lesssim(1+\|A\|_{W^{k-1 ,\infty}(\Omega)})\left(\| f^1\|_{k-2}
+\|  f^2\|_{k-1}\right).
\end{align}
  Thanks to the  regularity $\nabla \partial_1^{l+1} p\in L^2$, we further see from \eqref{202sadfa012191952} for $0\leqslant l\leqslant k-2$ that \eqref{202sadfa012191952} holds with $l+1$ in place of $l$, if $l+1\leqslant k-2$.

(2) Next we shall further derive that $ \partial_2^{m+1} \partial_1^{l+1-m}p\in L^2$ for $1\leqslant m\leqslant l+1$.  Since $p$ satisfies  \eqref{202012191952}, thus
\begin{align}\label{202012192036}
\int  A_{22}\partial_2 p\partial_2\psi\mm{d}y
=\int  \left(f^1 \psi
 -\sum_{i,j\neq 2} A_{ij}\partial_j p\partial_i\psi
\right)\mm{d}y
+\int_{\partial\Omega} f^2\vec{\mathbf{n}}_2\psi\mm{d}y_1 .
\end{align}

Let  $\varphi\in C_0^{\infty}(\Omega)$.  Noting that $A_{22} \geqslant \theta$, we can take $\psi=   A_{22}^{-1}\varphi$ in \eqref{202012192036} to get that
\begin{align}\label{202012192038}
\int \partial_2 p\partial_2\varphi\mm{d}y
=\int_{\Omega}A_{22}^{-1} \left(f^1+ \sum_{i,j\neq 2} \partial_i (A_{ij}\partial_j p)
 + \partial_2A_{22} \partial_2p \right)\varphi\mm{d}y ,
\end{align}
which implies
\begin{align}
\label{202102261331}
\partial_2^2p=-A_{22}^{-1} \left(f^1+ \sum_{i,j\neq 2} \partial_j (A_{ij}\partial_i p)
 + \partial_2A_{22} \partial_2p
\right).
\end{align}
Since we have assumed that \eqref{20201022271158} holds for  $0\leqslant l\leqslant k-2$, we easily see from
\eqref{202012192022} and \eqref{202102261331} that  \eqref{20201022271158}  also holds with $l+1$ in place of $l$.

Finally, to complete the proof of Lemma \ref{lem:08181945}, next it suffices to verify \eqref{neumann}.

By the relation \eqref{202102261331}, we can compute out that
\begin{align}\label{202012192057}
-\mm{div}\left(A\nabla p \right)= f^1\mbox{ a.e. in } \Omega.
\end{align}
Multiplying the above identity \eqref{202012192057}  by
$\psi\in {H}^{1}$ in $L^2$, and then using the integral by parts, we get
\begin{align}\label{202012192101}
-  \int \left(A\nabla p \right)\cdot\nabla \psi\mm{d}y
+\int_{\partial\Omega} (A\nabla p )\cdot\vec{\mathbf{n}}\psi\mm{d}y_1
=\int  f^1\psi\mm{d}y.
\end{align}
Comparing \eqref{202012191952} with \eqref{202012192101}, we further have
\begin{align}\nonumber
\int_{\partial\Omega} \left((A\nabla p) -f^2\right)\cdot\vec{\mathbf{n}}\psi\mm{d}y_1=0
 \mbox{ for any }\psi\in {H}^{1},
\end{align}
which implies
\begin{align}\nonumber
(A\nabla p)\cdot\vec{\mathbf{n}}= f^2 \vec{\mathbf{n}}_2\mbox{  on }\partial\Omega
\end{align}in the sense of trace.
Hence $p$ solves problem \eqref{neumann}. The proof of Lemma \ref{lem:08181945} is complete.
\hfill $\Box$
\end{pf}
}

\begin{lem}\label{pro:1221}Diffeomorphism mapping theorem:
There exists a sufficiently small constant $\iota\in(0,1]$, depending on $\Omega$, such that,
for any ${\varsigma}\in H_{\mm{s}}^3$
satisfying  $\|{\varsigma}\|_4\leqslant \iota$,
$\psi:=\varsigma+y$ (after possibly being redefined on a set of measure zero with respect to variable $y$) satisfies the diffeomorphism  properties  as $\zeta$ in \eqref{20210301715x}, \eqref{20210301715} and  $\inf_{y\in \Omega}\det(\nabla \varsigma  +I)\geqslant 1/4$.
\end{lem}
{\begin{pf}
In  \cite{JFJSARMA2019}, the authors had proved that Lemma \ref{pro:1221} holds for ${\varsigma}\in H_0^4$, where $  H_0^4:=\{w\in H^4~|~w|_{\partial\Omega}=0\}$, see  \cite[Lemma 4.2]{JFJSARMA2019}. Next we further verify Lemma \ref{pro:1221} based on  \cite[Lemma 4.2]{JFJSARMA2019}.

Let $\delta=h/5$  and  $\Omega_\delta=\mathbb{T}\times \delta$.
By Lemma \ref{20201003302206}, there exists an extension function $\tilde{{\varsigma}}$ such that
$$\tilde{\varsigma}|_{\Omega}=\varsigma,\ \tilde{\varsigma}|_{\partial\Omega_\delta} =0\mbox{ and }\| \tilde{\varsigma}\|_{H^4( \Omega_{\delta})}\lesssim  \| {\varsigma}\|_4.$$

Thus, thanks to \cite[Lemma 4.2]{JFJSARMA2019}, we have
\begin{align}\label{2020xx12211750}
&\tilde{\psi}:=\tilde{{\varsigma}}+y:\ \overline{\Omega_\delta}\to \overline{\Omega_\delta}  \mbox{ is a }C^2\mbox{-diffeomorphism mapping},\\
&  \tilde{\psi}\mbox{ is an identity map on the boundary }\partial\Omega_\delta. \nonumber
\end{align}

Noting that $\varsigma=\tilde{{\varsigma}}|_{\overline{\Omega}}$, thus
\begin{align}\nonumber
\psi :={\varsigma}+y:\overline{\Omega}\to \overline{\Omega} \mbox{ is injective}.
\end{align}
To complete the proof of Lemma \ref{pro:1221}, next it suffices to verify that
\begin{align}\nonumber
\psi  :\ \overline{\Omega}\to\overline{\Omega} \mbox{ is surjective}.
\end{align}

To this purpose, we let $x^0\in \overline{\Omega} $, then there exists $y^0\in\overline{\Omega_\delta} $ such that
$\tilde{\psi}(y^0)=x^0$ by  \eqref{2020xx12211750}.
Since ${\varsigma}_2|_{\partial\Omega}=0$,  it is easy to see that, for sufficiently small $\iota$,
\begin{align}\label{202012211907n0}
\tilde{\psi} :={\varsigma}+y:\Xi_i\rightarrow\Xi_i\mbox{ is bijective},
\end{align}
where $\Xi_i:=\mathbb{T}\times\{i\}$ and $i=0$, $h$. By \eqref{2020xx12211750} and \eqref{202012211907n0}, we further see that
\begin{align}\nonumber
y^0\in {\Omega_\delta}\backslash \partial\Omega .
\end{align}
Now we claim that
\begin{align}
\label{20201022251957}
y^0\in \Omega .
\end{align}

In fact, if \eqref{20201022251957} fails, then  $y^0\in  \Omega_- $ or $\Omega_+$, see \eqref{2022104122056} for the definition $\Omega_- $ and $\Omega_+$. We assume that $y^0\in \Omega_-$, then
$\tilde{\psi}$ maps the segment $l_-:=\{y\in \Omega_- ~|~y_1=y^0_1,\ -\delta\leqslant y_2\leqslant y^0_2\}$ to a continuous curve, which lies in $ \Omega_\delta$. Noting that $\tilde{\psi}( y^0_1,-\delta)=(y^0_1,-\delta) $ and  $\tilde{\psi}(y^0_1, y^0_2)\in \Omega $, thus there exists a unique point $y^{c}\in l_-$, such that
 $\psi (y^{c})\in \Xi_0$, which contradicts with  \eqref{202012211907n0}. Hence $y^0 \notin  \Omega_-$.
 Similarly  $y^0 \notin  \Omega_+$. Consequently, \eqref{20201022251957} holds.
This completes the proof.
\hfill$\Box$
\end{pf}
}
\begin{lem}\label{20221033311357}
New version of theorem of continuity in the mean of integral:  Let $f\in C^0(\overline{I_T},L^2)$ and  $\varsigma  \in C^0(\overline{\Omega_T} )$ satisfies the diffeomorphism properties  as \eqref{20210301715x} and \eqref{20210301715}.  Then for any given $\varepsilon>0$ and for any given $s\in \overline{I_T}$, there exists a $\delta>0$ such that for any $t\in \overline{I_T}$ satisfying $|t-s|<\delta$,
\begin{align}
\nonumber
\| f(\varphi(y,t),s)-f( \varphi(y,s),s)\|_0\leqslant \varepsilon.
\end{align}
\end{lem}
\begin{pf}
It is well-know that, for any $\chi\in L^2(\mathbb{R}^n)$ with $n\geqslant 1$,
$$\lim_{h\to 0}\int_{\mathbb{R}^n} |\chi(x+h)-\chi(x)|\mm{d}x=0, $$
see the theorem of continuity in the mean of integral in \cite[Theorem 4.21]{2020221032301238}.

Recalling the proof of the above assertion, we easily   see that, for any $\chi\in L^2(\mathbb{R}^n)\cap L^\infty(\mathbb{R}^n)$ with $n\geqslant 1$,
\begin{align}
\label{20210303311249}
\lim_{h\to 0}\int_{\mathbb{R}^n} |\chi(x+h)-\chi(x)|^2\mm{d}x=0.
\end{align}

Then following the argument of \eqref{20210303311249} with further using the assumptions of $\varsigma$ and the embedding inequality \eqref{esmmdforinfty}, we easily get desired conclusion in Lemma \ref{20221033311357}.
\hfill$\Box$
\end{pf}

\begin{lem}\label{20021032019018} Some results for functions with values in Banach spaces:
  Let $T>0$,  integers $i$, $j\geqslant 0$ be given and  $1\leqslant p\leqslant \infty$.
\begin{enumerate}
\item[(1)] Assume $f\in L^p_TH^{i}$, $\partial_1^k\nabla^{i}f\in L^p_TL^2$ for any $1\leqslant k\leqslant j$, then $f\in L^p_TH^{j,i}$, where
$H^{j,i}:=\{w\in H^{ i}~|~\partial_1^k w\in H^i\mbox{ for any }1\leqslant k\leqslant j\}$.
\item[(2)] Let $X$ be a separable Banach space and $T>0$.
If $w\in C^0(\overline{I_T},X)$, then  $w$: $I_T\to X$ is a strongly measurable function and
\begin{align}
\|w\|_X\in C^0(\overline{I_T}).  \label{2022104121338}
\end{align}
\item[(3)] Let $  j\geqslant i+1$ and $c$ be a constant. If $f\in L^p(I_T,H^i)$ with $1<p\leqslant \infty$,  $\|f\|_{j}\leqslant c \|g\|_j$ holds for a.e. $t\in I_T$ and $g\in L^p(I_T,H^j)$, then
\begin{align}
f\in L^p_TH^i.
\label{202104032134}
\end{align}
\item[(4)] Assume $f\in L^\infty_TH^{i}$ and $\{D^\tau_1f\}_{|\tau|\in (0,1)}$ is uniformly bounded in $L^\infty_TH^{i}$.
\begin{enumerate}
     \item[(a)]  There exists a subsequence (still denoted by $\{D^\tau_1f\}_{|\tau|\in (0,1)}$ ) of $\{D^\tau_1f\}_{|\tau|\in (0,1)}$ such that  \begin{align}
D_1^\tau f \rightharpoonup \partial_1f  \mbox{ weakly-* in }L^\infty_TH^i\mbox{ as }\tau\to 0.\label{202104141653asfda}
\end{align}
Moreover, $f \in L^\infty_TH^{1,i}$.
  \item[(b)] If additionally  $f\in L^p_TH^{i}$  and $\{D^{-\tau}_1D^\tau_1f\}_{|\tau|\in (0,1)}$  is uniformly bounded in $L^p_TH^{i}$ with $1<p<\infty$, then (a subsequence)
   \begin{align}D_1^{-\tau}D_1^\tau(f) \rightharpoonup \partial_1^2 f \mbox{ weakly in }L^p_TH^i\mbox{ as }\tau\to 0 \mbox{ (a subsequence)}.
\label{202104141653} \end{align}
Moreover, $f \in L^\infty_TH^{2,i}$.  \item[(c)]  If additionally $D_1^\tau\partial^\alpha f$ in $C^0(\overline{I_T},L^2_{\mm{weak}})$ and  $D_1^\tau\partial^\alpha f$ is uniformly continuous in $H^{-1}$, where  the multiindex $\alpha$ satisfying $|\alpha|=i$, then
\begin{align}
D_1^\tau\partial^\alpha f \to \partial^\alpha \partial_1 f\mbox{ in } C^0(\overline{I_T},L^2_{\mm{weak}}) \mbox{ (a subsequence)}. \label{202104adsa141653}
\end{align}
   \end{enumerate}
\item[(5)] Let $\varphi =\varsigma+y$ and $\varsigma\in \mathbb{A}_{T,\iota}^{4,1/4}$ defined by \eqref{2022104161633}.
\begin{enumerate}
\item[(a)] If $f\in C^0(\overline{I_T},H^i)$ or $f\in L^p_TH^i$ with $0\leqslant i\leqslant 4$,  then
 \begin{align}
\label{2020103250855}
&F:=f(\varphi ,t)\in C^0(\overline{I_T},H^i)\mbox{ or }L^p_TH^i
\end{align}
and
\begin{align} \mathcal{F}:=f(\varphi^{-1},t)\in C^0(\overline{I_T},H^i)\mbox{ or } L^p_TH^i.
\label{202104031901} \end{align}
Moreover,
\begin{align}
 & \|F\|_{L^p_TH^i}\lesssim
P(\|\varsigma\|_{L^\infty_TH^4}) \|f \|_{L^p_TH^i},
\label{2021sfa04031901} \\
&  \|\mathcal{F}\|_{L^p_TH^i}\lesssim
P(\|\varsigma\|_{L^\infty_TH^4})\|f \|_{L^p_TH^i}. \label{2022104101908} \end{align}
 \item[(b)]  If $\varsigma$ additionally satisfies $ \varsigma_t \in L^\infty_TH^3$,  then, for any $f\in L^p(I_T,H^i)$ satisfying $f_t\in L^p(I_T,H^{i-1})$ with $1\leqslant i\leqslant 4 $, \begin{align}F_t=(f_t(x,t) +\varsigma_t\cdot \nabla f(y,t))|_{x=\varphi}\in  L^p_TH^{i-1}
\label{202104032132}
\end{align}
and
\begin{align}\mathcal{F}_t=(f_t(y,t) -(\nabla \varphi)^{-1}\varsigma_t \cdot \nabla f(y,t))|_{y=\varphi^{-1}}\in  L^p_TH^{i-1}.
\label{2021040sdaf32132}
\end{align}
\end{enumerate}
\end{enumerate}
\end{lem}
\begin{pf}
(1) Since $f$ is a Bochner integrable by $f\in L^p_TH^{i}$, thus there exists a sequence $\{f^n\}_{n=1}^\infty$ of simple functions such that
\begin{align}
\label{201302201832}
\lim_{n\to \infty}\int_0^T\|f^n(t)-f(t)\|_{ i}\mm{d}t=0.
\end{align}

Similarly to \eqref{202121020281617}, we define the mollifications of $f^n$, resp. $f$ with respect to  $y_1$ as follows:
 \begin{align}\nonumber
 S_\varepsilon^1 (f^n):=\chi^{\varepsilon} * f^n,\mbox{ resp. }S_{ \varepsilon}^1(f):=\chi^{\varepsilon} * f,
\end{align}
where $\varepsilon\in (0,1)$.
Then $S_{\varepsilon}^1 (f^n)$, $S_{\varepsilon}^1 (f)\in L^p_TH^{j,i}$. By \eqref{201302201832}, we have, for  given $\varepsilon $,
\begin{align}
\label{202103201834}
\lim_{n\to \infty}\int_0^T\| S_{\varepsilon}^1(f^n(t)-f(t))\|_{H^{j,i}}\mm{d}t=0.
\end{align}
By the regularity of $f$, we see that
\begin{align}\nonumber
 \int_0^T\|  f (t) \|_{H^{j,i}}\mm{d}t<\infty.
\end{align}
Thus
\begin{align}\label{202103201834x}
\lim_{\varepsilon\to 0}\int_0^T\| S_{\varepsilon}^1 (f(t))-f(t)\|_{H^{j,i}}\mm{d}t=0.
\end{align}

Exploiting \eqref{202103201834} and \eqref{202103201834x}, there exists a sequence $\{S^1_{1/m}(f^{n_m}(t))\}_{m=1}^\infty$ of simple functions such that
\begin{align}\nonumber
 \int_0^T\|S^1_{1/m}(f^{n_m}(t)) -f(t)\|_{H^{j,i}} \leqslant 1/m,
\end{align}
which yields that
\begin{align} \nonumber
\lim_{m\to \infty}\int_0^T\| S^1_{1/m}(f^{n_m}(t))-f(t)\|_{H^{j,i}}\mm{d}t=0.
\end{align}
Thus
\begin{align}\nonumber
 \| S^1_{1/m}(f^{n_m}(t))-f(t)\|_{H^{j,i}}\to 0\mbox{ for a.a. }t\in I_T \mbox{ (a subsequence)}.
\end{align}
Hence $f(t)$: $I_T\to H^{j,i}$ is also strongly measurable. Thanks to this fact and the regularity of $f$, i.e.,
$$
 \infty>\begin{cases}
\int_0^T\|  f (t) \|_{H^{j,i}}^p\mm{d}t &\hbox{ for }p\in [1,\infty); \\
 \mm{ess}\sup_{t\in I_T}\|f\|_{H^{j,i}}  & \hbox{ for }p=\infty,
      \end{cases}$$
we immediately get $f\in L^p_TH^{j,i}$.

(2) The conclusion is obvious by Pettis theorem (see Theorem 7 in APPENDIX E in \cite{ELGP}) and the separability of $X$. In addition, \eqref{2022104121338} is obvious due to the triangle inequality of norm.

(3)
By Lemma \ref{20201003302206}, there exists a function  $\tilde{f}\in H^i(\mathbb{T}\times \mathbb{R})$ such that
$$\tilde{f}|_{\Omega}=f,\ \| \tilde{f}\|_{H^i(\mathbb{T}\times \mathbb{R})}\lesssim  \| {f}\|_i  .$$
Let $\varepsilon\in (0,1)$. Then we denote  $S_\varepsilon(\tilde{f})$  the mollifications of  $\tilde{f}$ with respect to the 2D  variable $(y_1,y_2)$. Obviously, $S_\varepsilon(\tilde{f})\in L^p_TH^j$
\begin{align}
&S_\varepsilon(\tilde{f}) \to {f} \mbox{ strongly in }L^p_TH^i\mbox{ for }p>1,\label{20221033021274} \\
& \|S_\varepsilon(\tilde{f} )\|_i\lesssim  \| f \|_j\mbox{ for any } t\in I_T,  \nonumber
\end{align}
which implies that
\begin{align}\nonumber
 \|S_\varepsilon(\tilde{f} )\|_{L^p_TH^j}\lesssim  \| g\|_{L^p_TH^j}.
\end{align}
Thus there exists $\chi\in {L^p_TH^j}$ such that
\begin{align}
 S_\varepsilon(\tilde{f} )
                               \rightharpoonup  f\mbox{ weakly in }L^p_TH^j   \hbox{ for }p>1  ,
\end{align}
which, together with \eqref{20221033021274}, yields \eqref{202104032134} for $p>1$.

Thanks to the above result for $p>1$, we easily see that \eqref{202104032134} also holds for $p=\infty$.

(4)--(a) Since $\{D^h_1f\}_{h\in (0,1)}$ is uniformly bounded in $L^\infty_TH^{j,i}$, then, for any multindex $\alpha$ satisfying $|\alpha|=i$,
\begin{align}\nonumber
D^\tau_1\partial^\alpha f \rightharpoonup  \omega    \mbox{ weakly-* in }L^\infty_TL^2 \mbox{ (a subsequence)}.
\end{align}
This mean that,  for any $\chi\in H^1 $ and for any $\phi\in C_0^\infty(I_T)$,
\begin{align}
&-\int_0^T\int  \partial^\alpha f \partial_1\chi \mm{d}y \phi\mm{d}t=- \lim_{\tau\to 0}\int_0^T\int   \partial^\alpha f D^\tau_1\chi  \mm{d}y\phi\mm{d}t \nonumber \\
&=\lim_{\tau\to 0}\int_0^T\int D^\tau_1 \partial^\alpha f \chi \mm{d}y \phi\mm{d}t= \int_0^T\int \omega \chi  \mm{d}y\phi\mm{d}t \mbox{ (a subsequence)}.
\label{202103200safas1saf925}
\end{align}
Since $H^1 $ is a separable space,  thus we further derive from \eqref{202103200safas1saf925}
 that $ \omega= \partial_1\partial^\alpha f=\partial^\alpha \partial_1f $.
Thus, by the first assertion in Lemma \ref{20021032019018}, we get $f\in  L^\infty_TH^{1,i}$.

(4)--(b) If additionally  $\{D^{-\tau}_1D^\tau_1f\}_{|\tau|\in (0,1)}$ further  is uniformly bounded in $L^2_TH^{i}$ with respect to  $\tau\in (0,1)$, then,  for any  multiindex $\alpha$ satisfying $|\alpha|=i$,
\begin{align}
D^{-\tau}_1D^\tau_1 \partial^\alpha f \to \psi^\alpha    \mbox{ in }L^2_TL^2\mbox{ and }
D^\tau_1\nabla^i f \to \nabla^i \partial_1 f    \mbox{ in }L^2_TL^2 \mbox{ (a subsequence)}.
\label{2022103202014}
\end{align}

By \eqref{2022103202014},  for any multiindex $\alpha$ satisfying $|\alpha|=i$,  for any $ {\chi}\in C^2(\Omega) $, and for any $\phi\in C_0^\infty(I_T)$,
 \begin{align}
 &-\int_0^T \int  \partial^\alpha \partial_1 f \partial_1 {\chi}\mm{d}y\phi\mm{d}t = -\lim_{h\to 0} \int_0^T \int D^h_1\partial^\alpha  f  D^{h}_1 {\chi}\mm{d}y\phi\mm{d}t  \nonumber  \\
&= \lim_{h\to 0} \int_0^T \int D^{-h}_1D^h_1 \partial^\alpha f  {\chi} \mm{d}y\phi\mm{d}t
 = \int_0^T \int  \psi^\alpha  {\chi}\mm{d}y\phi\mm{d}t \mbox{ (a subsequence)}.\nonumber
\end{align}

By a density argument, we further derive from the above identity that, for any $\chi\in H^1$ and for any $\phi\in C_0^\infty(I_T)$,
 \begin{align}
-\int_0^T \int  \partial^\alpha \partial_1 f \partial_1 \chi\mm{d}y\phi\mm{d}t =
\int_0^T \int  \psi^\alpha \chi\mm{d}y\phi\mm{d}t ,\nonumber
\end{align}
which implies that $\partial_1^2\nabla^i f=\nabla^i \partial_1^2f\in L^2_TL^2$. Moreover $f\in L^2_TH^{2,i}$.

(4)--(c) If additionally $D_1^\tau\partial^\alpha f^\varepsilon$ in $C^0(\overline{I_T},L^2_{\mm{weak}})$ and  $D_1^\tau\partial^\alpha f^\varepsilon$ is uniformly continuous in $H^{-1}$,  then we have
\begin{align}
 D_1^\tau\partial^\alpha f^\varepsilon \to \vartheta \mbox{ in } C^0(\overline{I_T},L^2_{\mm{weak}})\mbox{ (a subsequence)}, \nonumber
\end{align}
which implies that, for any $\chi\in L^2$ and for any $\phi\in C_0^\infty(I_T)$,
\begin{align}
\lim_{ \tau\to 0} \int_0^T\int D_1^\tau\partial^\alpha f^\varepsilon \chi \mm{d}y\phi\mm{d}t= \int_0^T\int \vartheta \chi   \mm{d}y\phi\mm{d}t \mbox{ (a subsequence)}
\label{202210320asfdsa2014}
\end{align}
Thus we get $ \vartheta= \partial^\alpha \partial_1 f$ from \eqref{202103200safas1saf925} and \eqref{202210320asfdsa2014}.

(5)--(a)  Let us first consider the case $f\in C^0(\overline{I_T},H^i)$. Let $s\in \overline{I}_T$ be any given. Thanks to \eqref{20221033311357}, it is easy to check that for any $\varepsilon>0$, there exists a $\delta$ such that for any $t\in \overline{I}_T$ satisfying $|t-s|<\delta $,
$$
\|F(y,t)-F(y,s)\|_i =\|f(\varphi(y,t),t)-f(\varphi(y,s),s)\|_i  \leqslant \varepsilon.
$$
Thus
\begin{align}
\label{2022104031855}
F:=f(\varphi ,t)\in C^0(\overline{I_T},H^i).
\end{align}

Now we further consider the case $f\in L^p_TH^i$.  Let $\tilde{f}=f$ in $\Omega_T $ and $\tilde{f}=0$ outside $\Omega\times (\mathbb{R}\backslash I_T)$.
Let $S^t_\varepsilon(\tilde{f})$ be defined as in \eqref{202210205162222}.
  We can check that $S^{t}_\varepsilon (\tilde{f} ) \in  C^0 (\mathbb{R},H^i )  $ and
\begin{align}
&  S^{t}_\varepsilon (\tilde{f}(x,t))|_{x=\varphi}\to F( y,s) \mbox{ strongly in }{L^p_TH^i}\mbox{ for }p>1. \label{20201042512113}
\end{align}
In addition, it is easy to verify that, for a.e. $t\in I_T$,
\begin{align}
&\|F\|_{ i}\lesssim P(\|\varsigma\|_{ 3}) \|f \|_{i}.
\label{2020104saf2512113}
\end{align}
Thus we easily see from \eqref{2022104031855}, \eqref{20201042512113} and \eqref{2020104saf2512113} that \eqref{2020103250855} and
 \eqref{2021sfa04031901} hold.
Next we turn to deriving \eqref{202104031901} and \eqref{2022104101908}.

Thanks to the regularity $ \varsigma \in C^0(\overline{I_T}, H^3) $, we have (after possibly being redefined on a set of measure zero)
 \begin{align}
& \tilde{\varphi} (\tilde{y}):\overline{ \Omega_T} \to \tilde{\varphi}( \overline{\Omega_T})\mbox{ is a  homeomorphism  mapping},\label{201808151safda856sdfaxxxxx}
\end{align}
where $\tilde{\varphi} (y,t):=(\varphi (y,t),t)$, please refer to  (8.12) and (8.13) in \cite{JFJSZWC}. In particular, for given $t$,
\begin{align}
&  \varphi^{-1}(y,t)\in C^0(\overline{\Omega})\mbox{ and  }\varphi^{-1}(y,t): \Omega   \to \Omega  \mbox{ is a  homeomorphism  mapping} \label{20sdfa2104031847}
 \end{align}
where $\varphi^{-1}(y,t)$ denotes the inverse mapping of $\varphi$ with respect to $y $.

It is easy to check that
\begin{align}
&\nabla (\varphi^{-1})=  (\nabla \varphi)^{-\mm{T}}|_{y=\varphi^{-1} } .\label{20220104031929x}
\end{align}
Thanks to Lemma \ref{20221033311357}, \eqref{201808151safda856sdfaxxxxx}, \eqref{20sdfa2104031847} and \eqref{20220104031929x}, similarly to \eqref{2022104031855},   we can verify that
\begin{align}
\label{2022104131754}
\varphi^{-1}\in C^0(\overline{I_T},H^3 )
\end{align}
Thus \eqref{202104031901} obviously holds by following the argument of \eqref{2022104031855} again.

(5)--(b)
Let $\phi\in C_0^\infty(I_T)$ and  $\psi\in C^\infty_0(\Omega)$.
Let   $S_\delta^t$, resp. $S_\varepsilon$ denote the 1D, resp. 2D mollifiers with respect to variables $t$, resp. $(y_1,y_2)$. Let $S_\nu^t$ is defined as $S_\delta^t$ with $\nu$ in place of $\delta$. Then we can compute out that, for sufficiently small $\delta$,  $\varepsilon$,    and $\nu$,
$$
\begin{aligned}
& -\int_0^T\int S_\delta^t (S_\delta( f( x,t)))|_{x=y+ S_\nu^t( \varsigma)  }) \psi\phi_t\mm{d}y\mm{d}t\\
&=\int_0^T\int ( S_\delta^t (S_\varepsilon( f_t(x,t)))|_{x=y+ S_\nu^t( \varsigma)  }+S_\nu^t(\varsigma_t)\cdot S_\delta^t (S_\varepsilon(\nabla f(x,t)))|_{x=y+ S_\nu^t( \varsigma)  }) \psi\phi\mm{d}y\mm{d}t,
\end{aligned}
$$
which implies that
$$
\begin{aligned}
 -\int_0^T\int  {F}  \psi\phi_t\mm{d}y\mm{d}t =\int_0^T\int   (f_t +\varsigma_t\cdot \nabla f)|_{x=\varphi} \psi\phi\mm{d}y\mm{d}t.
\end{aligned}
$$

Noting that $C^\infty_0(\Omega)$ is density in $L^2$, 
thus we have, for any $\varphi \in L^2$ and for any $\phi\in C_0^\infty(I_T)$,
$$
\begin{aligned}
 -\int_0^T\int F\varphi\phi_t\mm{d}y\mm{d}t=\int_0^T\int (f_t+\varsigma_t\cdot\nabla f)|_{x=\varphi}\varphi\phi\mm{d}y\mm{d}t,
\end{aligned}
$$
which immediately implies that $F_t=(f_t(x,t) +\varsigma_t\cdot \nabla f(x,t))|_{x=\varphi}\in  L^p_TL^2$.  Exploiting the third assertion in Lemma \ref{20021032019018} and \eqref{2021sfa04031901}, we further have $F_t \in  L^p_TH^i$.

Thanks to the regularity $(\varsigma,\varsigma_t)$, we have (after possibly being redefined on a set of measure zero)
 \begin{align}
& \tilde{\varphi} (\tilde{y}):\overline{ \Omega_T} \to \tilde{\varphi}( \overline{\Omega_T})\mbox{ is a  homeomorphism  mapping},\nonumber \\
& \tilde{\varphi} (\tilde{y}): \Omega_T \to \tilde{\varphi} ( {\Omega_T})\mbox{ is a }C^1\mbox{-diffeomorphic mapping},\nonumber
\end{align}
where $\tilde{\varphi} (y,t):=(\varphi (y,t),t)$, please refer to  (8.12) and (8.13) in \cite{JFJSZWC}. Moreover,
\begin{align}
&\tilde{\varphi}^{-1}(\tilde{x})
=(\varphi^{-1}(y,t),t ),\ \nabla_{\tilde{x}}\tilde{\varphi}^{-1}
=(\nabla_{\tilde{y}}\tilde{\varphi})^{-1}|_{\tilde{y}=\tilde{\varphi}^{-1}},\nonumber
 \end{align}
where  $\tilde{x}=(x,t)$.
In particular, we compute out that
\begin{align}
&\partial_t\varphi^{-1} =-((\nabla  \varphi )^{-1} \varsigma_t )|_{y=\varphi^{-1} }, \label{20220104031929}
\end{align}
Thus we immediately get \eqref{2021040sdaf32132} by  \eqref{2022104131754} and  \eqref{20220104031929}. This completes the proof.\hfill $\Box$
\end{pf}

\vspace{4mm} \noindent\textbf{Acknowledgements.}
The research of Fei Jiang was supported by NSFC (Grant Nos. 12022102) and the Natural Science Foundation of Fujian Province of China (2020J02013), and the research of Song Jiang by National Key R\&D Program (2020YFA0712200), National Key Project (GJXM92579), and
NSFC (Grant No. 11631008), the Sino-German Science Center (Grant No. GZ 1465) and the ISF-NSFC joint research program (Grant No. 11761141008).

\renewcommand\refname{References}
\renewenvironment{thebibliography}[1]{%
\section*{\refname}
\list{{\arabic{enumi}}}{\def\makelabel##1{\hss{##1}}\topsep=0mm
\parsep=0mm
\partopsep=0mm\itemsep=0mm
\labelsep=1ex\itemindent=0mm
\settowidth\labelwidth{\small[#1]}
\leftmargin\labelwidth \advance\leftmargin\labelsep
\advance\leftmargin -\itemindent
\usecounter{enumi}}\small
\def\newblock{\ }
\sloppy\clubpenalty4000\widowpenalty4000
\sfcode`\.=1000\relax}{\endlist}
\bibliographystyle{model1b-num-names}

\end{document}